\documentclass[11pt]{amsart}
\usepackage{amsmath,amssymb,amsthm}
\usepackage{amsmath,amsthm,amscd,amssymb,amsfonts, amsbsy, enumerate}
\usepackage{eufrak}
\usepackage{bm}
\usepackage{xcolor}
\SetMathAlphabet{\mathfrak}{bold}{U}{euf}{b}{n}
\usepackage{latexsym}
\usepackage{txfonts}
\usepackage[colorlinks,citecolor=red,pagebackref,hypertexnames=false]{hyperref}

\numberwithin{equation}{section}
\newtheorem{theorem}{Theorem}[section]
\newtheorem{lemma}[theorem]{Lemma}
\newtheorem{proposition}[theorem]{Proposition}
\newtheorem{corollary}[theorem]{Corollary}
\theoremstyle{definition}
\newtheorem{definition}[theorem]{Definition}
\newtheorem{remark}[theorem]{Remark}
\newcommand{\T}{\mathbb T}
\newcommand{\R}{\mathbb R}
\newcommand{\E}{\mathbb E}
\newcommand{\B}{\mathcal B}
\newcommand{\dd}{\,\mathrm d}
\newcommand{\supp}{\operatorname{supp}}
\newcommand{\dist}{\operatorname{dist}}
\newcommand{\Lip}{\operatorname{Lip}}
\newcommand{\one}{\mathbf 1}

\newcommand{\la}{\lambda}

\title[Richardson dispersion and anomalous dissipation]{Richardson dispersion and anomalous dissipation  for the random velocity model with controlled sweeping}
\author{Charles Fefferman}
\address{Department of Mathematics, Princeton University, Fine Hall, Washington Road, Princeton, NJ 08544, USA}
\email{cf@math.princeton.edu}

\author{Svitlana Mayboroda}
\address{Department of Mathematics, ETH Z\"urich, R\"amistrasse 101, 8092 Z\"urich, Switzerland}
\email{svitlana.mayboroda@math.ethz.ch}

\subjclass[2020]{Primary 76F25; Secondary 76F05, 60G60, 34F05}
\keywords{Richardson dispersion, pair dispersion, turbulence, random velocity fields, sweeping, Kubo number}
\date{\today}
\begin{document}
\maketitle

\begin{abstract} The present paper establishes the Richardson dispersion law, predicting the particle separation with time according to $\E R_t^2 \sim t^3$, as well as anomalous dissipation for the associated tracer, for a random H\"older continuous velocity vectorfield with controlled sweeping.  \end{abstract}

\section{Introduction}\label{sec:intro}
\subsection{History and motivation}

The present paper establishes the Richardson dispersion law and the anomalous dissipation of a tracer for a random H\"older continuous velocity vectorfield. The construction is homogeneous in time and geometrically dense in scales. In addition, it exhibits sweeping control long conjectured to be the key to important aspects of the turbulent behavior: small eddies are materially transported by the larger ones. Its physics companion is simultaneously released as \cite{JMDLSDDE26}. The Obukhov-Corrsin-Yaglom law for the same vectorfield is established in the simultaneously released \cite{DL26}.

Richardson's 1926 prediction for turbulent relative dispersion is one of the oldest predictions in classical turbulence \cite{Ri26}. Yet it has evaded rigorous mathematical treatment to this date. Richardson proposed that the effective diffusivity of a particle pair at separation $r$ should grow like $r^{4/3}$, or, equivalently, that the square of the separation between the particles grows as $t^3$ with time:
\[
        \E R_t^2 \sim t^3.
\]
This scaling is consistent with Kolmogorov’s later prediction of smoothness of a turbulent vectorfield being roughly 1/3 (also open to this date and actually requiring an adjustment due to intermittency effects). More generally, if velocity increments at scale $r$ are of order $r^\alpha$, dimensional reasoning suggests the separation scale $r(t)\sim t^{1/(1-\alpha)}$ and hence
\[
        \E R_t^2\sim t^{2/(1-\alpha)}.
\]
We refer to \cite{Fr95, Saw01,FGV01} for classical reviews.

The intuitive reason behind the scaling is as follows. If velocity increments across a distance $r$ are of size $r^\alpha$, a pair at separation $r$ drifts apart at speed $r^\alpha$, and the ordinary differential equation $\dot r=r^\alpha$ gives $r^{1-\alpha}\approx t$, that is, $r^2\approx t^{2/(1-\alpha)}$; for the Kolmogorov exponent $\alpha=1/3$ this is $t^3$. The heuristic is not a proof, since the increments have no preferred sign: they average to zero, and only their fluctuations, accumulated over many independent eddies, push the particles apart. 

Rigorous results targeting a quantitative Richardson law are scarce and quite far from expected physical reality. The classical example which allows for explicit computations is the Kraichnan model \cite{Kr68}. Here, however, the advecting velocity is Gaussian and white in time. This is much rougher than expected regularity in time for real turbulence,
%, and moreover, the choice of the spatial covariance that yields the Richardson law is not the Kolmogorov spatial scaling choice; 
see, for instance, \cite{FGV01,LR02}. We also mention the work of Fannjiang \cite{Fan03}, which proves generalized Richardson laws for colored-in-time synthetic velocities in an inertial-scale limit. In the regime yielding a closed pair-dispersion law, however, the limiting dynamics are again of white-in-time Kraichnan type, and the finite-correlation Kolmogorov endpoint is not itself covered by the limit theorem. 

In a different direction, much more recently, Hess-Childs and Rowan \cite{HCR25b} constructed an explicit deterministic divergence-free velocity $V\in L^\infty_t C^\alpha_x(\mathbb T^2)$ for which two independently diffusive tracers satisfy the sharp two-sided estimate
\[
        \E R_t^2\sim R_0^2+\mu t+t^{2/(1-\alpha)},
\]
where $\mu\to 0$ is the molecular diffusivity constant. In particular, at $\alpha=1/3$ they obtain Richardson's $t^3$ law, uniformly as $\mu$ tends to zero. This is the only genuinely colored-in-time construction yielding Richardson, however, as the authors themselves emphasize, their velocity is a deterministic, explicitly engineered two-dimensional flow, rather than a stationary random turbulence model with finite correlation time. In particular, is not time-homogeneous: the Richardson estimate is proved for trajectories launched from a distinguished initial time, rather than uniformly with respect to the time origin.  Without going into technical details, we will just say that the
complexity here is concentrated near $t=0$ and thins out afterwards,
rather than being homogeneous in time as physical expectations would
require.

More generally, the present paper builds on  a recent line of rigorous work on anomalous dissipation of passive scalars and on
the related breakdown of uniqueness for rough transport.
% \cite{DEIJ22,CCS23,BD23,ACM19,ACM19b,BCCDS23,AV25,BSW23,JS26,HCR25,R24,BBDM26,CHR26,HCR25b}. We will discuss them below. 
%All these velocity fields are either deterministic and assembled scale by scale, or random but white in time or autonomous; the refreshing ensembles of \cite{CHR26}, the closest in spirit to ours, are studied there for uniqueness and for upper bounds on separation, and, as those authors point out, random fields without self-advection should not be expected to behave like fluids. 
The latter, as we mentioned above, has a rich and growing history that our work builds on. Drivas, Elgindi, Iyer and
Jeong \cite{DEIJ22} opened the subject with the first velocity field,
in $L^1_tC^{1-}_x$, that provably dissipates a passive scalar, built
from alternating shears refined in time in the spirit of Pierrehumbert;
Colombo, Crippa and Sorella \cite{CCS23} reached the full supercritical
Obukhov--Corrsin regime and showed that vanishing
diffusivity selects no unique limit. Bru\`e and De~Lellis \cite{BD23}, using the mixing flows of Alberti,
Crippa and Mazzucato \cite{ACM19,ACM19b}, obtained anomalous dissipation
for the forced three-dimensional Navier--Stokes equations. Together
with Colombo, Crippa and Sorella \cite{BCCDS23}, they reached
$L^3_tC^\alpha_x$ for every $\alpha<1/3$; both constructions use
viscosity-dependent forcing.  Armstrong and Vicol \cite{AV25}
introduced fractal homogenization and obtained anomalous diffusion
for every $H^1$ datum, with no distinguished time and with the small
scales advected by the large ones, and Burczak, Sz\'ekelyhidi and Wu
\cite{BSW23} obtained the zeroth law of scalar turbulence for a scalar
carried by a $C^{1/3-}$ weak solution of the Euler equations.
Johansson and Sorella \cite{JS26} achieved anomalous dissipation with
an autonomous two-dimensional $C^\alpha$ field, $\alpha<1$, and
Hess-Childs and Rowan \cite{HCR25} built a universal total
dissipator, a single $L^\infty_tC^\alpha_x$ field under which every
datum loses all of its variance. On the random side, Rowan \cite{R24}
proved anomalous diffusion for the Kraichnan ensemble and constructed correlated-in-time variants with no anomalous dissipation at any fixed positive correlation time; Bagnara, Boutros, De~Lellis and the second author
\cite{BBDM26} established regularity thresholds above which anomalous
dissipation, Richardson diffusion and the Obukhov--Corrsin law are
impossible for certain classes of autonomous random fields; and Colombo, Hess-Childs and Rowan \cite{CHR26} established uniqueness above Hölder exponent \(1/2\), away from zeros, for autonomous random multiscale fields, and constructed nonuniqueness examples below this threshold in \(d\ge4\). For refreshing ensembles they obtained corresponding thresholds and quantitative separation estimates, including lower bounds for the spread of trajectories in their nonuniqueness constructions.  %All these velocity fields are either deterministic and assembled scale by scale, or random but white in time or autonomous; the refreshing ensembles of \cite{CHR26}, the closest in spirit to ours, are studied there for uniqueness and for upper bounds on separation, and, as those authors point out, random fields without self-advection should not be expected to behave like fluids. 
With the exception of \cite{HCR25b} discussed above, none of these yields a Richardson law.

Anomalous dissipation is not unrelated to Richardson dispersion, but
it is far from implying it, implicitly or explicitly. By the Lagrangian fluctuation--dissipation relation of Drivas and Eyink \cite{DE17}, anomalous dissipation for passive scalars is intimately tied to spontaneous stochasticity. Informally, spontaneous stochasticity asserts that initially nearby or coincident trajectories retain a nontrivial separation in the zero-noise limit, without prescribing a quantitative growth law. Richardson dispersion is substantially stronger: it specifies the exact scaling power \(\E R_t^2\sim t^3\) over the inertial range. And indeed, the numerical evidence and scaling arguments show that in synthetic example the dissipation exponent can be different from 3/2.

This is connected to an important issue emphasized in the physics literature, the role of \emph{sweeping} in synthetic turbulence models. In standard kinematic simulations the Fourier modes representing small eddies evolve without being materially transported by the larger scales. Thomson and Devenish \cite{TD05} argued that this absence of sweeping gives small-scale velocity increments along particle trajectories an incorrect decorrelation time and can change the predicted pair-dispersion exponent; see also the subsequent discussion in \cite{OVSH06,DT09} and the analytical treatment of sweeping decorrelation by Eyink and Benveniste \cite{EB13}. This issue is specific to synthetic Eulerian velocity models: in an actual flow, a small eddy and two nearby tracers are to leading order transported together by the larger-scale velocity. To the best of our knowledge, aside from the above, the only mathematical result taking sweeping into account is the aforementioned new model of Armstrong and Vicol \cite{AV25}, cf. also \cite{AV25survey,BSW23} which brought powerful homogenization techniques to the table. Their velocity however is still a carefully engineered, deterministic, space-time periodic fractal field built from alternating shear flows, rather than a random stationary turbulence model, and their scale hierarchy is extremely lacunary: the active length scales are separated double-exponentially, rather than forming anything like a dense geometric inertial range. Their theorem concerns anomalous scalar diffusion only, not the Richardson pair-dispersion power law. In a different but related multiscale homogenization problem, the authors of \cite{ABK26, ABK-superdiffusion} remove scale separation for a stationary random incompressible drift and prove one-particle superdiffusion, $\mathbb{E}|X_t|^2 \sim t(\log t)^{1/2}$ in the critical case and $t^{2/(2-\gamma)}$ for self-similar drifts with small Hurst exponent $\gamma$; these are not Richardson pair-dispersion theorems either.

In the present paper and its companion in physics, \cite{JMDLSDDE26}, we introduce a random velocity with prescribed spatial Hölder regularity and scale-dependent renewal times. The same hierarchy of scales is renewed throughout the time interval, with smaller eddies transported by larger ones. We prove Richardson dispersion and anomalous dissipation for this field.  The tracer particles  satisfy deterministic ordinary differential equations in a common realization of the velocity, all randomness comes from the velocity ensemble. We prove in Theorem~\ref{thm:main}, uniformly in the ultraviolet cutoff and for every initial separation $r_0\ge\ell_J$, that 
\[
c\bigl(r_0^{\,1-\alpha}+\kappa t\bigr)^{\frac{2}{1-\alpha}}
\le \E R_t^2
\le C\bigl(r_0^{\,1-\alpha}+t\bigr)^{\frac{2}{1-\alpha}}.
\]
Here, $\kappa$ is the Kubo number, a dimensionless constant measuring the ratio between the correlation time of an eddy and its characteristic turnover time. In our model, an eddy of size $\ell$ has lifetime $\tau_\ell\approx\kappa\,\ell^{1-\alpha}$, so $\kappa$ controls how long structures persist relative to the natural scale-dependent advective time. We need $\kappa$ to be small for the proof -- perhaps, a limitation of the argument rather than a necessity as the numerical results in \cite{JMDLSDDE26} show, and the left and right bound differ with $\kappa$. However, it does not have to vanish and depends on the parameters of the system only. 

In particular, when $\alpha=1/3$ and the Kubo number is fixed, in the admissible time range,
\[
        c\kappa^3 t^3\le \E R_t^2\le Ct^3.
\]
Thus the model exhibits the Richardson cubic law with constants independent of the small-scale cutoff, for a random velocity with finite, scale-dependent correlation times and with the smaller eddies transported by the larger ones. Furthermore, we establish anomalous dissipation.

With this, let us pass to the model.

\subsection{The model}\label{sec:model}

Fix an integer $d\ge2$ and $0<\alpha<1$, and put $\beta=1-\alpha$.
The physical space is $\T^d=\R^d/\mathbb Z^d$, equipped with its
flat distance. Let $R>0$ and let $\psi\in C_c^4(\R^d)$ be radial,
not identically zero, and supported in $\overline{B(0,R)}$.
Let $e_1,\ldots,e_d$ be the coordinate vectors and set
\[
 \mathcal A=\{(p,q):1\le p<q\le d\},\qquad
 J^{pq}=e_q e_p^{\mathsf T}-e_p e_q^{\mathsf T},\qquad
 W^{pq}=J^{pq}\nabla\psi.
\]
Every $W^\nu$, $\nu\in\mathcal A$, belongs to
$C_c^3(\R^d;\R^d)$ and is divergence-free, because $J^\nu$ is
skew-symmetric. In dimension two there is a single mode,
$W^{12}=(-\partial_2\psi,\partial_1\psi)=\nabla^\perp\psi$. 

The averaged velocity increments of one scale of the model are
governed by the structure function of this family. For $z\in\R^d$ define
\begin{equation}\label{eq:structure}
 D(z)=\sum_{\nu\in\mathcal A}\int_{\R^d}
 [W^\nu(y)-W^\nu(y-z)]\otimes[W^\nu(y)-W^\nu(y-z)]\dd y.
\end{equation}
Lemma~\ref{lem:structure} shows that, for $z\ne0$,
\begin{equation}\label{eq:isotropy}
 D(z)=D_L(|z|)\widehat z\otimes\widehat z
       +D_N(|z|)(I-\widehat z\otimes\widehat z),\qquad
 \widehat z=z/|z|,
\end{equation}
where $D_N$ is the eigenvalue in each of the $d-1$ transverse
directions, and that $D_N=D_L+rD_L'/(d-1)$, the K\'arm\'an--Howarth
identity \cite{KH38}, which is the divergence-free condition written
for an isotropic tensor. We are assuming the following structural hypothesis on the wavelets
\begin{equation}\tag{ND}\label{eq:ND}
 D_L\text{ is nondecreasing on }(0,\infty),
\end{equation}
or, equivalently $D_N\ge D_L$. Proposition~\ref{prop:explicitwavelet}
verifies \eqref{eq:ND} for $\psi(x)=(1-|x|^2)_+^5$ in every dimension
$d\ge2$. 

We defer the detailed discussion of the condition \eqref{eq:ND} to the next section, but here let us just mention that it ensures the correct sign of the renewal gains -- the main contribution to our observable, and that it is rather generic for standard families of wavelets.

Recall that $d\ge2$, $0<\alpha<1$, and $\beta=1-\alpha$.
Let $\kappa>0$ be the Kubo number and define
\begin{equation}\label{eq:scales}
 \ell_j=2^{-j},\qquad A_j=\ell_j^\alpha,
 \qquad \la_j=\ell_j^{-\beta},\qquad
 \tau_j=\kappa\,2^{-\lceil\beta j\rceil}.
\end{equation}
Set
\[
 j_0=\min\{j\ge0:R\ell_j\le1/4\},\qquad J\ge j_0,
 \qquad I_j=\{0,\ldots,2^j-1\}^d.
\]
Thus $A_j$ is the amplitude of scale $j$, $\la_j=A_j/\ell_j$ the size of its
gradient, and $\tau_j$ its lifetime. The ratio of lifetime to
turnover time lies in $(\kappa/2,\kappa]$; the dyadic rounding
aligns the renewal times across scales (Lemma~\ref{lem:nested}). The integer $J$ is finite but otherwise arbitrary; all
constants below are independent of $J$. Put
$y_i^j=\ell_j i\pmod{\mathbb Z^d}$ and $Q_j=[0,\ell_j)^d$.
The periodized wavelets are
\[
 W_j^\nu(x)=\sum_{n\in\mathbb Z^d}W^\nu((x+n)/\ell_j),\qquad x\in\T^d.
\]
Let $K_j=\lceil\tau_j^{-1}\rceil$ and
$I_{jk}=[k\tau_j,(k+1)\tau_j)\cap[0,1)$ for
$0\le k<K_j$. For each such $(j,k)$, let $\xi_k^j$ be a random
variable uniformly distributed on $Q_j$. For $i\in I_j$ and
$\nu\in\mathcal A$, let $\eta_{i,k}^{j,\nu}$ be real random variables
satisfying
\begin{equation}\label{eq:coefficients}
 \E\eta_{i,k}^{j,\nu}=0,\qquad
 \E(\eta_{i,k}^{j,\nu})^2=\sigma^2>0,\qquad
 |\eta_{i,k}^{j,\nu}|\le M .
\end{equation}
The last condition means that the law of each coefficient is
supported in $[-M,M]$; we regard each coefficient as a coordinate of
the data space with values in $[-M,M]$, so that the bound holds for
every realization.
All the variables $\xi_k^j$ and $\eta_{i,k}^{j,\nu}$ are independent.
Independent signs $\pm1$ with probability $1/2$ are admissible with
$M=\sigma=1$. We realize the data for all scales $j\ge j_0$
on one countable product probability space; $\E$ and $\mathbb P$
denote expectation and probability on this space. Clearly, each finite
truncation uses only finitely many of these variables.
Thus the data space is the product of the cubes $Q_j$, one for each
window, and of the intervals $[-M,M]$, one for each coefficient.

The fields are defined inductively in $j$. Write
$u_{<j}=\sum_{m=j_0}^{j-1}u_m$, with $u_{<j_0}=0$.
On $I_{jk}$, the centers solve
\begin{equation}\label{eq:centers}
 \dot P_{i,k}^j(t)=u_{<j}(P_{i,k}^j(t),t),\qquad
 P_{i,k}^j(k\tau_j)=y_i^j+\xi_k^j.
\end{equation}
That is, every $k\tau_j$ the centers are placed roughly on a lattice (with some amount of randomness, each uniformly in its own cube) and then they are carried forward by the flow of larger eddies until the next reset after the time $\tau_j$.

Define
\begin{equation}\label{eq:velocity}
 u_j(x,t)=A_j\sum_{i\in I_j}\sum_{\nu\in\mathcal A}
        \eta_{i,k}^{j,\nu}W_j^\nu(x-P_{i,k}^j(t))\quad(t\in I_{jk}),
 \qquad u=u^{(J)}:=\sum_{j=j_0}^J u_j.
\end{equation}
Values at $t=1$ are assigned by left limits; they do not affect the
trajectories. Centers of generation $j$ are reset at its
renewal times $k\tau_j$. Every renewal of a coarser scale is also
a scheduled renewal of all finer scales (Lemma~\ref{lem:nested}).

Here and below $\mathbb N=\{0,1,2,\ldots\}$. By Lemma~\ref{lem:nested},
the union of all renewal times
for the truncation at $J$ is $\tau_J\mathbb N\cap[0,1)$.
Each $u_j$ depends only on the data of scales $j_0,\ldots,j$ and
does not depend on the truncation index $J$. Thus all the fields
$u^{(J)}$ above are consistently defined, using the same data at their shared scales.

Each $u^{(J)}$ is the truncation at scale $J$ of the random field
\begin{equation}\label{eq:ulimit}
 u^{(\infty)}:=\sum_{j\ge j_0}u_j ,
\end{equation}
defined on the same probability space: $u_j$
does not depend on $J$, and $|u_j|\le C_*A_j$ is summable
(Lemma~\ref{lem:regularity}), so
the series converges uniformly in $(x,t)$, and $u^{(\infty)}$ is
bounded, Borel and divergence-free, with
$\|u^{(\infty)}(\cdot,t)\|_{C^\alpha}\le C_*$ for every $t$. The field
$u^{(\infty)}$ is the object of Theorem~\ref{thm:ADfixed} and of
Corollary~\ref{cor:deterministic} in Subsection~\ref{sec:fixedfield}.

For deterministic initial positions $x_0,\widetilde x_0\in\T^d$, let
\begin{equation}\label{eq:particles}
 \dot x_t=u(x_t,t),\qquad
 \dot{\widetilde x}_t=u(\widetilde x_t,t),\qquad
 R_t=\dist(x_t,\widetilde x_t),\qquad r_0=R_0.
\end{equation}
When $R_t\le1/4$ we write
$z_t$ for the unique shortest representative in $\R^d$ of
$x_t-\widetilde x_t$; then $|z_t|=R_t$.

\subsection{Main results}\label{sec:results}

\begin{theorem}\label{thm:main}
Fix an integer $d\ge2$ and consider the model of
Section~\ref{sec:model}. Assume \eqref{eq:ND} and
\eqref{eq:coefficients}. There are constants
$\kappa_0\in(0,e^{-1})$ and $c,C,C_H>0$, depending only on
$d,\alpha,\psi,R,M,\sigma$ ($C_H$ is the constant of
\eqref{eq:radiusdet} below), with the following property. For every
$0<\kappa\le\kappa_0$, every $J\ge j_0$, and every deterministic pair
of initial positions with
\[
 \ell_J\le r_0\le1/8,
\]
one has $\delta_\kappa:=C\kappa(1+|\log\kappa|)\le1/4$, and for every
$t\in[0,1]$ satisfying
\begin{equation}\label{eq:localtime}
 r_0^\beta+C_Ht\le(1/4)^\beta,
\end{equation}
\begin{equation}\label{eq:maincorrected}
 (1+\delta_\kappa)\E R_t^\beta
 \ge r_0^\beta-\delta_\kappa\ell_J^\beta+c\kappa t
 \ge(1-\delta_\kappa)r_0^\beta+c\kappa t.
\end{equation}
Moreover, after changing $c$ and $C$,
\begin{equation}\label{eq:momentbounds}
 c(r_0^\beta+\kappa t)^{2/\beta}
 \le\E R_t^2\le C(r_0^\beta+t)^{2/\beta}.
\end{equation}
All constants are independent of $J$, $r_0$, and $\kappa$ in the
specified ranges.
\end{theorem}

\begin{corollary}\label{cor:cubic}
Let $d\ge2$ and $\alpha=1/3$. For times satisfying
\eqref{eq:localtime} and $t\ge r_0^{2/3}/\kappa$,
\[
 c\kappa^3t^3\le\E R_t^2\le Ct^3.
\]
\end{corollary}

Thus, for each fixed sufficiently small Kubo number, the mean-square
separation has a cubic lower bound on this time range, uniformly in
the small-scale cutoff; the time range is nonempty when $r_0$ is small
compared with a power of $\kappa$. The upper bound in
\eqref{eq:momentbounds} is a consequence of the H\"older continuity of
$u$ alone. The exponent $2/\beta$ is therefore sharp.

\medskip

\subsubsection{The passive scalar}
Fix $\varepsilon_0\in(0,1]$ and set, for $J\ge j_0$,
\[
 \mu_J:=\varepsilon_0\ell_J^{1+\alpha}.
\]
Thus the P\'eclet number of the finest scale,
$A_J\ell_J/\mu_J=1/\varepsilon_0$, does not depend on $J$. For a real
$\theta_0\in L^2(\T^d)$ and a realization of the data, let
$\theta=\theta^{(J)}$ be the energy solution (Definition~\ref{def:energy})
on $[0,1]$ of
\begin{equation}\label{eq:scalar}
 \partial_t\theta+u^{(J)}\cdot\nabla\theta=\mu_J\Delta\theta,
 \qquad\theta(0)=\theta_0 .
\end{equation}
It exists and is unique, its restriction to $[0,t_0]$ is the energy
solution on $[0,t_0]$, $\int\theta(t)=\int\theta_0$, and for
$t_0\in(0,1]$ the energy identity holds (Lemma~\ref{lem:weakbasic}):
\begin{equation}\label{eq:energy}
 \|\theta^{(J)}(t_0)\|^2+\mathbf D_J(\theta_0)=\|\theta_0\|^2,\qquad
 \mathbf D_J(\theta_0):=2\mu_J\int_0^{t_0}\|\nabla\theta^{(J)}(t)\|^2\dd t .
\end{equation}
Here and below $\|\cdot\|$ is the norm of $L^2(\T^d)$ (or of
$L^2(\T^{2d})$). The quantity $\mathbf D_J(\theta_0)\in[0,\|\theta_0\|^2]$
is the variance dissipated by time $t_0$.

\subsubsection{Anomalous dissipation}\label{sec:ADresults}

The mean dissipation is computed mode by mode under a symmetry assumption
on the coefficient laws.

\begin{definition}[Hypothesis (S)]\label{def:S}
For every scale $j$ and window $k$, the coefficients
$\eta^{j,\nu}_{i,k}$, $i\in I_j$, $\nu\in\mathcal A$, are identically
distributed, and their common law is symmetric
($\eta\overset{d}{=}-\eta$).
\end{definition}

Independent signs $\pm1$ with probability $\frac12$ satisfy
\eqref{eq:coefficients} and (S).

Retain the notation $\kappa_*$ from Lemma~\ref{lem:regularity}.

\begin{theorem}[Anomalous dissipation]\label{thm:ADgeneral}
Let $d\ge2$ and assume \eqref{eq:ND}, \eqref{eq:coefficients} and (S).
There are $\kappa_1\in(0,\kappa_*]$ and $c>0$, depending
only on $d,\alpha,\psi,R,M,\sigma$, and for every $\kappa\in(0,\kappa_1]$ a number
$\varepsilon_1(\kappa)\in(0,1]$, with the following property. Consider
the model with $\kappa\in(0,\kappa_1]$, let
$\varepsilon_0\in(0,\varepsilon_1(\kappa)]$ and $t_0\in(0,1]$, and put
$c_{\rm AD}:=c\,(\kappa t_0)^{2/\beta}$. Then for every real
$\theta_0\in L^2(\T^d)$ with $\int\theta_0=0$:
\begin{enumerate}[(a)]
\item $\liminf_{J\to\infty}\E\,\mathbf D_J(\theta_0)\ge c_{\rm AD}\|\theta_0\|^2$;
\item if $\theta_0\ne0$, then
 $\mathbb P\big(\limsup_{J\to\infty}\mathbf D_J(\theta_0)>0\big)\ge c_{\rm AD}$;
\item for every $\gamma>0$ and $A\ge1$ there is
 $J_*=J_*(\gamma,A,\kappa,\varepsilon_0,t_0)$, depending also on the
 fixed model parameters, such that
 $\E\,\mathbf D_J(\theta_0)\ge\frac12c_{\rm AD}\|\theta_0\|^2$
 for all $J\ge J_*$ and every real mean-zero
 $\theta_0\in H^\gamma(\T^d)$ satisfying
 $\|\theta_0\|_{\dot H^\gamma}\le A\|\theta_0\|$.
\end{enumerate}
\end{theorem}

The distinction between (a) and (c) is natural: $L^2$ is the natural
space for scalar variance, while a bound on the Sobolev-to-$L^2$ norm
ratio gives uniformly small high-frequency tails relative to the initial
variance and hence a threshold uniform over the corresponding initial data. Here the homogeneous Sobolev seminorm is
\[
 \|\theta_0\|_{\dot H^\gamma}^2
 :=\sum_{k\in\mathbb Z^d}(2\pi|k|)^{2\gamma}
       |\widehat\theta_0(k)|^2.
\]

Recall from Section~\ref{sec:model} that $u^{(J)}$ is the
truncation at scale $J$ of the random field $u^{(\infty)}$ of
\eqref{eq:ulimit}, which is bounded, Borel, divergence-free and
$C^\alpha$ in $x$, uniformly in $t$. For $\mu>0$
and a real $\theta_0\in L^2(\T^d)$ let $\theta^\mu$ be the energy
solution on $[0,1]$ of
\begin{equation}\label{eq:scalarfixed}
 \partial_t\theta^\mu+u^{(\infty)}\cdot\nabla\theta^\mu=\mu\Delta\theta^\mu,
 \qquad\theta^\mu(0)=\theta_0,
\end{equation}
and for $t_0\in(0,1]$ put
\[
 \mathbf D^\mu_\infty(\theta_0):=2\mu\int_0^{t_0}\|\nabla\theta^\mu(t)\|^2\dd t .
\]
The following theorem concerns this single field, which depends neither
on $\mu$ nor on $\varepsilon_0$.

\begin{theorem}[Anomalous dissipation for the fixed field]\label{thm:ADfixed}
Assume the hypotheses of Theorem~\ref{thm:ADgeneral}, fix
$\kappa\in(0,\kappa_1]$ and $\varepsilon_0\in(0,\varepsilon_1(\kappa)]$,
let $t_0\in(0,1]$, and let $c_{\rm AD}$ be as in
Theorem~\ref{thm:ADgeneral}. There is $a_*\in(0,1)$, depending only on
$d,\alpha,\psi,R,M$ and $\varepsilon_0$ (explicitly, \eqref{eq:astar}),
such that for every real $\theta_0\in L^2(\T^d)$ with $\int\theta_0=0$:
\begin{enumerate}[(a)]
\item $\liminf_{\mu\downarrow0}\E\,\mathbf D^\mu_\infty(\theta_0)\ge a_*c_{\rm AD}\|\theta_0\|^2$;
\item if $\theta_0\ne0$, then
 $\mathbb P\big(\limsup_{\mu\downarrow0}\mathbf D^\mu_\infty(\theta_0)>0\big)\ge c_{\rm AD}$;
\item for every $\gamma>0$ and $A\ge1$ there is
 $\mu_*=\mu_*(\gamma,A,\kappa,\varepsilon_0,t_0)>0$, depending also
 on the fixed model parameters, such that
 $\E\,\mathbf D^\mu_\infty(\theta_0)
 \ge\frac12a_*c_{\rm AD}\|\theta_0\|^2$
 for all $\mu\in(0,\mu_*]$ and every real mean-zero
 $\theta_0\in H^\gamma(\T^d)$ satisfying
 $\|\theta_0\|_{\dot H^\gamma}\le A\|\theta_0\|$.
\end{enumerate}
\end{theorem}

As in Theorem~\ref{thm:ADgeneral}, part (c) gives a threshold
$\mu_*$ uniform over all nonzero real mean-zero Sobolev data with
$\|\theta_0\|_{\dot H^\gamma}/\|\theta_0\|\le A$, whereas (a) also
covers data belonging to no $H^\gamma$ with $\gamma>0$.
Again the factor $\frac12$ in (c) may be replaced by any
$1-\epsilon$, $0<\epsilon<1$, with a threshold depending also on
$\epsilon$.

The theorem is proved in Subsection~\ref{sec:fixedfield} by comparing the
dissipation for $u^{(\infty)}$ with the dissipation for the truncation
$u^{(J)}$ at the scale $J$ for which $\mu_J\le\mu<2^{1+\alpha}\mu_J$; the
scales finer than $J$ change the dissipation by at most a constant
factor (Lemma~\ref{lem:disscomparison}). Part (b) says that, on an event
of probability at least $c_{\rm AD}$, a single realization of
$u^{(\infty)}$ dissipates the variance of $\theta_0$ anomalously, in
the sense of \cite{DEIJ22}. For the explicit model of
Proposition~\ref{prop:explicitwavelet} this gives a deterministic
statement, Corollary~\ref{cor:deterministic} in
Subsection~\ref{sec:fixedfield}: for every $d\ge2$, $0<\alpha<1$ and
$t_0\in(0,1]$ there is a divergence-free
$u\in L^\infty([0,1];C^\alpha(\T^d))$ for which the energy solutions of
$\partial_t\theta^\mu+u\cdot\nabla\theta^\mu=\mu\Delta\theta^\mu$ with
$\theta^\mu(0)=\cos(2\pi x_1)$ satisfy
$\limsup_{\mu\downarrow0}\mu\int_0^{t_0}\|\nabla\theta^\mu\|^2\dd t>0$.

\subsection{Outline of the proof}\label{sec:strategy}

The field $u$ is reset randomly at certain times depending on the scale but otherwise every next generation is deterministically carried by the flow of previous generations. Inside a window the flow is a
deterministic ordinary differential equation with a $C^3$ right-hand
side, and the mean-zero property of the coefficients can be exploited
only at the instants at which a window is refreshed. Probability
enters the proof through finite-dimensional averages: at each
refresh time, and there, we integrate over the fresh coefficients
and lattice shift of the refreshed scale with everything else held
fixed. No martingales are used; where
one would condition on the past, we integrate over the variables of a
single window and apply Fubini's theorem
(Section~\ref{sec:renewals}). 

\subsubsection{Richardson dispersion} We start with

\medskip

\noindent \emph{The symbol of the model} (Section~\ref{sec:regularity}). This is the main contribution to the observable. 
Averaging the increment of a freshly drawn scale-$j$ field between two
points at distance $r$ produces $\sigma^2A_j^2D(z/\ell_j)$, where $D$
is the structure function \eqref{eq:structure}. Accordingly, for a
radial test function $\Phi$ of the separation we set
\[
 \Gamma_j(r):=\frac{\sigma^2A_j^2}{2}\,D(z/\ell_j):D^2\Phi(z),\qquad |z|=r,
\]
which depends on $z$ only through $r=|z|$ by the isotropy
\eqref{eq:isotropy}; $\tau_j^2\Gamma_j(r)$ is, to leading order, the mean
change of $\Phi(z)$ produced by one window of scale $j$ while the
particles are at distance $r$, and exactly the averaged jump, at a
refresh of scale $j$, of the corrected observable introduced below
(Proposition~\ref{prop:renewal}). In the body of the paper
$\Phi$ is replaced by its scale-adapted version $\Phi_j$ of
Section~\ref{sec:testfunctions}, which gives the precise definition
\eqref{eq:Gamma}. One computes $D:D^2\Phi$ explicitly; for
$\Phi=r^{\beta}$ it equals
$\beta r^{\beta-2}[(d-1)D_N-\alpha D_L]$. The concavity of $r^\beta$
makes the longitudinal term negative, the $d-1$ transverse directions
make the other term positive, and \eqref{eq:ND}, in the form
$D_N\ge D_L$, makes the sum positive. This is the place where
 \eqref{eq:ND} is heavily used. Counting
powers at $r\approx\ell_j$, one window of scale $j$ changes the
average of $R^\beta$ by
$\tau_j^2\Gamma_j(r)\approx\tau_j^2A_j^2r^{\beta-2}\approx\kappa^2\ell_j^\beta\approx\kappa\tau_j$,
that is, at rate $\kappa$ per unit time independently of the scale;
$\beta=1-\alpha$ is the exponent for which the gain per window is
scale-free.

\medskip\noindent
\emph{Test functions} (Section~\ref{sec:testfunctions}). We would like to estimate $\Phi=r^\beta$, but factually, it
is modified twice, due to the low frequency and high frequency cutoffs. Below the finest scale it is continued as a
logarithm: there the field is Lipschitz, $\log r$ is the observable
whose derivative along a Lipschitz flow is bounded, and the
continuation keeps the gain at scale $J$ bounded below uniformly in
$r$. This step is the source of the uniformity in $J$. Then, for each
scale $b$, a copy $\Phi_b$ is formed in which the scales much coarser
than the current separation are switched off, at a ratio $L$; this
keeps every sum over scales of length $O(\log L)$ at the price of an
error $O(L^{-\beta})$. The cutoff is nondecreasing in $r$, so that the
sign computation survives it.

\medskip\noindent
\emph{Commutator estimates} (Section~\ref{sec:operators}). The state
of the system, a configuration, consists of the two particles and all
wavelet centers; between refreshes it solves an autonomous ODE whose
generator $\B=\sum_a\B_a$ splits by scale. Here $\B_a$ is the
first-order differential operator
\[
 \B_aF=u_a(x)\cdot\nabla_xF+u_a(\widetilde x)\cdot\nabla_{\widetilde x}F
 +\sum_{m>a}\sum_{i\in I_m}u_a(P_i^m)\cdot\nabla_{P_i^m}F
\]
acting on functions $F$ of the particles $x,\widetilde x$ and of the
centers $P_i^m$, with the coefficients and shifts of the current
windows frozen: $\B_a$ transports the two particles and the centers
of all finer generations $m>a$ with the generation-$a$ field $u_a$, and
moves nothing else (Definition~\ref{def:operators}, where a query
point is added so that $\B_a$ also acts on the frozen fields $u_b$
themselves). The key analytic fact,
Lemma~\ref{lem:material}, is that differentiating the scale-$b$ field
along the scale-$a$ flow costs one full gradient of the coarser of the
two scales and never a gradient of the finer one. When $a<b$ the
coarse field transports the fine wavelet and the observation point
together, and $\B_au_b$ involves only the increment of $u_a$ across
the support of the wavelet, a gain of a factor $\ell_b/\ell_a$ over
the convective estimate. This is the content of ``controlled
sweeping'', and it is what makes the sums over scales converge.

\medskip\noindent
\emph{A corrected observable} (Section~\ref{sec:corrector}). Once the initial computation are carried out, it becomes clear that that it is more transparent not to study
$\Phi(z_t)$ directly but rather the corrected observable
\[
 Q=\Phi+\sum_bh_b\,\B_b\Phi_b+\sum_{a,b}q_{ab}\,\B_a\B_b\Phi_b,
\]
where $h_b$ is the time remaining in the current scale-$b$ window and
$q_{ab}$ an explicit quadratic function of these phases. The weights
are carefully chosen so that between refreshes $\frac{d}{dt}Q$ cancels exactly
to second order, mixed scales included, leaving a third-order
remainder of size $O(\kappa^2\log^2L)$ and the cutoff error. The
estimates are Littlewood--Paley sums: scales coarser than the
separation $r$ contribute $r\ell_j^{-\beta}$, finer ones contribute
$\ell_j^\alpha$, and both series are geometric away from $r$.

\medskip\noindent
\emph{The gains} (Sections~\ref{sec:renewals} and~\ref{sec:trajectories}). This brings us back to the symbol of the model. At a
refresh of scale $j$ the observable jumps. Averaging over the fresh
data annihilates the first-order term and every mixed term, which are
linear, respectively bilinear, in independent centered families, and
leaves the diagonal term $\tau_j^2\Gamma_j$, which is the sign
computation above. At every time the scale matching the current
separation contributes at rate $\kappa$, all other scales contribute
nonnegatively, and integration in time gives $c\kappa t$. Choosing $L$
as a power of $1/\kappa$ makes the errors a fraction of the gain, and
Jensen's inequality converts the bound on $\E R_t^\beta$ into one on
$\E R_t^2$.

\subsubsection{Anomalous dissipation}\label{sec:ADstrategy}

\emph{The dissipation as a two-point quantity.} Let $\overleftarrow u$
be the time reversal of the field on $[0,t_0]$
(Subsection~\ref{sec:reversal}), and let $\varrho^\varsigma$ solve the
advection--diffusion equation on $\T^d\times\T^d$ in which
$\overleftarrow u$ acts on each of the two variables, with diffusivity
$\mu_J$, starting from the normalized indicator of the
$\varsigma$-neighbourhood of the diagonal. For bounded $\theta_0$, duality and an
approximate identity give (Lemma~\ref{lem:pairidentity})
\[
 \mathbf D_J(\theta_0)=\lim_{\varsigma\downarrow0}\frac12
 \int\!\!\int\big(\theta_0(y)-\theta_0(y')\big)^2
 \varrho^\varsigma(t_0,y,y')\dd y\dd y' .
\]
This is an analytic counterpart of the Lagrangian
fluctuation--dissipation relation in \cite{DE17}.
Indeed, for a fixed velocity realization, let $X_s^x,\widetilde X_s^x$
solve the backward tracer equation with independent Brownian motions:
\[
 \mathrm dX_s^x=-u^{(J)}(X_s^x,t_0-s)\,\mathrm ds
       +\sqrt{2\mu_J}\,\mathrm dW_s,\qquad X_0^x=x,
 \qquad 0\le s\le t_0.
\]
Writing $\mathbb E_W$ for expectation over these Brownian motions only,
the scalar representation and dissipation identity are
\[
 \theta^{(J)}(t_0,x)=\mathbb E_W\theta_0(X_{t_0}^x),\qquad
 \mathbf D_J(\theta_0)=\frac12\int_{\T^d}
 \mathbb E_W\bigl[\bigl(\theta_0(X_{t_0}^x)
              -\theta_0(\widetilde X_{t_0}^x)\bigr)^2\bigr]\dd x,
\]
where the pair density $\varrho^\varsigma_{t_0}$ in the analytic identity
above approximates, as $\varsigma\downarrow0$, the joint law of their
origins, averaged over the common starting point $x$.
Under (S), the mean dissipation of a single Fourier mode $e^{2\pi\mathrm ik\cdot x}$
is, up to a factor, the mean of $\int\sin^2(\pi k\cdot(y-y'))\varrho^\varsigma$
(Proposition~\ref{prop:spectral}). Theorem~\ref{thm:ADgeneral} is thereby
reduced to showing that $\varrho^\varsigma$ moves a definite amount of mass
away from the diagonal by time $t_0$, uniformly in $J$ and $\varsigma$, mode by mode.

\emph{Richardson growth.} Theorem~\ref{thm:pair} is the analogue of
Theorem~\ref{thm:main} for such solutions: if $R=\dist(y,y')$, then
$\E\int R^\beta\varrho(t)\ge c\kappa t-o(1)$ as $J\to\infty$. Its proof
uses the same corrector and jump identities as the trajectory proof:
the equation is
tested against the corrected observable of Section~\ref{sec:corrector},
so that the evolution identity for weak solutions takes the place of the
chain rule along trajectories, and the gains are produced at the renewal
times by averaging over the fresh data. For the diffusing pair, three further ingredients are needed.

(i) \emph{Time reversal.} The construction of the vectorfield $\vec u$ is such that it not time-reversible per se, but its backwards version carries essentially the same properties.  Writing $\tilde{\zeta}$ for a renewal time
in the time variable of the chosen evolution, the averaging there requires
that $\varrho(\tilde{\zeta})$ be independent of the fresh data. For the reversed
field this holds when the lifetimes are nested
(Lemma~\ref{lem:reversed}(c)).

(ii) \emph{The diagonal.} On the diagonal the two velocities coincide,
so the drift cannot separate the pair. We have to show that mass leaves the diagonal by
diffusion. Accordingly, $\Phi$ is replaced by a function $\Phi^\sharp$
with a quadratic core whose radius $r_\sharp$ is a fixed multiple of $\ell_J$
(Subsection~\ref{sec:core}), chosen so that $2\mu\Delta\Phi^\sharp=d\kappa/2$
on the core: there the diffusion produces a gain of the same size
$\kappa$ as the drift produces at every larger scale. Away from the core
the main diffusion term $2\mu\Delta\Phi^\sharp$ is still nonnegative
(Lemma~\ref{lem:core}(e)) and is discarded in the lower bound. The
diffusion acting on the correctors, and the terms created by the spatial
cutoff, are controlled as errors, because the P\'eclet number of the
finest scale is fixed.

(iii) \emph{Localization in scale.} Just as in the argument for the Richardson dispersion, the gain of scale $j$ is collected
at its renewal times. It is always nonnegative, but it is bounded below
by a multiple of $\kappa$ only when the separation at that time lies in
a range comparable to $\ell_j$. This is analogous to the aforementioned estimates for $\Gamma_j$. One therefore needs to
know, in addition,  that the separation at a renewal time of scale $j$ is comparable to
the separation during the following window. For trajectories in the Richardson dispersion argument, this follows from the speed bound \eqref{eq:radiusdet}. Here diffusion
precludes a deterministic speed bound; instead, the H\"older control of
the relative drift allows us to construct exponential supersolutions
and apply the parabolic comparison principle (Lemma~\ref{lem:expmoment}).
Consequently, over times shorter than a fixed multiple of $a^\beta$,
only an exponentially small fraction of mass can move from separation
at most $a$ to separation at least $2a$, or conversely.

\subsection{Organization and conventions}\label{sec:conventions}

Section~\ref{sec:regularity} collects the regularity of the velocity,
its matrix potentials, and the structure function. Section~\ref{sec:testfunctions}
constructs both the logarithmic and quadratic-core observables and proves
their gain estimates. The scale operators and corrector estimates are
developed in Sections~\ref{sec:operators} and~\ref{sec:corrector}, including
the effect of diffusion on the correctors. Section~\ref{sec:renewals}
groups the fresh-data averages, time reversal, and jump identities.
Section~\ref{sec:trajectories} completes the proof of Richardson
dispersion for trajectories. We then turn to anomalous dissipation:
Section~\ref{sec:weak} supplies the energy and duality estimates
(their standard proofs are collected in Appendix~\ref{sec:appendix}), and
Section~\ref{sec:pair} proves separation growth for diffusing pairs.
Section~\ref{sec:AD} deduces anomalous dissipation for the truncations
and, using the comparison of dissipations in
Subsection~\ref{sec:fixedfield}, for the single limiting field.

Constants denoted by $C$ and $c$ may change from line to line. They may
depend on the dimension, on the model parameters $\alpha$, $\psi$,
$R$, $M$, $\sigma$, and on the auxiliary cutoff functions and
numbers fixed once and for all in the course of the proofs: $\vartheta$,
$\chi$ and $\rho$ in Section~\ref{sec:testfunctions}, $\Upsilon$ in
Subsection~\ref{sec:core}, and $\chi_\circ$, $\lambda_\uparrow$,
$\lambda_{\downarrow,\theta}$ in Section~\ref{sec:pair}. These auxiliary
objects are themselves chosen as functions of $d,\alpha,\psi,R$ only
(the number $\rho$ through the threshold $v_*$ of
Lemma~\ref{lem:structure}, the cutoff profiles through $\alpha$, and the
remaining functions independently of all parameters), so that a
constant depending on them depends only on $d,\alpha,\psi,R,M,\sigma$;
this is the dependence stated in the theorems. The constants do not
depend on $J$, on $\kappa$, on the
cutoff parameter $L$ of Section~\ref{sec:testfunctions}, on the
configuration, or on the realization of the random data. Tensor norms
are Euclidean operator norms.

\subsection{Acknowledgements}\label{sec:ackn} The model has been conceived by Raul Jimenez at the workshop in Castro Urdiales mathematics institute in discussions with  C.~De Lellis, D.~N.~Spergel, G.~David, and M.~Engelstein, and the second author of this paper. The authors are immensely thankful to aforementioned collaborators for all the discussions that led to the model and some understanding of its dynamics. The companion physics paper \cite{JMDLSDDE26} is an outgrowth of this effort and the H\"older continuity proof for the small Kubo number in Lemma~\ref{lem:regularity} is due to Camillo De Lellis. The authors used AI for some aspects of writing and to pick a particular wavelet with desired properties; we checked by hand that indeed that  wavelet behaves as asserted. Both authors have been supported by the Simons Foundation, through the Geometry of Flows Initiative.

\section{Velocity estimates and the structure function}\label{sec:regularity}

This section collects the deterministic properties of the model and the properties of the structure function of the wavelet family. Lemma~\ref{lem:nested} records the comparability and nesting of the lifetimes. Lemma~\ref{lem:regularity} establishes H\"older regularity of $u$, uniform in $J$.
Lemma~\ref{lem:structure} establishes the isotropy of $D$, the
K\'arm\'an--Howarth identity and the small-argument expansion of the
structure function for every admissible wavelet;
Proposition~\ref{prop:explicitwavelet} exhibits an explicit wavelet
satisfying \eqref{eq:ND} in every dimension. The matrix-potential estimate of Lemma~\ref{lem:tailpotential}
controls the velocity scales omitted by a finite truncation.
We begin with the elementary properties of the lifetimes
$\tau_j=2^{-\lceil\beta j\rceil}\kappa$ introduced in Section~\ref{sec:model}.

\begin{lemma}[Comparability and nesting]\label{lem:nested}
For every $j\ge j_0$,
\begin{equation}\label{eq:taucomparable}
 \tfrac12\kappa\ell_j^\beta<\tau_j\le\kappa\ell_j^\beta=\kappa/\la_j ,
\end{equation}
and $\tau_j/\tau_{j+1}\in\{1,2\}$. Consequently, for $a\le m$,
$\tau_a/\tau_m$ is a nonnegative integer power of $2$, and every
renewal time $k\tau_a$ of scale $a$ is a renewal time of every finer
scale $m\ge a$.
\end{lemma}

\begin{proof}
Since $\lceil\beta j\rceil\in[\beta j,\beta j+1)$, we have
$2^{-\lceil\beta j\rceil}\in(2^{-\beta j}/2,\,2^{-\beta j}]$, and
$2^{-\beta j}=\ell_j^\beta$. Since $0<\beta<1$, for every real $x$ one
has $\lceil x\rceil\le\lceil x+\beta\rceil\le\lceil x\rceil+1$, so
$\lceil\beta(j+1)\rceil-\lceil\beta j\rceil\in\{0,1\}$. Iterating gives
$\tau_a/\tau_m=2^{\lceil\beta m\rceil-\lceil\beta a\rceil}$, and then
$k\tau_a=(k\tau_a/\tau_m)\tau_m$ with $k\tau_a/\tau_m\in\mathbb N$.
\end{proof}

\begin{lemma}[Packing of wavelets and H\"older regularity of $u$]\label{lem:regularity}
There are $\kappa_*>0$ and $C_*>0$, depending only on
$d,\alpha,\psi,R,M$, such that, for $\kappa\le\kappa_*$,
the centers of each active generation are $\ell_j/2$-separated
for all $t\in[0,1]$.
The number of its centers within distance $R\ell_j$ of any point is
bounded by a constant depending only on $d,R$. For $0\le q\le3$,
\begin{equation}\label{eq:derivativebounds}
 A_j\sum_i\sum_{\nu\in\mathcal A}|\eta_{i,k}^{j,\nu}|
       |D^qW_j^\nu(x-P_{i,k}^j(t))|
 \le C_*A_j\ell_j^{-q}.
\end{equation}
As a result,
\begin{equation}\label{eq:holder}
 |u(x,t)-u(y,t)|\le C_*\dist(x,y)^\alpha,
 \qquad \Lip(u(\cdot,t))\le C_*\la_J.
\end{equation}
For every finite $J$, all center and particle equations have unique
Carath\'eodory solutions, measurable in the random data. Distinct
particles remain distinct. As long as  $R_t\le1/4$,
\begin{equation}\label{eq:radiusdet}
 \left|\frac{\mathrm d}{\mathrm dt}R_t^\beta\right|\le C_H
 \quad\text{for almost every }t.
\end{equation}
Condition \eqref{eq:localtime} ensures $R_s\le1/4$ for all $s\le t$.
\end{lemma}

\begin{proof}
A set with mutual torus distance at least $\ell_j/2$ has at most
$N_{d,R}=\lceil(8R+1)^d\rceil$ points in a ball of radius
$R\ell_j$: disjoint balls of radius $\ell_j/8$ around these points
lie in a ball of radius $(R+1/8)\ell_j\le3/8$; all these balls
have radius less than $1/2$, so their volumes are those of the
corresponding Euclidean balls, and comparing volumes gives the bound.
Set
\[
 K:=N_{d,R}\,\#\mathcal A\,M\max_{0\le q\le3}\max_{\nu\in\mathcal A}
 \|D^qW^\nu\|_\infty .
\]
Since $R\ell_j\le1/4$, for each $x$ and $i$ at most one term of the
periodization $W_j^\nu(x-P_{i,k}^j(t))$ is nonzero, and
$D^qW_j^\nu=\ell_j^{-q}(D^qW^\nu)(\cdot/\ell_j)$ term by term. Hence,
at every time $t$ at which the centers of generation $j$ are
$\ell_j/2$-separated, the packing bound gives
\eqref{eq:derivativebounds} for generation $j$ at time $t$ with $C_*$
replaced by $K$.

We prove by induction on $j\ge j_0$ that, for $\kappa\le\kappa_*$, the
centers of generation $j$ are $\ell_j/2$-separated at every
$t\in[0,1]$, so that \eqref{eq:derivativebounds} holds for generation
$j$ with constant $K$ at every time. Generation $j_0$ is stationary
between its renewals, since $u_{<j_0}=0$, and a shifted grid is
$\ell_{j_0}$-separated. Suppose the claim holds for all generations
$m<j$. Then, for every $t$,
\[
 \Lip(u_{<j}(\cdot,t))\le K\sum_{m<j}\la_m
 \le Kc_\beta\la_j,
 \qquad c_\beta=\sum_{n\ge1}2^{-n\beta}<\infty.
\]
Over an interval of length $\tau_j\le\kappa/\la_j$, the flow of
$u_{<j}$ and its
inverse have Lipschitz constants at most $\exp(Kc_\beta\kappa)$.
Choose $\kappa_*$ so that this number is at most $2$. A newly shifted
grid is $\ell_j$-separated, hence remains $\ell_j/2$-separated
throughout its window, which closes the induction. This proves
\eqref{eq:derivativebounds} with $C_*=K$; the constant is enlarged
below as needed.

Summing the gradient bounds gives the Lipschitz estimate in
\eqref{eq:holder}. Summing
\[
 |u_j(x)-u_j(y)|\le C\min(\ell_j^\alpha,
                     \dist(x,y)\ell_j^{-\beta})
\]
over scales above and below $\dist(x,y)$ gives its H\"older estimate.
There are finitely many windows and variables for fixed $J$. Induction
and the usual Picard construction on the intervals between renewal
times give the asserted existence, uniqueness, and measurability.
The flow and its inverse are Lipschitz with constants at most
$e^{C\la_Jt}$, so $R_t\ge r_0e^{-C\la_Jt}>0$ for a distinct pair.
The chain rule and the H\"older estimate give
\eqref{eq:radiusdet}. Integration up to a first exit from the local
region proves the final assertion.
\end{proof}

\begin{lemma}[The field between renewal times]\label{lem:betweenrenewals}
Let $\kappa\le\kappa_*$, fix $J$ and a realization of the data, and let
$(a,b)$ be an interval between consecutive points of
$(\tau_J\mathbb N\cap[0,1))\cup\{1\}$. On $(a,b)$ all shifts and coefficients are
constant, all centers are $C^1$ in time, and $u^{(J)}$ is $C^3$ in $x$
and $C^1$ in $t$; these derivatives have one-sided limits at $a$ and
$b$. Moreover $\nabla\cdot u^{(J)}=0$ and
$|u^{(J)}|\le C_*\sum_jA_j\le C_*(1-2^{-\alpha})^{-1}$.
\end{lemma}

\begin{proof}
By Lemma~\ref{lem:nested} no window of any scale starts in $(a,b)$.
By induction on scales, the centers of scale $j$ solve
$\dot P=u_{<j}(P,t)$, whose right-hand side is $C^3$ in $P$ and
continuous in $t$ on $[a,b]$ (with one-sided limits); hence the centers
are $C^1$ on $[a,b]$, and
$u_j=A_j\sum_{i,\nu}\eta_i^\nu W_j^\nu(x-P_i(t))$ is $C^3$ in $x$ (as
$W^\nu\in C^3_c$) and $C^1$ in $t$. Each $W^\nu$ is divergence-free,
and the bound follows from \eqref{eq:derivativebounds} with $q=0$.
\end{proof}

\subsection{Matrix potentials of the velocity tail}

For a matrix field $B$, write
$(\operatorname{div}B)_i=\sum_m\partial_mB_{im}$ and use the operator
norm pointwise. The wavelet representation gives a potential for the
tail with one additional power of length.

\begin{lemma}[Matrix potential of the velocity tail]\label{lem:tailpotential}
For the model with $\kappa\le\kappa_*$ and every $J\ge j_0$,
there is a skew-symmetric matrix field $B_{>J}$ on $[0,1]\times\T^d$,
Borel in $(x,t)$ and $C^1$ in $x$, such that
\begin{equation}\label{eq:tailpotential}
\begin{gathered}
 u^{(\infty)}-u^{(J)}=\operatorname{div}B_{>J},\\
 \|B_{>J}\|_{L^\infty_{t,x}}\le C_{\rm pot}\ell_J^{1+\alpha},\qquad
 \|\nabla_xB_{>J}\|_{L^\infty_{t,x}}\le C_{\rm pot}\ell_J^{\alpha}.
\end{gathered}
\end{equation}
The constant $C_{\rm pot}$ depends only on $d,\alpha,\psi,R,M$, and is
independent of $J$, the realization, and $\kappa\le\kappa_*$.
\end{lemma}

\begin{proof}
The wavelet representation provides a matrix potential at each
scale. We sum these potentials after using the packing estimate to
bound them and their gradients. Put
\[
 \psi_j(x):=\sum_{n\in\mathbb Z^d}\psi((x+n)/\ell_j),\qquad
 B_j(x,t):=A_j\ell_j\sum_{i,\nu}\eta_{i,k}^{j,\nu}
 J^\nu\psi_j(x-P^j_{i,k}(t))\quad(t\in I_{jk}).
\]
Since $J^\nu$ is skew-symmetric, so is $B_j$; and since
$\ell_j\nabla\psi_j(x)=\sum_n\nabla\psi((x+n)/\ell_j)$ and
$W^\nu=J^\nu\nabla\psi$,
\[
 \operatorname{div}B_j=A_j\ell_j\sum_{i,\nu}\eta_{i,k}^{j,\nu}
 J^\nu\nabla\psi_j(x-P^j_{i,k}(t))
 =A_j\sum_{i,\nu}\eta_{i,k}^{j,\nu}W_j^\nu(x-P^j_{i,k}(t))=u_j .
\]
By Lemma~\ref{lem:regularity}, at most $N_{d,R}$ centers of
generation $j$ lie within distance $R\ell_j$ of any point, and
$R\ell_j\le\frac14$; hence
\begin{gather*}
 \|B_j\|_\infty\le N_{d,R}\,\#\mathcal A\,M\|\psi\|_\infty\,A_j\ell_j
 =C\ell_j^{1+\alpha},\\
 \|\nabla_x B_j\|_\infty\le N_{d,R}\,\#\mathcal A\,M\|\nabla\psi\|_\infty\,A_j
 =C\ell_j^\alpha .
\end{gather*}
Both series $\sum_{j>J}\ell_j^{1+\alpha}\le\ell_J^{1+\alpha}$ and
$\sum_{j>J}\ell_j^\alpha\le C\ell_J^\alpha$ converge, so
$B_{>J}:=\sum_{j>J}B_j$ converges uniformly together with its
$x$-gradient, is skew-symmetric, and has divergence
$\sum_{j>J}u_j=u^{(\infty)}-u^{(J)}$.
\end{proof}

\subsection{The structure function and an explicit wavelet}

Next, let us discuss the shape of the wavelets and the condition \eqref{eq:ND}. Fix unit vectors $e,n\in\mathbb R^d$ with $n\perp e$. By isotropy,
the quantities
\[
D_L(r)=e^\top D(re)e,
\qquad
D_N(r)=n^\top D(re)n
\]
are independent of the choice of these directions. Explicitly,
\[
D_L(r)=\sum_\nu\int_{\mathbb R^d}
\left|[W^\nu(y)-W^\nu(y-re)]\cdot e\right|^2\,dy,
\]
and
\[
D_N(r)=\sum_\nu\int_{\mathbb R^d}
\left|[W^\nu(y)-W^\nu(y-re)]\cdot n\right|^2\,dy.
\]
Thus $D_L(r)$ measures the integrated squared velocity increment
parallel to the separation vector $re$, whereas $D_N(r)$ measures
the corresponding increment in one perpendicular direction.
There are $d-1$ transverse directions.

For a fresh scale-$j$ random field, multiplication by
$\sigma^2 A_j^2$ and evaluation at $r/\ell_j$ give the corresponding
velocity-increment variances at separation $r$.

The condition $D_L'(r)\ge0$, $r>0,$ says that the longitudinal increment variance does not decrease
as separation increases. Incompressibility and isotropy imply
\[
D_N(r)=D_L(r)+\frac{r}{d-1}D_L'(r).
\]
Consequently, \eqref{eq:ND} is equivalent to
\[
D_N(r)\ge D_L(r),\qquad r>0.
\]

The proof uses \eqref{eq:ND}  to make the averaged second-order renewal gain
nonnegative. Let $F$ be a radial test function, writing
$F(z)=F(r)$ with $r=|z|>0$. Then
\[
D(z):D^2F(z)
=D_L(r)F''(r)+(d-1)D_N(r)\frac{F'(r)}r,
\]
with the colon denoting matrix inner product. 
Taking $F(r)=r^\beta$, where $\beta=1-\alpha\in(0,1)$, gives
\[
D(z):D^2|z|^\beta
=\beta r^{\beta-2}
\left[(d-1)D_N(r)-\alpha D_L(r)\right].
\]
The longitudinal contribution is negative because $r^\beta$
is concave, while the transverse contributions are positive.
Under \eqref{eq:ND}, and using $D_L(r)>0$,
\[
(d-1)D_N(r)-\alpha D_L(r)
\ge(d-1-\alpha)D_L(r)>0.
\]
This is the key sign estimate: the transverse contributions
more than compensate for the longitudinal concavity.

In addition, condition \eqref{eq:ND}  is used for the logarithmically continued
and cutoff test functions. Writing
\[
H_F(r)=rF'(r),
\]
one obtains
\[
D(z):D^2F(z)
=D_L(r)\frac{H_F'(r)}r
+\left[(d-1)D_N(r)-D_L(r)\right]\frac{H_F(r)}{r^2}.
\]
The construction ensures $H_F(r)\ge0$ and $H_F'(r)\ge0$.
Moreover, for $d\ge2$,
\[
(d-1)D_N(r)-D_L(r)
=(d-2)D_L(r)+(d-1)\bigl[D_N(r)-D_L(r)\bigr]\ge0.
\]
Thus both terms in the contraction are nonnegative.
Quantitative nondegeneracy at a comparable scale supplies
the strictly positive gain.

All in all, \eqref{eq:ND} ensures the correct
sign of the renewal gains.

\begin{lemma}[Structure-function properties]\label{lem:structure}
The tensor \eqref{eq:structure} is isotropic, and for $r>0$ it
satisfies the identity
\begin{equation}\label{eq:KH}
 D_N(r)=D_L(r)+\frac{r}{d-1}D_L'(r),\qquad D_L(r)>0.
\end{equation}
In particular, \eqref{eq:ND} is equivalent to $D_N\ge D_L$.

There is a constant $d_0>0$ such that 
\begin{equation}\label{eq:smallstructure}
 D_L(r)=d_0r^2+o(r^2),\qquad
 D_N(r)=\frac{d+1}{d-1}d_0r^2+o(r^2)\quad\textrm{as } r\downarrow0.
\end{equation}
Consequently there are $v_*>0$ and $c_*>0$ such that
\begin{equation}\label{eq:strain}
 D_N(r)-D_L(r)\ge c_*r^2,\quad0<r\le v_*.
\end{equation}
\end{lemma}

\begin{proof} The lemma is standard isotropic incompressible covariance theory (see e.g. \cite{FGV01} which resembles some parts of its argument), but since our model is not close to what has appeared in the literature, we provide a self-contained proof. 
Throughout the proof, derivatives of $C$, $D$, and $H$ are all with
respect to the separation variable $z$. 

First, by definition,
\[
\begin{aligned}
 \sum_{p<q}(J^{pq})_{ak}(J^{pq})_{b\ell}
 &=\frac12\sum_{p,q}
 (\delta_{aq}\delta_{kp}-\delta_{ap}\delta_{kq})
 (\delta_{bq}\delta_{\ell p}-\delta_{bp}\delta_{\ell q})\\
 &=\frac12\bigl(
 \delta_{ab}\delta_{k\ell}-\delta_{a\ell}\delta_{bk}
 -\delta_{a\ell}\delta_{bk}+\delta_{ab}\delta_{k\ell}\bigr)\\
 &=\delta_{ab}\delta_{k\ell}-\delta_{a\ell}\delta_{bk}.
\end{aligned}
\]
Multiplying by $x_ky_\ell$ and summing over $k,\ell$ gives, for
$x,y\in\R^d$,
\[
 \sum_{\nu\in\mathcal A}(J^\nu x)_a(J^\nu y)_b
 =\delta_{ab}(x\cdot y)-x_by_a.
\]
In particular,
\[
 \sum_\nu W_a^\nu(y)W_b^\nu(y-z)
 =\delta_{ab}\sum_k\partial_k\psi(y)\partial_k\psi(y-z)
   -\partial_b\psi(y)\partial_a\psi(y-z).
\]

Define
\[
 H(z)=\int_{\R^d}\psi(y)\psi(y-z)\dd y,
 \qquad
 C_{ab}(z)=\sum_{\nu\in\mathcal A}
       \int_{\R^d}W_a^\nu(y)W_b^\nu(y-z)\dd y.
\]
Since $\psi$ is even, $\psi(y-z)=\psi(z-y)$, so $H=\psi*\psi$.
Differentiating twice in $z$ and integrating
by parts then gives
\[
 \partial_a\partial_bH(z)
 =-\int\partial_b\psi(y)\partial_a\psi(y-z)\dd y,
 \]
yielding the first identity in
\begin{equation}\label{eq:covariancepotential}
 C_{ab}(z)=-\delta_{ab}\Delta H(z)+\partial_a\partial_bH(z),
 \qquad D(z)=2(C(0)-C(z)).
\end{equation}
In particular, $C_{ab}=C_{ba}$. Also $H$ is even, so its second
derivatives, and hence $C$, are even. To verify the second identity,
expand the increments:
\[
\begin{aligned}
 D_{ab}(z)
 ={}&\sum_\nu\int W_a^\nu(y)W_b^\nu(y)\dd y
     +\sum_\nu\int W_a^\nu(y-z)W_b^\nu(y-z)\dd y\\
 &-\sum_\nu\int W_a^\nu(y)W_b^\nu(y-z)\dd y
     -\sum_\nu\int W_a^\nu(y-z)W_b^\nu(y)\dd y\\
 ={}&2(C_{ab}(0)-C_{ab}(z)).
\end{aligned}
\]
 
For any orthogonal matrix $U$, the change of variables $y=Ux$ and
radiality of $\psi$ give
\[
 H(Uz)=\int\psi(Ux)\psi(U(x-z))\dd x
      =\int\psi(x)\psi(x-z)\dd x=H(z).
\]
Thus $H$ is radial; write $H(z)=h(r)$, where $r=|z|$.
For $z\ne0$, set $\widehat z=z/r$. Direct differentiation gives
\[
 \partial_a r=\widehat z_a,
 \qquad
 \partial_a\widehat z_b
 =\frac{\delta_{ab}-\widehat z_a\widehat z_b}{r}.
\]
Consequently,
\[
 \partial_a\partial_bH
 =h''(r)\widehat z_a\widehat z_b
   +\frac{h'(r)}r(\delta_{ab}-\widehat z_a\widehat z_b),
 \qquad
 \Delta H=h''(r)+(d-1)\frac{h'(r)}r.
\]
The covariance formula therefore reads
\[
 C(z)=C_L(r)\widehat z\otimes\widehat z
       +C_N(r)(I-\widehat z\otimes\widehat z),
\]
where
\[
 C_L(r)=-(d-1)\frac{h'(r)}r,
 \qquad C_N(r)=-h''(r)-(d-2)\frac{h'(r)}r.
\]
At the origin, rotational invariance gives $C(0)=cI$ for a scalar
$c$: indeed the Hessian of $H$ at zero is invariant under conjugation
by every orthogonal matrix, so it, and hence $C(0)$, is a scalar
matrix. It follows from \eqref{eq:covariancepotential} that

\[
\begin{gathered}
 D(z)=D_L(r)\widehat z\,\otimes\,\widehat z
       +D_N(r)(I-\widehat z\,\otimes\,\widehat z),\\
 D_L(r)=2(c-C_L(r)),\qquad D_N(r)=2(c-C_N(r)).
\end{gathered}
\]

This proves \eqref{eq:isotropy}, with $D_N$ the eigenvalue in each
of the $d-1$ transverse directions.

Now, each wavelet is divergence-free, since the matrix $J^\nu$ is skew-symmetric and the Hessian of $\psi$ is
symmetric. Therefore,
\[
 \sum_a\partial_{z_a}C_{ab}(z)
 =-\sum_{\nu,a}\int W_a^\nu(y)
                 \partial_{y_a}W_b^\nu(y-z)\dd y=0.
\]
There is no boundary term because the wavelets are compactly
supported. Hence $\sum_a\partial_aD_{ab}=0$. 

Write $A(r)=D_L(r)-D_N(r)$, so that
$D_{ab}=D_N(r)\delta_{ab}+A(r)\widehat z_a\,\widehat z_b$.
The first term has divergence $D_N'(r)\widehat z_b$.
For the second, the product rule gives
\[
\begin{aligned}
 \sum_a\partial_a(A\,\widehat z_a\,\widehat z_b)
 ={}&A'\sum_a(\widehat z_a)^2\,\widehat z_b
 +A\left(\sum_a\partial_a\widehat z_a\right)\widehat z_b
 +A\sum_a\widehat z_a\partial_a\widehat z_b\\
 ={}&A' \,\widehat z_b
 +A\left(\frac{d-1}{r}\right)\widehat z_b
 +0.
\end{aligned}
\]
Hence,
\[
 0=\sum_a\partial_aD_{ab}
 =\left[D_L'(r)+\frac{d-1}{r}(D_L(r)-D_N(r))\right]\widehat z_b.
\]
Since $\widehat z$ is a unit vector, the scalar coefficient is zero.
This proves the identity in \eqref{eq:KH}.
Moreover, $r/(d-1)>0$, so $D_N\ge D_L$ is equivalent to $D_L'\ge0$.
Since $D_L$ is continuously differentiable on $(0,\infty)$, the
latter condition is equivalent to \eqref{eq:ND}.

To prove strict positivity, use isotropy to write
\[
 D_L(r)=\sum_\nu\int_{\R^d}
       [W_1^\nu(y)-W_1^\nu(y-re_1)]^2\dd y.
\]
If this vanished for some $r>0$, nonnegativity and continuity would
give $W_1^\nu(y)=W_1^\nu(y-re_1)$ for every $y,\nu$.
Iteration gives $W_1^\nu(y)=W_1^\nu(y-nre_1)$ for all positive
integers $n$, and compact support then forces $W_1^\nu\equiv0$.
In particular, $W_1^{12}=-\partial_2\psi=0$. Thus $\psi$ is constant
along lines parallel to $e_2$, and compact support forces
$\psi\equiv0$, a contradiction. Therefore $D_L(r)>0$ for all $r>0$.

Finally, put $f_\nu=W_1^\nu$, $g_\nu=\partial_1f_\nu$, and
$\Delta_rf_\nu=f_\nu-f_\nu(\cdot-re_1)$. The fundamental theorem
of calculus and translation invariance of the $L^2$ norm imply
\[
 \Delta_rf_\nu=\int_0^r g_\nu(\cdot-\varrho e_1)\dd\varrho,
 \qquad
 \|g_\nu(\cdot-\varrho e_1)-g_\nu\|_2
 \le\varrho\|\partial_1g_\nu\|_2.
\]
Consequently,
\[
 \Delta_rf_\nu=rg_\nu+E_{\nu,r},
 \qquad \|E_{\nu,r}\|_2\le\tfrac12r^2\|\partial_1g_\nu\|_2.
\]
Expanding the squared norm and using Cauchy--Schwarz gives
\[
 D_L(r)=\sum_\nu\|\Delta_rf_\nu\|_2^2=d_0r^2+O(r^3),
 \qquad
 d_0=\sum_\nu\|g_\nu\|_2^2
    =\sum_\nu\int_{\R^d}|\partial_1W_1^\nu(y)|^2\dd y.
\]
For the derivative, differentiate the integral itself:
\[
\begin{aligned}
 D_L'(r)
 &=2\sum_\nu\langle\Delta_rf_\nu,g_\nu(\cdot-re_1)\rangle\\
 &=2\sum_\nu\langle rg_\nu+O_{L^2}(r^2),
                         g_\nu+O_{L^2}(r)\rangle
 =2d_0r+O(r^2).
\end{aligned}
\]
If $d_0=0$, continuity gives $\partial_1W_1^\nu\equiv0$ for every
$\nu$. Each $W_1^\nu$ is then constant along lines parallel to $e_1$
and hence zero by compact support. As above, $W_1^{12}=0$ forces
$\psi\equiv0$, a contradiction. Thus $d_0>0$.

Substitution in \eqref{eq:KH} now gives
\[
 D_N(r)=d_0r^2+\frac{2d_0}{d-1}r^2+O(r^3)
       =\frac{d+1}{d-1}d_0r^2+o(r^2),
\]
which proves \eqref{eq:smallstructure}. In particular,
\[
 D_N(r)-D_L(r)=\frac{2d_0}{d-1}r^2+o(r^2).
\]
Taking $c_*=d_0/(d-1)$ and choosing $v_*>0$ so that the remainder
has absolute value at most $c_*r^2$ for $0<r\le v_*$ proves
\eqref{eq:strain}.
\end{proof}

The next proposition shows that \eqref{eq:ND} is satisfied by an
explicit wavelet in every dimension. The proof expresses $\psi*\psi$
through the convolution of two nonnegative radial functions in $d+4$
variables; integrating out the four auxiliary variables produces the
two radial derivatives that enter $D_L'$.

\begin{proposition}[Explicit admissible wavelets in every dimension]\label{prop:explicitwavelet}
For each integer $d\ge2$, set
\[
 \psi(x)=(1-|x|^2)_+^5,\quad x\in\R^d,\qquad
 W^\nu=J^\nu\nabla\psi,\quad\nu\in\mathcal A,
\]
where $a_+=\max(a,0)$. This family has the required regularity,
support, and incompressibility ($\psi\in C^4_c$, $W^\nu\in C^3_c$,
$R=1$, hence $j_0=2$), and its combined structure function
satisfies
\[
 D_L'(r)>0\quad(0<r<2),\qquad D_L'(r)=0\quad(r\ge2).
\]
More precisely, let
\[
 G(X)=\frac{20}{\pi^2}(1-|X|^2)_+^3,\quad X\in\R^{d+4},
 \qquad K=G*_{\R^{d+4}}G,
\]
and write $K(X)=k(|X|)$. Then
\begin{equation}\label{eq:explicitND}
 D_L'(r)=8(d-1)\pi^2r\,k(r)\quad(r>0).
\end{equation}
The dimension $d+4$ is auxiliary and is not the dimension of the flow.
\end{proposition}

\begin{proof}
The function $q\mapsto q_+^5$ is $C^4$. Thus $\psi\in C_c^4(\R^d)$,
and
\[
 W^{pq}(x)=-10(1-|x|^2)_+^4(x_pe_q-x_qe_p)
\]
is a nonzero $C_c^3$ field supported in the closed unit ball.
Skew-symmetry of $J^{pq}$ gives zero divergence.

For an integrable function $F$ on $\R^{d+4}=\R^d\times\R^4$,
let $(\mathcal PF)(x)=\int_{\R^4}F(x,y)\dd y$. Since the area
of the unit sphere in $\R^4$ is $2\pi^2$, polar integration gives,
for $a=1-|x|^2>0$,
\begin{multline*}(\mathcal PG)(x)
=\frac{20}{\pi^2}\int_{\mathbb R^4}(a-|y|^2)_+^3\,dy=\frac{20}{\pi^2}\,|S^3|
  \int_0^{\sqrt a}(a-r^2)^3r^3\,dr\\
=40\int_0^{\sqrt a}r^3(a-r^2)^3\dd r
                 =20\int_0^a q(a-q)^3\dd q=a^5.
\end{multline*}
Both sides vanish if $a\le0$ ($|x|\geq 1$). Hence $\mathcal PG=\psi$ and  Fubini's theorem implies
\[
 \mathcal P(G*_{\R^{d+4}}G)
       =(\mathcal PG)*_{\R^d}(\mathcal PG)=\psi*_{\R^d}\psi.
\]

Write $(\psi*\psi)(z)=\tilde h(|z|^2)$; thus $H(z)=\tilde h(|z|^2)$, where
$H=\psi*\psi$ is the function introduced before
\eqref{eq:covariancepotential} (and $\tilde h(r^2)=h(r)$ in the
notation of the proof of Lemma~\ref{lem:structure}). The function $K$ is continuous,
radial, nonnegative, and supported in the closed ball of radius two.
For $q=|z|^2$, the identity $\mathcal PK=\psi*\psi$ and radiality
of $K$ give
\[
\tilde h(q)=\int_{\mathbb R^4}K(z,y)\dd y
=\int_{\mathbb R^4}k\bigl(\sqrt{q+|y|^2}\bigr)\dd y=2\pi^2\int_0^\infty r^3k\bigl(\sqrt{q+r^2}\bigr)\dd r.
\]
Setting $b=q+r^2$, so that
$r^3\dd r=\tfrac12(b-q)\dd b$, 
\begin{equation}\label{eq:marginalconvolution}
\tilde h(q)=\pi^2\int_q^\infty(b-q)k(\sqrt b)\dd b,\qquad q\ge0.
\end{equation}
The fundamental theorem of calculus then assures
\begin{equation}\label{eq:hconvex}
 \tilde h'(q)=-\pi^2\int_q^\infty k(\sqrt b)\dd b,
 \qquad \tilde h''(q)=\pi^2k(\sqrt q)\ge0.
\end{equation}
Since $H(z)=\tilde h(|z|^2)$, the chain rule gives
\[
\partial_i^2H(z)=2\tilde h'(|z|^2)+4z_i^2\tilde h''(|z|^2).
\]
On the other hand, at $z=re_1$, the covariance
identity \eqref{eq:covariancepotential} yields
\[
C_{11}(re_1)
=(-\Delta H+\partial_1^2H)(re_1)
=-\sum_{i=2}^d\partial_i^2H(re_1)
=-2(d-1)\tilde h'(r^2).
\]
Since $D_L(r)=D_{11}(re_1)$ and $D(z)=2(C(0)-C(z))$, it follows that
\[
D_L(r)=2\bigl(C_{11}(0)-C_{11}(re_1)\bigr)
=2C_{11}(0)+4(d-1)\tilde h'(r^2).
\]
Differentiating and using \eqref{eq:hconvex} proves
\eqref{eq:explicitND}. For $0\le r<2$, the open unit balls centered
at $0$ and $re_1$ in $\R^{d+4}$ have a nonempty open intersection.
Both factors of the convolution are strictly positive on that set,
so $k(r)>0$. For $r\ge2$ their interiors are disjoint and $k(r)=0$.
\end{proof}

\section{Auxiliary observables and their gains}\label{sec:testfunctions}

\subsection{Logarithmic continuation and scale cutoffs}

The observable $r^\beta$, $r=|z|$, is modified in two ways before
use. Below the finest scale $\ell_J$ it is continued as a logarithm.
There the field is Lipschitz with constant $\la_J$, and the
continuation keeps the quantity $r\Phi'(r)$ bounded below by a
constant, which is what the gain at scale $J$ requires
(Lemma~\ref{lem:gain}, \eqref{eq:finestgain}); this is the source of
the uniformity in $J$, at the price of the term $\ell_J^\beta$ in
\eqref{eq:maincorrected}. Roughly speaking, we want \(\Phi\) to behave like
$$ \begin{cases} r^\beta,&\text{above a small threshold},\\ \text{constant}+\beta\varepsilon^\beta\log r, &\text{below that threshold}. \end{cases} $$
The transition must preserve the sign of the averaged gain (cf. Lemma~\ref{lem:gain}) and the derivative bounds used later.

Second, for each scale $b$ a function
$\Phi_b$ is formed in which the scales somewhat coarser than $r$, starting roughly $rL$, are switched off; without this, the corrected observable of
Section~\ref{sec:corrector} would contain one term of size
$\kappa r^\beta$ for every coarse scale and would not be uniformly
bounded. The cutoff is nondecreasing in $r$, as the sign computation
of Lemma~\ref{lem:gain} requires.

It is easier to construct $H=r\Phi'$ first. To this end, retain constant $v_*$ from Lemma~\ref{lem:structure}. Fix $\rho>0$ so small that
\begin{equation}\label{eq:rho}
 2\rho<\min(v_*,1/4),\qquad \varepsilon=\rho\ell_J.
\end{equation}
Let $\vartheta$ be a nonnegative even smooth function with integral one,
supported in the interval $(-\log2,\log2)$, and set
$g(v)=\int\vartheta(a)\max(v-a,0)\dd a$. Then
$g(v)\ge\max(v,0)$, $g'\ge0$, $g(v)=0$ for $v\le-\log2$,
and $g(v)=v$ for $v\ge\log2$. Define
\begin{equation}\label{eq:Hphi}
 H(r)=\beta\varepsilon^\beta
       \exp\big(\beta g(\log(r/\varepsilon))\big),
 \qquad
 \Phi(r)=(2\varepsilon)^\beta+
                 \int_{2\varepsilon}^r\frac{H(\varrho)}\varrho\dd \varrho, \quad r>0.
\end{equation}
We also write $\Phi(z)=\Phi(|z|)$ on the punctured local separation
region. Then $H'\ge0$,
\begin{equation}\label{eq:phiproperties}
 \Phi(r)=r^\beta\ (r\ge2\varepsilon),\qquad
 \Phi(r)\le r^\beta\ (r>0),\qquad r\Phi'(r)=H(r)>0.
\end{equation}
Indeed $H(r)\ge\beta r^\beta$, and integrating backwards from
$2\varepsilon$ proves the inequality. All in all, 
\[
\begin{cases}
\Phi(r)=r^\beta, & r\ge 2\varepsilon,\\
\Phi(r)=C_\varepsilon+\beta\varepsilon^\beta\log r,
& 0<r\le \varepsilon/2,
\end{cases}
\]
as desired.
It has no finite extension at $r=0$, since
$\Phi(r)\to-\infty$ as $r\downarrow0$.

Let us also set the notation
\[
 S(r)=\max(r,\ell_J),\qquad \widetilde S(r)=\sqrt{r^2+\ell_J^2};
 \qquad S\le \widetilde S\le\sqrt2S,
\]
to be used below. 
The explicit construction gives, for $k=1,2,3$,
\begin{equation}\label{eq:phiderivatives}
 |\nabla^k\Phi(z)|\le CS(r)^\beta r^{-k}.
\end{equation}
This follows by differentiating $H(r)/r$: all derivatives of $g$
are bounded, since $\vartheta\in C_c^\infty$, $H(r)\le C\max(r,\varepsilon)^\beta$, and each
logarithmic differentiation contributes a factor $r^{-1}$.

Choose a smooth nonincreasing $\chi:[0,\infty)\to[0,1]$ equal to one
on $[0,1]$ and zero on $[2,\infty)$. For $L\ge2$, let
\begin{equation}\label{eq:cutoffs}
 w_j(r)=\chi(\ell_j/(L\widetilde S(r))),\qquad
 \Phi_j(r)=-\int_r^1 w_j(\varrho)\Phi'(\varrho)\dd \varrho.
\end{equation}
The functions $w_j$
are nondecreasing in $r$; they depend on $J$ (through $\widetilde S$)
and on $L$, which the notation suppresses. Their first two derivatives obey
$|w_j^{(k)}(r)|\le Cr^{-k}$ for $k=1,2$, with constants independent
of $j,J,L$. Hence
\begin{equation}\label{eq:phijderivatives}
 |\nabla^k\Phi_j(z)|\le CS^\beta r^{-k}\one_{\{\ell_j\le3LS\}},
 \qquad k=1,2,3.
\end{equation}
The factor $3$ accommodates $2\widetilde S\le2\sqrt2S$. Where the indicator
vanishes, $w_j$ and all its derivatives vanish.

The quantity $\Gamma_j(r)$ introduced in the following lemma is the
averaged jump of the corrected observable when scale $j$ is refreshed while the
particles are at separation $r$ (Proposition~\ref{prop:renewal}). It is one of the core ingredients in the argument.

\begin{lemma}[Positivity and a quantitative one-scale gain]\label{lem:gain}
Define, for $0<r\le1/4$,
\begin{equation}\label{eq:Gamma}
 \Gamma_j(r)=\frac{\sigma^2A_j^2}{2}\,
                    D(z/\ell_j):D^2\Phi_j(z),\qquad |z|=r.
\end{equation}
The functions $ \Gamma_j= \Gamma_j(r)$ are well-defined independently of the direction of $z$ and
$\Gamma_j\ge0$. Moreover, the following lower bounds hold.

There is $c_1>0$, independent of $J,L,\kappa$,
such that
\begin{equation}\label{eq:finestgain}
 \tau_J\Gamma_J(r)\ge c_1\kappa
 \quad\text{if }\quad 0<r\le\min(2\ell_J,1/4).
\end{equation}
For each fixed $0<a<b<\infty$ there is $c_{a,b}>0$ such that
\begin{equation}\label{eq:comparablegain}
 \tau_j\Gamma_j(r)\ge c_{a,b}\kappa
\end{equation}
whenever $r\ge2\varepsilon$, $r/\ell_j\in[a,b]$, $r\le1/4$,
and $w_j(r)=1$.
\end{lemma}

The hypothesis $w_j(r)=1$ in \eqref{eq:comparablegain} holds
whenever $\ell_j\le Lr$, since $\widetilde S(r)\ge r$ and
$\chi=1$ on $[0,1]$. In the application in
Lemma~\ref{lem:count}, we will use it for the selected scale which satisfies
$\ell_j\le2r\le Lr$, so this hypothesis is fulfilled.

\begin{proof}
For a radial function $F(z)=F(r)$, $r=|z|>0$, the Hessian
acts by $F''(r)$ on the radial direction $\widehat z$ and by
$F'(r)/r$ on its $(d-1)$-dimensional orthogonal complement,
as follows by differentiating $\nabla F(z)=F'(r)\,\widehat z$.
Since $z/\ell$ has the same direction as $z$ for $\ell>0$,
\eqref{eq:isotropy} (proved with Lemma~\ref{lem:structure}) gives
\[
 D(z/\ell):D^2F(z)
 =D_L(r/\ell)F''(r)+(d-1)D_N(r/\ell)\frac{F'(r)}r.
\]
The factor $d-1$ counts the transverse directions. Setting
$H_F(r)=rF'(r)$, we have
\[
 \frac{F'(r)}r=\frac{H_F(r)}{r^2},
 \qquad
 F''(r)=\frac{H_F'(r)}r-\frac{H_F(r)}{r^2}.
\]
Then
\begin{equation}\label{eq:radialcontraction}
 D(z/\ell):D^2F(z)
 =D_L(r/\ell)\frac{H_F'(r)}r
  +((d-1)D_N(r/\ell)-D_L(r/\ell))\frac{H_F(r)}{r^2}.
\end{equation}

For $F=\Phi_j$, we have $H_F=r\,\Phi_j'=w_jH\ge0$ and
$H_F'=w_j'H+w_jH'\ge0$, since $w_j,H$ are nonnegative
and nondecreasing by construction.
Moreover,
\[
 (d-1)D_N-D_L=(d-2)D_L+(d-1)(D_N-D_L)\ge0,
\]
so \eqref{eq:ND} and $d\ge2$ prove $\Gamma_j\geq 0$.

We next establish the quantitative bounds. Whenever $r\ge2\varepsilon$ and $w_j(r)=1$, we have
$w_j'(r)=0$, since $w_j$ is differentiable and bounded above
by one. Since $\Phi$ agrees with $r^\beta$ on $[2\varepsilon,\infty)$
and is smooth on $(0,\infty)$, it follows that
\[
 \Phi_j'(r)=w_j(r)\Phi'(r)=\beta r^{\beta-1},
 \qquad
 \Phi_j''(r)=w_j'(r)\Phi'(r)+w_j(r)\Phi''(r)
 =-\alpha\beta r^{\beta-2}.
\]
Substituting into the radial contraction formula gives
\begin{multline}\label{eq:powergain}
 \tau_j\Gamma_j(r)
=\frac{\sigma^2\tau_jA_j^2}{2}
\left[D_L(r/\ell_j)\Phi_j''(r)
+(d-1)D_N(r/\ell_j)\frac{\Phi_j'(r)}r\right]\\
 =\frac{\sigma^2\beta\tau_j\ell_j^{-\beta}}{2}
       (r/\ell_j)^{\beta-2}
       \bigl[(d-1)D_N(r/\ell_j)-\alpha D_L(r/\ell_j)\bigr].
\end{multline}
Define
\[
 K(q)=\frac{\sigma^2\beta}{2}q^{\beta-2}
       \bigl[(d-1)D_N(q)-\alpha D_L(q)\bigr],
 \qquad q>0.
\]
By Lemma~\ref{lem:structure} and \eqref{eq:ND},
\[
 (d-1)D_N(q)-\alpha D_L(q)
 \ge(d-1-\alpha)D_L(q)>0.
\]
Thus $K$ is continuous and strictly positive on $(0,\infty)$.
For any fixed $0<a<b<\infty$, set
\[
 c_{a,b}=\tfrac12\min_{q\in[a,b]}K(q)>0.
\]
Every radius satisfying the hypotheses of
\eqref{eq:comparablegain} has $r\ge2\varepsilon$,
$w_j(r)=1$, and $r/\ell_j\in[a,b]$. Hence
\eqref{eq:powergain} gives
$\tau_j\Gamma_j(r)=(\tau_j\ell_j^{-\beta})K(r/\ell_j)\ge\tfrac12\kappa K(r/\ell_j)\ge c_{a,b}\kappa$,
proving \eqref{eq:comparablegain}.

For \eqref{eq:finestgain}, note that $w_J\equiv1$, since
$\ell_J/(L\widetilde S(r))\le1/L\le1/2$.
For $2\varepsilon\le r\le\min(2\ell_J,1/4)$,
\eqref{eq:powergain} gives
\[
 \tau_J\Gamma_J(r)\ge
 \frac\kappa2\min_{q\in[2\rho,2]}K(q)>0.
\]
For $0<r\le2\varepsilon$, we have $r/\ell_J\le2\rho<v_*$.
Using \eqref{eq:taucomparable}, \eqref{eq:radialcontraction}, $H'\ge0$,
$(d-1)D_N-D_L\ge(d-1)(D_N-D_L)$, and
\eqref{eq:strain}, we obtain
\[
 \tau_J\Gamma_J(r)
 \ge c\kappa\ell_J^{-\beta}H(r)
 \ge c\kappa\rho^\beta,
\]
where the last inequality uses $H(r)\ge\beta\varepsilon^\beta$
and $\varepsilon=\rho\ell_J$.
Taking
$c_1:=\min\big(\tfrac12\min_{q\in[2\rho,2]}K(q),\,c\rho^\beta\big)$,
which is independent of $J,L,\kappa$ because $\rho$ is fixed, the two
bounds prove \eqref{eq:finestgain}.\end{proof}

\subsection{Quadratic continuation at the diagonal}\label{sec:core}

For the diffusing pair we modify $\Phi$ near the diagonal. Its logarithmic singularity is useful for the
trajectory argument, but must be removed before we apply the
advection--diffusion equation. We replace it by a quadratic function
on a ball of radius $r_\sharp$, comparable to $\ell_J$.

The choice of this radius serves two purposes. Inside the quadratic
region, diffusion gives a positive contribution of order $\kappa$;
outside it, the drift gives the gain already established in the
preceding gain estimate. The modification preserves the derivative bounds
needed for the correctors. We record these facts in
Lemma~\ref{lem:core} and then transfer the estimates for the corrected
observable in Corollary~\ref{cor:coretransfer}.

Let $\varepsilon_0\in(0,1]$ and $J\ge j_0$, and set
\begin{equation}\label{eq:diffusivity}
 \mu:=\mu_J=\varepsilon_0\ell_J^{1+\alpha},\qquad
 r_\sharp:=\Big(\frac{4\beta\mu\varepsilon^\beta}{\kappa}\Big)^{1/2}
 =\Big(\frac{4\beta\rho^\beta\varepsilon_0}{\kappa}\Big)^{1/2}\ell_J,
\end{equation}
with $\varepsilon=\rho\ell_J$ as in \eqref{eq:rho}; the second expression
uses $1+\alpha+\beta=2$. We assume throughout
\begin{equation}\label{eq:eps0small}
 \varepsilon_0\le\frac{\kappa\rho^{1+\alpha}}{16\beta},
 \quad\text{equivalently}\quad r_\sharp\le\varepsilon/2 .
\end{equation}
In particular $r_\sharp\le\ell_J/16$, and the ratios $r_\sharp/\ell_J$
and $\varepsilon/r_\sharp$ do not depend on $J$.

To interpolate between the quadratic core and the original function,
fix a nondecreasing $\Upsilon\in C^\infty([0,\infty))$ such that
$\Upsilon(s)=s^2$ on $[0,\frac12]$ and $\Upsilon=1$ on $[1,\infty)$.
In particular, $0\le\Upsilon\le1$ and $\Upsilon(\frac14)=\frac1{16}$.

\begin{definition}\label{def:core}
For $r>0$ let $H^\sharp(r):=\Upsilon(r/r_\sharp)H(r)$, with $H$ from
\eqref{eq:Hphi}, and $H^\sharp(0):=0$; for $0\le r\le1$ let
\[
\begin{aligned}
 \Phi^\sharp(r)&:=(2\varepsilon)^\beta-\int_r^{2\varepsilon}
 \frac{H^\sharp(\varrho)}{\varrho}\dd\varrho,\\
 \Phi^\sharp_j(r)&:=-\int_r^1 w_j(\varrho)\frac{H^\sharp(\varrho)}{\varrho}
 \dd\varrho,
\end{aligned}
\]
with $w_j$ from \eqref{eq:cutoffs} (for $r>2\varepsilon$, $\int_r^{2\varepsilon}$
means $-\int_{2\varepsilon}^r$). Write $\Phi^\sharp(z)=\Phi^\sharp(|z|)$,

$\Phi^\sharp_j(z)=\Phi^\sharp_j(|z|)$, and define the corresponding
gain by
\[
 \Gamma^\sharp_j(r):=\frac{\sigma^2A_j^2}{2}
 D(z/\ell_j):D^2\Phi^\sharp_j(z),\qquad |z|=r.
\]

\end{definition}

For $r=|z|$ put $S=\max(r,\ell_J)$ as before and
\[
 {m_\sharp}={m_\sharp}(r):=\max(r,r_\sharp).
\]

The next lemma separates the properties of the modification from
their use in the evolution argument. Parts (a)--(d) concern the test
functions themselves; parts (e) and (f) give the positive contributions
of diffusion and of the drift, respectively.

\begin{lemma}[Properties of the cored test functions]\label{lem:core}
Assume \eqref{eq:eps0small}. There are constants $C,c_2>0$ depending
only on $d,\alpha,\psi,R,M,\sigma,\rho,\vartheta,\chi,\Upsilon$ (not on
$J,L,\kappa,\varepsilon_0$) such that the following hold.
\begin{enumerate}[(a)]
\item $H^\sharp\ge0$ and $(H^\sharp)'\ge0$ on $(0,\infty)$;
 $\Phi^\sharp=\Phi$ and $\Phi^\sharp_j=\Phi_j$ on $[r_\sharp,1]$;
 $\Phi^\sharp$ is nondecreasing and
 $\Phi^\sharp(r)\le\max(r,r_\sharp)^\beta\le r^\beta+\ell_J^\beta$.
\item For $0\le r\le r_\sharp/2$,
 $H^\sharp(r)=\beta\varepsilon^\beta r^2/r_\sharp^2$,
 \[
\begin{aligned}
 \Phi^\sharp(z)&=\Phi^\sharp(0)+\frac{\beta\varepsilon^\beta}{2r_\sharp^2}|z|^2,\\
 \Phi^\sharp_j(z)&=\Phi^\sharp_j(0)+\frac{\beta\varepsilon^\beta}{2r_\sharp^2}
 \int_0^{|z|^2}\varpi_j(q)\dd q,
\end{aligned}
 \]
 where $\varpi_j(q):=\chi\big(\ell_j/(L(q+\ell_J^2)^{1/2})\big)$.
 Consequently $\Phi^\sharp,\Phi^\sharp_j\in C^\infty(\{|z|<1\})$; in
 particular, they are smooth on the local separation region
 $\{|z|\le1/4\}$, including $z=0$, where they are used.
\item $\Phi^\sharp(0)\ge-\beta\varepsilon^\beta
 \big(1+4^\beta\log(2\varepsilon/r_\sharp)\big)$.
\item For $0\le|z|\le1/4$, $k=1,2,3,4$ and all $j$,
 \begin{equation}\label{eq:corederivs}
\begin{aligned}
 |\nabla^k\Phi^\sharp_j(z)|&\le CS^\beta{m_\sharp}^{-k}
 \one_{\{\ell_j\le3LS\}},\\
 |\nabla\Phi^\sharp_j(z)|&\le CS^\beta r{m_\sharp}^{-2}
 \one_{\{\ell_j\le3LS\}},
\end{aligned}
 \end{equation}
 and the same bounds without the indicator hold for $\Phi^\sharp$.
\item $\Delta\Phi^\sharp\ge0$ on $\{|z|\le1/4\}$, and
 $2\mu\Delta\Phi^\sharp(z)=d\kappa/2$ for $|z|\le r_\sharp/2$.
\item For every $j$ and $0\le r\le1/4$, $\Gamma^\sharp_j(r)\ge0$, and
 $\Gamma^\sharp_j=\Gamma_j$ on $(r_\sharp,1/4]$. Moreover
 \[
 \tau_J\Gamma^\sharp_J(r)\ge c_2\kappa\quad\text{for }
 r\in[r_\sharp/4,\,3\ell_J/2]\cap(0,1/4],
 \]
 and \eqref{eq:comparablegain} holds for $\Gamma^\sharp_j$
 with the same constant $c_{a,b}$.
\end{enumerate}
\end{lemma}

\begin{proof}
We begin with the elementary properties in (a)--(c), then prove the
derivative estimates and the lower bounds.

For (a), $\Upsilon$ and $H$ are nonnegative and nondecreasing ($H'\ge0$ by
\eqref{eq:phiproperties}), hence so is their product. Since
$\Upsilon=1$ on $[1,\infty)$, $H^\sharp=H$ on $[r_\sharp,\infty)$, which
gives the two identities by comparison with \eqref{eq:Hphi} and
\eqref{eq:cutoffs}. Monotonicity follows from
$(\Phi^\sharp)'=H^\sharp/r\ge0$. For $r\ge r_\sharp$,
$\Phi^\sharp(r)=\Phi(r)\le r^\beta$ by \eqref{eq:phiproperties}; for
$r<r_\sharp$, $\Phi^\sharp(r)\le\Phi^\sharp(r_\sharp)\le r_\sharp^\beta$;
finally $r_\sharp\le\ell_J$.

To prove (b), suppose that $r\le r_\sharp/2\le\varepsilon/4$. We have
$\log(r/\varepsilon)\le-\log2$, so $g=0$ and $H=\beta\varepsilon^\beta$;
and $\Upsilon(r/r_\sharp)=(r/r_\sharp)^2$. The formula for $\Phi^\sharp$
follows by integration. For $\Phi^\sharp_j$,
$w_j(\varrho)=\varpi_j(\varrho^2)$ by \eqref{eq:cutoffs}, and
$\int_0^rw_j(\varrho)\varrho\dd\varrho=\frac12\int_0^{r^2}\varpi_j$.
Since $\chi$ is smooth and $q\mapsto(q+\ell_J^2)^{-1/2}$ is smooth on
$(-\ell_J^2,\infty)$, $\varpi_j$ is smooth there, so both functions are
smooth functions of $|z|^2$ near $0$. Away from $0$ they are radial
with smooth profiles ($g,\vartheta,\chi,\Upsilon$ are smooth).

For the lower bound in (c), split the defining integral at
$r_\sharp$. On $(0,r_\sharp]$, $H=\beta\varepsilon^\beta$ (as $r_\sharp\le\varepsilon/2$), so
\[
\begin{aligned}
 \int_0^{r_\sharp}H^\sharp/\varrho\dd\varrho
 &=\beta\varepsilon^\beta\int_0^1\Upsilon(s)s^{-1}\dd s\\
 &\le\beta\varepsilon^\beta(\tfrac18+\log2)
 \le\beta\varepsilon^\beta.
\end{aligned}
\]

Since $g(v)\le\max(v,0)+\log2$, we have
$H(\varrho)\le\beta2^\beta\max(\varrho,\varepsilon)^\beta\le\beta4^\beta\varepsilon^\beta$
for $\varrho\le2\varepsilon$, hence
\[
 \int_{r_\sharp}^{2\varepsilon}H/\varrho\dd\varrho\le\beta4^\beta\varepsilon^\beta\log(2\varepsilon/r_\sharp).
\]

Combining these estimates and discarding the nonnegative term
$(2\varepsilon)^\beta$ proves (c).

We turn to the derivative estimates in (d). Write $F=\Phi^\sharp_j$,
$F(z)=f(|z|)$, and $\mathcal H:=w_jH^\sharp=rf'$. The proof for
$\Phi^\sharp$ is identical, with $w_j\equiv1$. Away from the origin we
use radial differentiation; inside the quadratic core we instead
differentiate the smooth function of $|z|^2$ supplied by (b).

\emph{Case $r\ge r_\sharp/2$.} Here ${m_\sharp}\le2r$. Radial
differentiation gives
\[
 \nabla^kF(z)=\sum_{m=1}^k f^{(m)}(r)T_{k,m}(z),\qquad
 |\nabla^kF(z)|\le C\sum_{m\le k}|f^{(m)}(r)|r^{m-k},
\]
where $T_{k,m}$ is smooth away from the origin and homogeneous of
degree $m-k$. These identities follow by induction, using
$\nabla f^{(m)}(r)=f^{(m+1)}(r)\widehat z$. Since
$f^{(m)}=(\mathcal H/r)^{(m-1)}$, it is enough to show that
\[
 |r^i\mathcal H^{(i)}(r)|\le CS^\beta,\qquad 0\le i\le3.
\]
We verify this bound for the three factors of $\mathcal H=w_j\Upsilon H$.
For $H$, the definition gives
\[
 (r\tfrac{d}{dr})^iH=H\cdot P_i(g',g'',\dots),\qquad
 H\le\beta2^\beta\max(r,\varepsilon)^\beta\le CS^\beta,
\]
where $P_i$ is a polynomial and the derivatives of $g$ are bounded.
For the core cutoff, write $s=r/r_\sharp$. Then
$r^i\frac{d^i}{dr^i}\Upsilon(r/r_\sharp)=s^i\Upsilon^{(i)}(s)$ is bounded:
for $i\ge1$, the derivatives of $\Upsilon$ are supported in $[0,1]$.
Finally, for $w_j$, put $\phi=\ell_j/(L\widetilde S)$. The identities
\[
 r\phi'=-\phi r^2/\widetilde S^2,\qquad
 r\frac{d}{dr}(r^2/\widetilde S^2)=2r^2\ell_J^2/\widetilde S^4
\]
show that $(r\frac d{dr})^iw_j$ is a sum of terms
$\chi^{(p)}(\phi)\phi^p$ multiplied by bounded factors, with $p\ge1$.
Since $\phi\in[1,2]$ on $\supp\chi^{(p)}$, these terms are bounded as
well. The Leibniz rule now proves the required estimate for
$\mathcal H$. Consequently,
\[
\begin{aligned}
 |\nabla^kF|&\le CS^\beta r^{-k}
       \le C2^kS^\beta{m_\sharp}^{-k},\\
 |\nabla F|&=\mathcal H/r\le CS^\beta r^{-1}
       \le4CS^\beta r{m_\sharp}^{-2}.
\end{aligned}
\]

\emph{Case $r<r_\sharp/2$.} Here ${m_\sharp}=r_\sharp$, $S=\ell_J$, and by
(b) $F(z)=G(|z|^2)$ with
$G'=\frac{\beta\varepsilon^\beta}{2r_\sharp^2}\varpi_j$. With
$\bar\phi(q)=\ell_j/(L(q+\ell_J^2)^{1/2})$ one has
$|\bar\phi^{(n)}(q)|\le C_n\bar\phi\,\ell_J^{-2n}$ for $q\ge0$, and
$\bar\phi\in[1,2]$ on $\supp\chi^{(p)}$, so $|\varpi_j|\le1$ and
$|\varpi_j^{(n)}|\le C\ell_J^{-2n}$ for $n\ge1$. Moreover
$|\nabla^k[G(|z|^2)]|\le C\sum_{\lceil k/2\rceil\le p\le k}|G^{(p)}(|z|^2)|\,r^{2p-k}$.
Using $r\le r_\sharp\le\ell_J$ and $\varepsilon\le\ell_J=S$,
we obtain the following estimates for the four derivative orders:
\begin{align*}
 |\nabla F|&\le C\varepsilon^\beta r/r_\sharp^2
       \le CS^\beta r{m_\sharp}^{-2},\\
 |\nabla^2F|&\le C\varepsilon^\beta r_\sharp^{-2}(1+r^2\ell_J^{-2})
       \le CS^\beta{m_\sharp}^{-2},\\
 |\nabla^3F|&\le C\varepsilon^\beta r_\sharp^{-2}\ell_J^{-2}r
       (1+r^2\ell_J^{-2})
       \le C\varepsilon^\beta r_\sharp^{-3},\\
 |\nabla^4F|&\le C\varepsilon^\beta r_\sharp^{-2}\ell_J^{-2}
       (1+r^2\ell_J^{-2}+r^4\ell_J^{-4})
       \le C\varepsilon^\beta r_\sharp^{-4}.
\end{align*}

It remains to check the support restriction in \eqref{eq:corederivs}.
If $\ell_j>3LS$, then $\ell_j/(L\widetilde S)>2$ on a
neighbourhood of $r$, so $w_j$ and all derivatives of $\Phi^\sharp_j$
vanish there, as in \eqref{eq:phijderivatives}.

For (e), the radial Laplacian gives, when $r>0$,  $\Delta\Phi^\sharp=(H^\sharp)'/r+(d-2)H^\sharp/r^2\ge0$ by
(a) and $d\ge2$. For $|z|\le r_\sharp/2$, (b) gives
$\Delta\Phi^\sharp=d\beta\varepsilon^\beta/r_\sharp^2=d\kappa/(4\mu)$ by
\eqref{eq:diffusivity}; by continuity $\Delta\Phi^\sharp\ge0$ at $0$.

Finally, we prove (f). Nonnegativity follows from the same argument
as in Lemma~\ref{lem:gain}: that proof uses only
$w_jH\ge0$ and $(w_jH)'\ge0$, which hold for $w_jH^\sharp$ by (a); at
$r=0$, $D(0)=0$. On $(r_\sharp,1/4]$, $\Phi^\sharp_j=\Phi_j$ on an open
interval around $r$, so the Hessians agree. For
$r\in[r_\sharp/4,2\varepsilon)$: $w_J\equiv1$ (proof of
Lemma~\ref{lem:gain}), $r/\ell_J<2\rho<v_*$, and by
\eqref{eq:radialcontraction}, $(H^\sharp)'\ge0$,
$(d-1)D_N-D_L\ge(d-1)(D_N-D_L)$ and \eqref{eq:strain},
\[
 \tau_J\Gamma^\sharp_J(r)\ge\frac{\sigma^2\tau_JA_J^2}{2}(d-1)c_*
 \ell_J^{-2}H^\sharp(r)
 \ge\frac{\sigma^2(d-1)c_*}{4}\kappa\ell_J^{-\beta}\cdot
 \frac{\beta\varepsilon^\beta}{16}
 =\frac{\sigma^2(d-1)c_*\beta\rho^\beta}{64}\kappa,
\]
using \eqref{eq:taucomparable}, $A_J^2\ell_J^{-2}=\ell_J^{-2\beta}$ and
$H^\sharp(r)\ge\Upsilon(\frac14)H(r)\ge\beta\varepsilon^\beta/16$. For
$r\in[2\varepsilon,3\ell_J/2]$, $\Gamma^\sharp_J=\Gamma_J$ and
\eqref{eq:powergain} gives
$\tau_J\Gamma_J(r)=(\tau_J\ell_J^{-\beta})K(r/\ell_J)
\ge\frac\kappa2\min_{[2\rho,3/2]}K$. Take $c_2$ the smaller of the two
constants. The last claim follows from $\Gamma^\sharp_j=\Gamma_j$ for
$r\ge2\varepsilon>r_\sharp$ and Lemma~\ref{lem:gain}.
\end{proof}

\section{Commutator estimates for the scale operators}\label{sec:operators}

The estimates of this section are deterministic: they concern one
fixed configuration of particles, centers, and coefficients, and
involve no averaging and no time evolution. They bound the derivative
of the scale-$b$ field along the flow of scale $a$, as seen by two
particles at distance $r$. Lemma~\ref{lem:material} shows that each
such derivative costs one full gradient of the coarser of the two
scales and never a gradient of the finer one; for $a<b$ this is a
commutator estimate, since the coarse flow moves the fine wavelet and
the observation point together and only the increment of $u_a$ across
the support of the wavelet, of size $\la_a\ell_b$, survives. We first
isolate the class of configurations to which the estimates apply;
Lemma~\ref{lem:regularity} shows that the process never leaves this
class (Remark~\ref{rem:attained}).

\begin{definition}[Admissible configurations, frozen fields]
\label{def:admissible}
Fix $d\ge2$, $\alpha\in(0,1)$, $\psi$, $R$, $M$, and $j_0\le J$ as in
Section~\ref{sec:model}. A \emph{configuration} is a triple
$\mathcal S=(x,\widetilde x,\Theta)$, where $x,\widetilde x\in\T^d$ and
\[
 \Theta=\Big(\big(P_i^j\big)_{i\in I_j},\
             \big(\eta_i^{j,\nu}\big)_{i\in I_j,\ \nu\in\mathcal A}
       \Big)_{j_0\le j\le J},
 \qquad P_i^j\in\T^d,\quad \eta_i^{j,\nu}\in\R .
\]
The configuration is \emph{admissible} if
\begin{enumerate}[(C1)]
\item $|\eta_i^{j,\nu}|\le M$ for all $j_0\le j\le J$, $i\in I_j$,
      $\nu\in\mathcal A$;
\item $\dist(P_i^j,P_{i'}^j)\ge\ell_j/2$ for all $j_0\le j\le J$ and all
      $i\ne i'$ in $I_j$;
\item $0<r:=\dist(x,\widetilde x)\le1/4$.
\end{enumerate}
For a configuration $\mathcal S$ and $j_0\le j\le J$, the
\emph{frozen field of generation $j$} is
\begin{equation}\label{eq:frozenfield}
 u_j(y)=u_j(y;\Theta)
 =A_j\sum_{i\in I_j}\sum_{\nu\in\mathcal A}
   \eta_i^{j,\nu}\,W_j^\nu\big(y-P_i^j\big),\qquad y\in\T^d .
\end{equation}
Since $W^\nu\in C^3_c(\R^d;\R^d)$, each $u_j$ is a $C^3$ vector field on
$\T^d$, and it depends on $\Theta$ only through the generation-$j$ data
$(P_i^j)_i$, $(\eta_i^{j,\nu})_{i,\nu}$, on which it depends in a $C^3$
manner as well. By (C2), the packing argument in the proof of
Lemma~\ref{lem:regularity} applies to the centers $(P_i^j)_i$ of every
admissible configuration: at most $N_{d,R}$ of them lie within distance
$R\ell_j$ of any given point. Together with (C1) and
$\|D^qW^\nu\|_\infty<\infty$ for $q\le3$, this gives the bound
\eqref{eq:derivativebounds} with $P_{i,k}^j(t)$ replaced by $P_i^j$,
and in particular
\begin{equation}\label{eq:frozenderivatives}
 \|D^q u_j\|_{L^\infty(\T^d)}\le C_*A_j\ell_j^{-q},\qquad 0\le q\le3,
\end{equation}
with $C_*$ depending only on $d,\alpha,\psi,R,M$.
\end{definition}

\begin{definition}[Scale operators and two-point differences]
\label{def:operators}
Let $\mathcal E=\T^d\times\T^d\times\T^d\times\prod_{j=j_0}^J(\T^d)^{I_j}$
be the space of tuples $(y;x,\widetilde x;P)$, $P=(P_i^j)_{j,i}$,
consisting of a \emph{query point} $y$ together with the positional
data of a configuration; the coefficients $(\eta_i^{j,\nu})$ are held
fixed and regarded as parameters. For $j_0\le a\le J$ define the
$C^3$ vector field $V_a$ on $\mathcal E$ by
\[
 V_a(y;x,\widetilde x;P)
 =\Big(u_a(y);\ u_a(x),\ u_a(\widetilde x);\
       \big(\one_{\{m>a\}}\,u_a(P_i^m)\big)_{m,i}\Big),
\]
and the first-order differential operator $\B_a:=V_a\cdot\nabla$ acting
on $C^1$ functions $F$ on $\mathcal E$ (componentwise on vector- or
tensor-valued $F$):

\begin{equation}\label{eq:Bdefinition}
\begin{aligned}
 \B_aF={}&u_a(y)\cdot\nabla_yF
      +u_a(x)\cdot\nabla_xF+u_a(\widetilde x)\cdot\nabla_{\widetilde x}F\\
 &+\sum_{m>a}\sum_{i\in I_m}u_a(P_i^m)\cdot\nabla_{P_i^m}F,
 \qquad \B:=\sum_{a=j_0}^J\B_a .
\end{aligned}
\end{equation}

Thus $\B_a$ transports the query point and the two particles with the
generation-$a$ field and transports the centers of every finer
generation $m>a$ with the same field; it does not move the centers of
generation $a$ or of any coarser generation. On functions that do not
depend on the query point, $\B$ is the infinitesimal generator of the
deterministic dynamics \eqref{eq:centers}--\eqref{eq:particles} between
two consecutive renewal times, with the active coefficients and
shifts frozen.

Each frozen field $u_b$ can be regarded as the function
$(y;x,\widetilde x;P)\mapsto u_b(y;\Theta)$ on $\mathcal E$, which is
$C^3$. Consequently $\B_au_b\in C^2(\mathcal E;\R^d)$ and
$\B_c\B_au_b\in C^1(\mathcal E;\R^d)$ are well defined for all
$a,b,c\in\{j_0,\ldots,J\}$, and $\B_a$ obeys the product rule
$\B_a(FG)=(\B_aF)G+F\,\B_aG$.

For a continuous function $G$ on $\mathcal E$ and a configuration
$\mathcal S=(x,\widetilde x,\Theta)$, the \emph{two-point difference}
of $G$ at $\mathcal S$ is
\begin{equation}\label{eq:twopoint}
 \delta G:=G(x;x,\widetilde x;P)-G(\widetilde x;x,\widetilde x;P),
\end{equation}
that is, the query point is evaluated at the two particles while all
other coordinates are held fixed. In particular
$\delta u_b=u_b(x)-u_b(\widetilde x)$, and $\delta(\B_au_b)$,
$\delta(\B_c\B_au_b)$ are the two-point differences of the functions
$\B_au_b$ and $\B_c\B_au_b$ on $\mathcal E$.
\end{definition}

For an admissible configuration with separation $r=\dist(x,\widetilde x)$
set
\begin{equation}\label{eq:rate-defs}
 \gamma_j=\frac{A_j}{r}\min(1,r/\ell_j)=\min(\la_j,A_j/r),
 \qquad
 K_{ab}=\la_{\min(a,b)}\,\gamma_{\max(a,b)},
\end{equation}
where  $j_0\le j,a,b\le J$. Evidently, $K_{ab}$ is one full
gradient of the coarser of the two scales times the increment factor $\gamma$ of
the finer one.
Writing the full expressions, we have
\[
\gamma_j
=
\begin{cases}
\lambda_j=2^{\beta j}=2^{(1-\alpha)j},
& r\le 2^{-j},\\[4pt]
\dfrac{A_j}{r}=\dfrac{2^{-\alpha j}}{r}
=\dfrac{2^{(\beta-1)j}}{r},
& r>2^{-j},
\end{cases}
\]
and
\[
K_{ab}
=
\begin{cases}
\lambda_a\lambda_b
=2^{\beta(a+b)}
=2^{(1-\alpha)(a+b)},
& a<b,\quad r\le 2^{-b},\\[4pt]
\dfrac{\lambda_a A_b}{r}
=\dfrac{2^{\beta a-\alpha b}}{r}
=\dfrac{2^{(1-\alpha)a-\alpha b}}{r},
& a<b,\quad r>2^{-b},\\[6pt]
\lambda_b\lambda_a
=2^{\beta(a+b)}
=2^{(1-\alpha)(a+b)},
& a\ge b,\quad r\le 2^{-a},\\[4pt]
\dfrac{\lambda_b A_a}{r}
=\dfrac{2^{\beta b-\alpha a}}{r}
=\dfrac{2^{(1-\alpha)b-\alpha a}}{r},
& a\ge b,\quad r>2^{-a}.
\end{cases}
\]

\begin{lemma}[Two material derivatives of a wavelet field]
\label{lem:material}
Let $\mathcal S=(x,\widetilde x,\Theta)$ be an admissible configuration
in the sense of Definition~\ref{def:admissible}, with
$r=\dist(x,\widetilde x)>0$. Let $u_j$ be the frozen fields
\eqref{eq:frozenfield}, let $\B_a$ and $\delta$ be as in
Definition~\ref{def:operators}, and define $\gamma_j,K_{ab}$ by
\eqref{eq:rate-defs} with this $r$. There is a constant $C>0$,
depending only on $d,\alpha,\psi,R,M$, such that, for all
$a,b,c\in\{j_0,\ldots,J\}$,
\begin{align}
 |\delta u_b|&\le C\,r\,\gamma_b,\label{eq:material0}\\
 |\delta(\B_a u_b)|&\le C\,r\,K_{ab},\label{eq:material1}\\
 |\delta(\B_c\B_a u_b)|
 &\le C\,r\,\la_{\min(a,b)}K_{c,\max(a,b)}.
 \label{eq:material2}
\end{align}
The operators act on the enlarged configuration space before the
query point is evaluated at $x$ and $\widetilde x$ to take $\delta$.
The constant is independent of $J,\kappa,L$, the configuration,
and the indices $a,b,c$.
\end{lemma}

\begin{remark}[Meaning of the iterated derivatives]\label{rem:Bderiv}
All coefficients are held fixed. In \eqref{eq:material1}, we regard
$u_b$ as the function
\[
 (y;x,\widetilde x;P)\longmapsto
 u_b(y;P^b)=A_b\sum_{i,\nu}\eta_i^{b,\nu}W_b^\nu(y-P_i^b),
 \qquad P^b=(P_i^b)_i.
\]
It depends on $y$ and $P^b$, but not on the separate particle
coordinates $x,\widetilde x$. The function $F_{ab}:=\B_a u_b$
is then defined on the same space. As shown in \eqref{eq:Fformula},
\[
 F_{ab}(y;P)=
 \begin{cases}
 A_b\displaystyle\sum_{i,\nu}\eta_i^{b,\nu}
 [u_a(y;P^a)-u_a(P_i^b;P^a)]\cdot\nabla W_b^\nu(y-P_i^b),&a<b,\\[3pt]
 (u_a(y;P^a)\cdot\nabla_y)u_b(y;P^b),&a\ge b.
 \end{cases}
\]
Thus $F_{ab}$ depends only on $y$ and the centers of generations
$a,b$, and remains independent of $x,\widetilde x$. For $a\ne b$,
\[
\begin{aligned}
 \B_cF_{ab}
 ={}&(u_c(y)\cdot\nabla_y)F_{ab}\\
 &+\one_{\{c<a\}}\sum_k
 (u_c(P_k^a)\cdot\nabla_{P_k^a})F_{ab}+\one_{\{c<b\}}\sum_i
 (u_c(P_i^b)\cdot\nabla_{P_i^b})F_{ab}.
\end{aligned}
\]
If $a=b$, there is only one family of centers, and its contribution
is included once. Each derivative acts on every occurrence of its
coordinate: for example, when $a<b$, the $P_i^b$ derivative acts
both on $u_a(P_i^b;P^a)$ and on the translated wavelet.
Only after forming $\B_cF_{ab}$ do we take its two-point difference,
by evaluating $y$ at $x$ and $\widetilde x$ with all other coordinates
fixed.

The two cases in \eqref{eq:Fformula} also explain the scale factor
$K_{ab}$. When $a<b$, the coarse field moves the query point and the
fine center together, leaving an increment across a distance of
order $\ell_b$. When $a\ge b$, the spatial derivative falls directly
on the coarser field $u_b$. The proof below makes these two
mechanisms explicit at the next derivative as well.
\end{remark}

\begin{proof}[Proof of Lemma~\ref{lem:material}]
Throughout the proof the admissible configuration $\mathcal S$ is
fixed. For a function $G$ on $\mathcal E$, $\|G\|_\infty$ denotes the
supremum over the query point $y\in\T^d$ with the positional data of
$\mathcal S$ fixed, $D^qG$ denotes derivatives in $y$, and the
Lipschitz constant of $G$ refers to the variable $y$. Since
$\dist(x,\widetilde x)=r\le1/4$, the two points are joined by a
minimizing geodesic of length $r$, and therefore every such $G$
satisfies
\begin{equation}\label{eq:twopointbound}
 |\delta G|\le\min\big(2\|G\|_\infty,\ r\,\|\nabla_yG\|_\infty\big).
\end{equation}
Applied to $D^qu_j$ and combined with \eqref{eq:frozenderivatives}
for $q$ and $q+1$, this gives, for $0\le q\le2$,
\begin{equation}\label{eq:spatialincrements}
 |\delta D^qu_j|\le CA_j\ell_j^{-q}\min(1,r/\ell_j).
\end{equation}
Equation \eqref{eq:material0} is then contained in
\eqref{eq:spatialincrements}.

Products are estimated using
$\delta(FG)=(\delta F)G(x)+F(\widetilde x)\delta G$.

Write $F_{ab}(y)=\B_a u_b(y)$. As a function on the enlarged configuration space, $u_b(y)$ depends
only on the query point $y$ and the scale-$b$ centers $P_i^b$; it
does not depend on the particle coordinates $x,\widetilde x$ or
the other centers. Moreover, the field $u_b$ depends on the
positional data only through the centers $P_i^b$, and by
\eqref{eq:Bdefinition} the operator $\B_a$ contains the derivatives
$\nabla_{P_i^b}$ if and only if $a<b$. Since
$\nabla_{P}\big[W_b^\nu(y-P)\big]=-\nabla W_b^\nu(y-P)$,
\eqref{eq:Bdefinition} gives
\begin{equation}\label{eq:Fformula}
 F_{ab}(y)=
 \begin{cases}
 A_b\displaystyle\sum_{i,\nu}\eta_i^{b,\nu}
 [u_a(y)-u_a(P_i^b)]\cdot\nabla W_b^\nu(y-P_i^b),&a<b,\\
 (u_a(y)\cdot\nabla)u_b(y),&a\ge b.
 \end{cases}
\end{equation}
Start with $a<b$. Write $v_i(y)=u_a(y)-u_a(P_i^b)$. On the support of each wavelet,
\[
|v_i|\le C\lambda_a\ell_b,\qquad
|Dv_i|\le C\lambda_a,\qquad
|D^2v_i|\le C\lambda_a/\ell_a
\le C\lambda_a/\ell_b,
\]
where the last inequality uses $\ell_a\ge\ell_b$.

By the product rule, a term in $D^q(v_i\cdot\nabla W_b^\nu)$
with $m$ derivatives on $v_i$ is bounded by
\[
C\lambda_a\ell_b^{1-m}\,
\ell_b^{-(q-m+1)}
=C\lambda_a\ell_b^{-q},
\qquad 0\le m\le q\le2.
\]
Multiplying by $A_b$, using the coefficient bound, and summing
over the uniformly bounded number of wavelets meeting each point
proves
\begin{equation}\label{eq:Fcoarse}
 \|D^qF_{ab}\|_\infty\le CA_b\la_a\ell_b^{-q},
 \quad q=0,1,2,\qquad a<b.
\end{equation}
The cases $q=0$ and $q=1$ of \eqref{eq:Fcoarse},
together with \eqref{eq:twopointbound}, give
\[
\begin{aligned}
|\delta F_{ab}|
&\le \min\!\left(2\|F_{ab}\|_\infty,\,
                  r\|\nabla F_{ab}\|_\infty\right)\\
&\le CA_b\lambda_a\min(1,r/\ell_b)\\
&=Cr\lambda_a\gamma_b
=CrK_{ab}, \quad a<b.
\end{aligned}
\]
The last equality uses $K_{ab}=\lambda_a\gamma_b$ when $a<b$.
This proves \eqref{eq:material1} in this case.

For $a\ge b$, the second line of \eqref{eq:Fformula} gives
\[
F_{ab}=(u_a\cdot\nabla)u_b.
\]
The product rule for differences and \eqref{eq:spatialincrements} yield
\[
\begin{aligned}
|\delta F_{ab}|
&\le |\delta u_a|\,\|\nabla u_b\|_\infty
   +\|u_a\|_\infty\,|\delta\nabla u_b|\\
&\le CA_a\lambda_b
\left[\min(1,r/\ell_a)+\min(1,r/\ell_b)\right]\\
&\le CA_a\lambda_b\min(1,r/\ell_a)
=Cr\lambda_b\gamma_a
=CrK_{ab},\quad a\geq b,
\end{aligned}
\]
and hence \eqref{eq:material1} for $a\geq b$. We moreover record that
\[
\|F_{ab}\|_\infty
\le \|u_a\|_\infty\|\nabla u_b\|_\infty
\le CA_a\lambda_b,\quad a\geq b.
\]
Hence,
\begin{equation}\label{eq:Fsup}
 \|F_{ab}\|_\infty\le C\la_{\min(a,b)}A_{\max(a,b)}.
\end{equation}

We now prove \eqref{eq:material2}. For this part of the proof, write
\[
 \mathsf m_j=\min(1,r/\ell_j),\qquad A_j\mathsf m_j=r\gamma_j.
\]
If $j\ge k$, then $\ell_j\le\ell_k$ and $\mathsf m_k\le\mathsf m_j$.
We repeatedly use \eqref{eq:twopointbound} to pass from a supremum
and a spatial gradient bound to a two-point estimate. All gradients
below are in the query variable, with the centers fixed.

\medskip\noindent
\emph{Case 1: $a<b$ and $c<b$.}
Here the desired bound is $Cr\la_a\la_c\gamma_b$.
By definition, $\B_c[u_a(y)]=F_{ca}(y)$. To differentiate
$u_a(P_i^b)$, however, we must account for both the evaluation
point and the centers defining $u_a$. This expression is independent
of the query coordinate, so only center derivatives contribute.
Since $c<b$, the center $P_i^b$ moves with velocity $u_c(P_i^b)$;
and since $a<b$, it is distinct from the scale-$a$ coordinates.
The chain rule therefore gives
\[
\begin{aligned}
 \B_c[u_a(P_i^b)]
 ={}&\left.(u_c(P_i^b)\cdot\nabla_y)u_a(y)\right|_{y=P_i^b}\\
 &+\one_{\{c<a\}}\sum_k
 \left.(u_c(P_k^a)\cdot\nabla_{P_k^a})u_a(y)\right|_{y=P_i^b}
 =F_{ca}(P_i^b).
\end{aligned}
\]
In the second line, the center derivatives are taken before
substituting $y=P_i^b$. Thus the spatial derivative in
$F_{ca}(P_i^b)$ is supplied by the motion of the evaluation point.
Similarly, since $c<b$, $\B_c[\nabla W_b^\nu(y-P_i^b)]$ carries contributions from both differentiating in $y$ and differentiating in centers, and 
\[
 \B_c[\nabla W_b^\nu(y-P_i^b)]
 =D^2W_b^\nu(y-P_i^b)[u_c(y)-u_c(P_i^b)],
\]
where the Hessian is taken componentwise. 

Applying the product rule
in \eqref{eq:Fformula}, with the coefficients held fixed, yields
\begin{align}\label{eq:secondcoarse}
 \B_cF_{ab}(y)=A_b\sum_{i,\nu}\eta_i^{b,\nu}\big\{&
 [F_{ca}(y)-F_{ca}(P_i^b)]\cdot\nabla W_b^\nu(y-P_i^b)\notag\\
 &+[u_a(y)-u_a(P_i^b)]\otimes[u_c(y)-u_c(P_i^b)]
 :D^2W_b^\nu(y-P_i^b)\big\}.
\end{align}
Here, for vectors $v,w\in\R^d$, $[v\otimes w]:D^2W^\nu_b$ denotes the
vector with components $\sum_{p,q}v_pw_q\,\partial_p\partial_q(W^\nu_b)_m$,
$1\le m\le d$, that is, the Hessian is contracted componentwise in the
vector index of $W^\nu_b$.
Indeed, since $a<b$ and $c<b$, Remark~\ref{rem:Bderiv} gives
\[
\B_cF_{ab}
=(u_c(y)\cdot\nabla_y)F_{ab}
+\one_{\{c<a\}}\sum_k
(u_c(P_k^a)\cdot\nabla_{P_k^a})F_{ab}
+\sum_k(u_c(P_k^b)\cdot\nabla_{P_k^b})F_{ab}.
\]
Consider one summand of $F_{ab}$, suppressing the fixed coefficient:
\[
T_i^\nu(y)
=[u_a(y)-u_a(P_i^b)]\cdot\nabla W_b^\nu(y-P_i^b).
\]
Write $v_i=u_a(y)-u_a(P_i^b)$ and $s_i=y-P_i^b$.

The query derivative contributes
\[
(u_c(y)\cdot\nabla_y)T_i^\nu=[(u_c(y)\cdot\nabla)u_a(y)]
       \cdot\nabla W_b^\nu(s_i)\\
+[v_i\otimes u_c(y)]:D^2W_b^\nu(s_i).
\]
It does not differentiate $u_a(P_i^b)$, which is independent
of the query coordinate.

The scale-$a$ center derivatives, when $c<a$, contribute
\[
\one_{\{c<a\}}\sum_k
\Big[
(u_c(P_k^a)\cdot\nabla_{P_k^a})u_a(y)-\left.
(u_c(P_k^a)\cdot\nabla_{P_k^a})u_a(z)
\right|_{z=P_i^b}
\Big]\cdot\nabla W_b^\nu(s_i).
\]
They act on the centers defining $u_a$, but not on the wavelet,
whose center belongs to generation $b$. In the second term,
the center derivative is taken with $z$ fixed, before evaluation.

The scale-$b$ center derivatives contribute only through $P_i^b$:
\[
\begin{aligned}
(u_c(P_i^b)\cdot\nabla_{P_i^b})T_i^\nu
={}&-\left[
\left.(u_c(P_i^b)\cdot\nabla_z)u_a(z)\right|_{z=P_i^b}
\right]\cdot\nabla W_b^\nu(s_i)\\
&-[v_i\otimes u_c(P_i^b)]:D^2W_b^\nu(s_i).
\end{aligned}
\]
The first term differentiates the evaluation point in $u_a(P_i^b)$;
the second differentiates the translated wavelet.

Adding these three contributions, all terms multiplying
$\nabla W_b^\nu(s_i)$ combine into
\begin{multline*}
\Big[(u_c(y)\cdot\nabla)u_a(y)+\mathbf1_{\{c<a\}}\sum_k(u_c(P_k^a)\cdot\nabla_{P_k^a})u_a(y)\\
-\left.(u_c(P_i^b)\cdot\nabla_z)u_a(z)\right|_{z=P_i^b}-\mathbf1_{\{c<a\}}\sum_k\left.(u_c(P_k^a)\cdot\nabla_{P_k^a})u_a(z)\right|_{z=P_i^b}\Big]\cdot\nabla W_b^\nu(s_i)
\\
=[F_{ca}(y)-F_{ca}(P_i^b)]\cdot\nabla W_b^\nu(s_i),
\end{multline*}
by the definition of $F_{ca}$. The two Hessian terms combine into
\[
[v_i\otimes(u_c(y)-u_c(P_i^b))]:D^2W_b^\nu(s_i).
\]
Multiplying by $A_b\eta_i^{b,\nu}$ and summing gives
\eqref{eq:secondcoarse}. Thus the three coordinate-derivative
contributions in the remark are regrouped into the two
product-rule terms in \eqref{eq:secondcoarse}.

We next estimate the two terms in braces in \eqref{eq:secondcoarse} and their first spatial
 derivatives. First,
\[
 \|\nabla F_{ca}\|_\infty\le C\la_c\la_a.
\]
For $c<a$, this is \eqref{eq:Fcoarse} with $q=1$ and indices $(c,a)$.
For $c\ge a$, differentiate $F_{ca}=(u_c\cdot\nabla)u_a$ to obtain
\[
 \|\nabla F_{ca}\|_\infty
 \le C\bigl(\la_c\la_a+A_c\la_a/\ell_a\bigr)
 \le C\la_c\la_a,
\]
since $A_c=\ell_c\la_c$ and $\ell_c\le\ell_a$.
On the support of a scale-$b$ wavelet, put
\[
 v_{a,i}=u_a(y)-u_a(P_i^b),\quad
 v_{c,i}=u_c(y)-u_c(P_i^b),\quad
 E_i=F_{ca}(y)-F_{ca}(P_i^b).
\]
The support condition and the preceding gradient bound imply
\[
\begin{aligned}
 |E_i|&\le C\ell_b\la_a\la_c,
 &|\nabla E_i|&\le C\la_a\la_c,\\
 |v_{a,i}|&\le C\ell_b\la_a,
 &|\nabla v_{a,i}|&\le C\la_a,\\
 |v_{c,i}|&\le C\ell_b\la_c,
 &|\nabla v_{c,i}|&\le C\la_c.
\end{aligned}
\]
Since $|\nabla W_b^\nu|\le C\ell_b^{-1}$ and
$|D^2W_b^\nu|\le C\ell_b^{-2}$, each term in braces in
\eqref{eq:secondcoarse} is bounded by $C\la_a\la_c$.
One spatial derivative costs at most one additional factor
$\ell_b^{-1}$: it either differentiates $E_i$ or one of the $v$'s,
or falls on a wavelet derivative, using
$|\nabla D^2W_b^\nu|\le C\ell_b^{-3}$.
Bounded overlap and the coefficient bound therefore give
\[
 \|\B_cF_{ab}\|_\infty\le CA_b\la_a\la_c,
 \qquad
 \|\nabla(\B_cF_{ab})\|_\infty
 \le CA_b\la_a\la_c/\ell_b.
\]
By \eqref{eq:twopointbound},
\[
 |\delta(\B_cF_{ab})|
 \le CA_b\la_a\la_c\mathsf m_b
 =Cr\la_a\la_c\gamma_b
 =Cr\la_aK_{cb}.
\]

\medskip\noindent
\emph{Case 2: $a<b$ and $c\ge b$.}
The function $F_{ab}$ depends only on $y,P^a,P^b$, and $\B_c$
moves neither family of centers in this case. Thus
$\B_cF_{ab}=(u_c\cdot\nabla)F_{ab}$.
The $q=1,2$ bounds in \eqref{eq:Fcoarse} imply
\[
 \|\nabla F_{ab}\|_\infty\le C\la_a\la_b,
 \qquad
 |\delta\nabla F_{ab}|\le C\la_a\la_b\mathsf m_b.
\]
The product rule for differences and $\mathsf m_b\le\mathsf m_c$ now give
\[
\begin{aligned}
 |\delta(\B_cF_{ab})|
 &\le |\delta u_c|\,\|\nabla F_{ab}\|_\infty
      +\|u_c\|_\infty\,|\delta\nabla F_{ab}|\\
 &\le C A_c\la_a\la_b(\mathsf m_c+\mathsf m_b)
 \le Cr\la_a\la_b\gamma_c
 =Cr\la_aK_{cb}.
\end{aligned}
\]

\medskip\noindent
\emph{Cases 3 and 4: $a\ge b$.}
We must now prove the bound $Cr\la_bK_{ca}$.
The product rule applied to $F_{ab}=(u_a\cdot\nabla)u_b$ gives
\begin{equation}\label{eq:secondfine}
 \B_cF_{ab}=F_{ca}\cdot\nabla u_b
                  +u_a\cdot(\B_c\nabla u_b).
\end{equation}
Here $\B_c\nabla u_b$ means that $\B_c$ acts directly on the
configuration-dependent spatial derivative $\nabla_yu_b$;
we do not interchange $\B_c$ and $\nabla_y$.
Using \eqref{eq:material1}, \eqref{eq:Fsup}, and
\eqref{eq:spatialincrements}, we obtain
\[
\begin{aligned}
 |\delta(F_{ca}\cdot\nabla u_b)|
 &\le |\delta F_{ca}|\,\|\nabla u_b\|_\infty
       +\|F_{ca}\|_\infty\,|\delta\nabla u_b|\\
 &\le CrK_{ca}\la_b+C\la_{\min(c,a)} A_{\max(c,a)}\la_b\mathsf m_b
 \le Cr\la_bK_{ca}.
\end{aligned}
\]
For the last inequality, ${\max(c,a)}\ge a\ge b$ implies $\mathsf m_b\le\mathsf m_{\max(c,a)}$,
and 
$$\la_{\min(c,a)} A_{\max(c,a)}\mathsf m_{\max(c,a)}=rK_{ca}.$$ It remains to estimate the
second term in \eqref{eq:secondfine}.

\medskip\noindent
\emph{Case 3: $a\ge b$ and $c<b$.} Since $c<b$, the operator $\B_c$ differentiates both the query
point and the scale-$b$ centers. As $\nabla u_b$ depends only
on these coordinates, with the coefficients held fixed,
\[
\B_c\nabla u_b
=(u_c(y)\cdot\nabla_y)\nabla u_b
+\sum_i(u_c(P_i^b)\cdot\nabla_{P_i^b})\nabla u_b.
\]
Using
$\nabla u_b(y)=A_b\sum_{i,\nu}\eta_i^{b,\nu}
\nabla W_b^\nu(y-P_i^b)$ and
$\nabla_{P_i^b}=-\nabla_y$ on each translated wavelet, we obtain
\[
G(y):=\B_c\nabla u_b(y)
=A_b\sum_{i,\nu}\eta_i^{b,\nu}
[u_c(y)-u_c(P_i^b)]\cdot D^2W_b^\nu(y-P_i^b).
\]
On each support (given by the wavelet), the velocity difference is bounded by
$C\ell_b\la_c$ and its spatial gradient by $C\la_c$.
Differentiating the displayed formula once and using bounded overlap
therefore gives
\[
 \|G\|_\infty\le C\la_b\la_c,
 \qquad \|\nabla G\|_\infty\le C\la_b\la_c/\ell_b,
 \qquad |\delta G|\le C\la_b\la_c\mathsf m_b.
\]
Since $\mathsf m_b\le\mathsf m_a$, the product difference satisfies
\[
\begin{aligned}
 |\delta(u_a\cdot G)|
 &\le |\delta u_a|\,\|G\|_\infty
       +\|u_a\|_\infty\,|\delta G|\\
 &\le C A_a\la_b\la_c(\mathsf m_a+\mathsf m_b)
 \le Cr\la_b\la_c\gamma_a
 =Cr\la_bK_{ca}.
\end{aligned}
\]

\medskip\noindent
\emph{Case 4: $a\ge b$ and $c\ge b$.}
Now $\B_c$ does not move the scale-$b$ centers, so
$\B_c\nabla u_b=(u_c\cdot\nabla)\nabla u_b$.
The second term in \eqref{eq:secondfine} is consequently
$(u_a\otimes u_c):D^2u_b$. Taking a difference in each of its
three factors and using \eqref{eq:spatialincrements} gives
\[
\begin{aligned}
 |\delta((u_a\otimes u_c):D^2u_b)|
 &\le C\frac{\la_b}{\ell_b}A_aA_c
             (\mathsf m_a+\mathsf m_c+\mathsf m_b)\\
 &\le C\frac{\la_b}{\ell_b}A_{\min(c,a)}A_{\max(c,a)}\mathsf m_{\max(c,a)}
 =Cr\la_b\frac{A_{\min(c,a)}}{\ell_b}\gamma_{\max(c,a)}
\\& \le Cr\la_b\la_{\min(c,a)}\gamma_{\max(c,a)}
 =Cr\la_bK_{ca}.
\end{aligned}
\]
Indeed, all three increment factors are at most
$\mathsf m_{\max(c,a)}$, and 
$$A_{\min(c,a)}/\ell_b=\la_{\min(c,a)}\ell_{\min(c,a)}/\ell_b\le\la_{\min(c,a)}.$$
Together with the estimate for the first term of
\eqref{eq:secondfine}, this completes both cases with $a\ge b$,
and hence proves \eqref{eq:material2} for every ordering of the indices.
\end{proof}

\begin{remark}[Configurations attained by the process are admissible]
\label{rem:attained}
Let $\kappa\le\kappa_*$ with $\kappa_*$ from Lemma~\ref{lem:regularity},
and fix a realization of the random data. For $t\in[0,1)$
and $j_0\le j\le J$, let $k_j(t)$ be the index of the scale-$j$ window
containing $t$, and set
\[
 \Theta_t=\Big(\big(P_{i,k_j(t)}^j(t)\big)_{i\in I_j},\
               \big(\eta_{i,k_j(t)}^{j,\nu}\big)_{i,\nu}\Big)_{j_0\le j\le J},
 \qquad \mathcal S_t=(x_t,\widetilde x_t,\Theta_t).
\]
At $t=1$, set $k_j(1)=K_j-1$ and interpret all active
center positions by their continuous extensions from the left.
Thus $\Theta_1=\Theta_{1-}$ and $\mathcal S_1=\mathcal S_{1-}$.
Then (C1) holds for $\Theta_t$ at every $t$ by \eqref{eq:coefficients},
(C2) holds at every $t$ by Lemma~\ref{lem:regularity},
and distinct particles remain distinct by the same lemma. Comparing
\eqref{eq:frozenfield} with \eqref{eq:velocity} shows
$u_j(\cdot;\Theta_t)=u_j(\cdot,t)$. Hence $\mathcal S_t$ is
admissible at every $t$ at which $R_t\le1/4$, and
Lemmas~\ref{lem:material} and~\ref{lem:operator} apply to it. At a renewal time both the
configuration immediately before it (defined by left limits of the
window data) and the configuration immediately after it are
admissible. This is the precise content of the phrase ``at every state
attained by the process'' used below.
\end{remark}

We now apply the operators of Definition~\ref{def:operators} to the test
functions of Section~\ref{sec:testfunctions}. Recall that on the set of configurations with
$0<\dist(x,\widetilde x)\le1/4$ the shortest representative
$z=z(x,\widetilde x)\in\R^d$ of $x-\widetilde x$ is unique, satisfies
$|z|=\dist(x,\widetilde x)$, and depends smoothly on $(x,\widetilde x)$,
with $\nabla_xz=I$ and $\nabla_{\widetilde x}z=-I$. For $J\ge j_0$ and
$L\ge2$, the functions $\Phi$ of \eqref{eq:Hphi} and $\Phi_j$ of
\eqref{eq:cutoffs}, built with $\varepsilon=\rho\ell_J$, are therefore
smooth functions
\[
 (x,\widetilde x,P)\longmapsto\Phi\big(z(x,\widetilde x)\big),\qquad
 (x,\widetilde x,P)\longmapsto\Phi_j\big(z(x,\widetilde x)\big)
\]
on that set, which depend neither on the query point nor on the centers
or coefficients. (Recall that $\rho>0$ is fixed independently of $J,\kappa,L$ so that
$2\rho<\min(v_*,1/4)$, where $v_*$ is the threshold in
\eqref{eq:strain}, and that $\varepsilon=\rho\ell_J$.) Consequently, for all $a,b,c\in\{j_0,\ldots,J\}$, the
iterated derivatives $\B_b\Phi_b$, $\B_a\B_b\Phi_b$, and
$\B_c\B_a\B_b\Phi_b$ are well-defined continuous functions of the
configuration, and the chain rule together with
\eqref{eq:Bdefinition} and \eqref{eq:twopoint} gives the exact
identities
\begin{equation}\label{eq:Bchain}
 \B_b\Phi_b=\nabla\Phi_b(z)\cdot\delta u_b,\qquad
 \B_a\B_b\Phi_b=D^2\Phi_b(z)[\delta u_a,\delta u_b]
               +\nabla\Phi_b(z)\cdot\delta(\B_au_b).
\end{equation}
Indeed, since $z=z(x,\widetilde x)$ is independent of the query point and
the centers, with $\nabla_xz=I$ and $\nabla_{\widetilde x}z=-I$,
the definition of $\B_a$ gives $\B_a z=\delta u_a$. Hence
\[
\B_b\Phi_b=\nabla\Phi_b(z)\cdot\B_b z=\nabla\Phi_b(z)\cdot\delta u_b.
\]
Next, apply $\mathcal B_a$ to this product:
\[
\mathcal B_a\mathcal B_b\Phi_b
={}\bigl(\mathcal B_a[\nabla\Phi_b(z)]\bigr)\cdot\delta u_b+\nabla\Phi_b(z)\cdot\mathcal B_a[\delta u_b].
\]
The chain rule handles the first term:
\[
\mathcal B_a[\nabla\Phi_b(z)]
=D^2\Phi_b(z)[\delta u_a],
\]
so its contribution is $D^2\Phi_b(z)[\delta u_a,\delta u_b]$.

For the second term, we must justify
\[
\mathcal B_a[\delta u_b]=\delta(\mathcal B_a u_b).
\]
Recall that $\delta u_b=u_b(x;P^b)-u_b(\widetilde x;P^b)$.
Differentiating the first evaluation gives
\[
\begin{aligned}
\mathcal B_a[u_b(x;P^b)]
={}&\left.(u_a(x)\cdot\nabla_y)u_b(y;P^b)\right|_{y=x}\\
&+\mathbf1_{\{a<b\}}\sum_i
\left.(u_a(P_i^b)\cdot\nabla_{P_i^b})
u_b(y;P^b)\right|_{y=x}.
\end{aligned}
\]
The first term accounts for the motion of $x$; the second accounts
for the motion of the centers defining $u_b$. Together they are
exactly
\[
(\mathcal B_a u_b)(x;x,\widetilde x;P).
\]
The same calculation at $\widetilde x$, followed by subtraction,
gives the claimed identity by the definition \eqref{eq:twopoint}
of $\delta$.

Similarly, for the third order, the five-term expansion
\begin{equation}\label{eq:Bchain3}
\begin{aligned}
 \B_c\B_a\B_b\Phi_b
 ={}&D^3\Phi_b(z)[\delta u_c,\delta u_a,\delta u_b]
   +D^2\Phi_b(z)[\delta(\B_cu_a),\delta u_b]\\
   &+D^2\Phi_b(z)[\delta u_a,\delta(\B_cu_b)]
   +D^2\Phi_b(z)[\delta u_c,\delta(\B_au_b)]\\
   &+\nabla\Phi_b(z)\cdot\delta(\B_c\B_au_b),
\end{aligned}
\end{equation}
where every $\delta(\cdot)$ is the two-point difference
\eqref{eq:twopoint} of the corresponding function on $\mathcal E$,
evaluated at the given configuration. To see this, apply $\B_c$ to the second identity \eqref{eq:Bchain} get the third-order
expansion. We use the following facts. First, locally, $z=x-\widetilde x+n$, where the integer vector $n$ is
constant. Thus
\[
(\B_c z)_i
=\sum_j u_{c,j}(x)\,\partial_{x_j}z_i
+\sum_j u_{c,j}(\widetilde x)\,
  \partial_{\widetilde x_j}z_i=u_{c,i}(x)-u_{c,i}(\widetilde x),
\]
so that 
\[
\B_c z=u_c(x)-u_c(\widetilde x)=\delta u_c.
\]
It is the change in separation generated by moving both particles
with the scale-$c$ velocity field.

Secondly,
\[
\B_c[\delta G]=\delta(\B_cG),
\qquad G=u_a,\ u_b,\ \B_a u_b.
\]
The latter identity follows from the same evaluation chain rule
used above: $\B_c$ transports each particle with the field assigned
to the query point, and the center derivatives act on both sides
in the same way.

The product and chain rules give three terms when applied to the
Hessian term:
\[
\begin{aligned}
\B_c\bigl(D^2\Phi_b(z)[\delta u_a,\delta u_b]\bigr)
={}&D^3\Phi_b(z)[\delta u_c,\delta u_a,\delta u_b]\\
&+D^2\Phi_b(z)[\delta(\B_cu_a),\delta u_b]\\
&+D^2\Phi_b(z)[\delta u_a,\delta(\B_cu_b)].
\end{aligned}
\]
These arise by differentiating the Hessian and its two vector
arguments, respectively. The gradient term contributes two more:
\[
\begin{aligned}
\B_c\bigl(\nabla\Phi_b(z)\cdot\delta(\B_a u_b)\bigr)
={}&D^2\Phi_b(z)[\delta u_c,\delta(\B_a u_b)]\\
&+\nabla\Phi_b(z)\cdot\delta(\B_c\B_a u_b).
\end{aligned}
\]
Adding yields \eqref{eq:Bchain3}. No commutation of the scale operators is used; their order
is preserved throughout.

\begin{lemma}[Operator bounds for the test functions]\label{lem:operator}
There is a constant $C>0$, depending only on $d,\alpha,\psi,R,M$ and on
the fixed auxiliary functions $\vartheta,\chi$ of Section~\ref{sec:testfunctions}, with the
following property. Let $J\ge j_0$ and $L\ge2$, and let $\Phi$ and $\Phi_j$
($j_0\le j\le J$) be the test functions defined in
\eqref{eq:Hphi} and \eqref{eq:cutoffs}, respectively,
with $\varepsilon=\rho\ell_J$. Let
$\mathcal S=(x,\widetilde x,\Theta)$ be any admissible configuration in
the sense of Definition~\ref{def:admissible}, let
$r=\dist(x,\widetilde x)$, $S=S(r)=\max(r,\ell_J)$, and let $\gamma_j$,
$K_{ab}$ be given by \eqref{eq:rate-defs} with this $r$. Then, with
all left-hand sides evaluated at $\mathcal S$, for all
$a,b,c\in\{j_0,\ldots,J\}$,
\begin{align}
 |\B_b\Phi_b|&\le CS^\beta \gamma_b,\label{eq:B1}\\
 |\B_a\B_b\Phi_b|&\le CS^\beta K_{ab},\label{eq:B2}\\
 |\B_c\B_a\B_b\Phi_b|&\le CS^\beta
 \big(\gamma_cK_{ab}+\gamma_aK_{cb}+\gamma_bK_{ca}
       +\la_{\min(a,b)}K_{c,\max(a,b)}\big).
 \label{eq:B3}
\end{align}
Moreover:
\begin{enumerate}[(i)]
\item If $\ell_b>3LS(r)$, then
      $\B_b\Phi_b$, $\B_a\B_b\Phi_b$, and $\B_c\B_a\B_b\Phi_b$ all
      vanish at $\mathcal S$, for every $a,c\in\{j_0,\ldots,J\}$.
\item For every $b\in\{j_0,\ldots,J\}$,
      $|\B_b\Phi|\le CS^\beta \gamma_b$ at $\mathcal S$.
\end{enumerate}
The constant $C$ does not depend on $J$, $\kappa$, $L$, the
configuration $\mathcal S$, or the indices $a,b,c$.
\end{lemma}

\begin{proof}
Insert the bounds of Lemma~\ref{lem:material} and the derivative
bounds \eqref{eq:phijderivatives} into the identities
\eqref{eq:Bchain}--\eqref{eq:Bchain3}. The first identity in
\eqref{eq:Bchain} gives
$|\B_b\Phi_b|\le CS^\beta r^{-1}\cdot Cr\gamma_b$, which is \eqref{eq:B1};
the second gives
$|\B_a\B_b\Phi_b|\le CS^\beta(\gamma_a\gamma_b+K_{ab})\le CS^\beta K_{ab}$,
which is \eqref{eq:B2}, because $\gamma_{\min(a,b)}\le\la_{\min(a,b)}$
implies $\gamma_a\gamma_b\le K_{ab}$. In \eqref{eq:Bchain3} the five terms are
bounded, in the displayed order, by $CS^\beta$ times
$\gamma_c\gamma_a\gamma_b$, $K_{ca}\gamma_b$, $\gamma_aK_{cb}$, $\gamma_cK_{ab}$, and
$\la_{\min(a,b)}K_{c,\max(a,b)}$; absorbing the first into the fourth
by $\gamma_a\gamma_b\le K_{ab}$ proves \eqref{eq:B3}. All constants come from
Lemma~\ref{lem:material} and \eqref{eq:phijderivatives} and therefore
depend only on the quantities listed in the statement.

For (i), if $\ell_b>3LS(r)$, then $\nabla^k\Phi_b(z)=0$ for $k=1,2,3$ by
\eqref{eq:phijderivatives}, so every term in
\eqref{eq:Bchain}--\eqref{eq:Bchain3} vanishes at $\mathcal S$,
whatever $a$ and $c$ are. For (ii), the first identity in
\eqref{eq:Bchain} holds verbatim with $\Phi$ in place of $\Phi_b$,
and \eqref{eq:phiderivatives} with $k=1$ together with
\eqref{eq:material0} gives $|\B_b\Phi|\le CS^\beta \gamma_b$.
\end{proof}

\section{Summation over scales and the corrected observable}\label{sec:corrector}

This section performs the summation over scales and introduces the
corrected observable. With $\gamma_j=\min(\la_j,A_j/r)$, the scales coarser
than $r$ contribute $r\la_j=r\ell_j^{-\beta}$ and the finer ones
contribute $A_j=\ell_j^\alpha$; the two geometric series meet at
$\ell_j\approx r$ and sum to a multiple of $r^\alpha$, as in the
Littlewood--Paley characterization of $C^\alpha$. The corrected
observable is obtained by adding to $\Phi$ lower-order terms chosen
so that its time derivative between refreshes cancels exactly to
second order.

The objects of this section are functions of an admissible
configuration $\mathcal S$ (Definition~\ref{def:admissible}) and
a vector of phases
\[
h=(h_j)_{j_0\le j\le J}\in\prod_{j=j_0}^J[0,\tau_j].
\]
Along the process,
\[
h_j(t)=(k+1)\tau_j-t,
\qquad k\tau_j\le t<(k+1)\tau_j,
\]
is the time remaining until the nominal end of the active window,
even if that endpoint exceeds $1$. At $t=1$, use the left-limit convention
$h_j(1)=h_j(1-)=K_j\tau_j-1$. Thus $\dot h_j=-1$ between
renewals, and $h_j$ jumps from $0$ to $\tau_j$ at each renewal
of scale $j$.

The correction $\sum_bh_b\B_b\Phi_b$ cancels the first-order
evolution up to the cutoff error. Between renewals, the configuration evolves with generator
$\B=\sum_a\B_a$, while $\dot h_b=-1$. Hence
\[
\frac{d}{dt}\Phi=\sum_b\B_b\Phi,
\qquad
\frac{d}{dt}(h_b\B_b\Phi_b)
=-\B_b\Phi_b+h_b\sum_a\B_a\B_b\Phi_b.
\]
Summing gives
\[
\frac{d}{dt}\left(\Phi+\sum_bh_b\B_b\Phi_b\right)
=\sum_b\B_b(\Phi-\Phi_b)
+\sum_{a,b}h_b\B_a\B_b\Phi_b.
\]
We will discuss this more in Proposition~\ref{prop:drift}.
The first sum is the cutoff error; the second consists of the
second-order terms to be canceled by the weights $q_{ab}$.
The choice of remaining time also makes the old scale-$b$
first-order correction vanish just before that scale is renewed.

To cancel the remaining second-order terms, we introduce weights
$q_{ab}$ whose derivatives along the phase flow equal $-h_b$.
We also require $q_{ab}$ to vanish whenever either of the two
windows ends, so that the corresponding old correction disappears
at a renewal. Both requirements are met by integrating the remaining
lifetime of window $b$ up to the first of the two window endpoints:
\begin{equation}\label{eq:qdefinition}
 q_{ab}
 =\int_0^{\min(h_a,h_b)}(h_b-s)\dd s
 =h_b\min(h_a,h_b)-\tfrac12\min(h_a,h_b)^2.
\end{equation}
Indeed,
\[
q_{ab}=
\begin{cases}
h_ah_b-\tfrac12h_a^2,&h_a\le h_b,\\[2pt]
\tfrac12h_b^2,&h_a\ge h_b.
\end{cases}
\]
The values and first derivatives of these two expressions agree
on $h_a=h_b$, so $q_{ab}$ is continuously differentiable on the
nonnegative quadrant.

Write
\begin{equation}\label{eq:phaseflow}
 \mathcal D_h :=-\sum_{j=j_0}^J\frac{\partial}{\partial h_j}
\end{equation}
for differentiation along the phase flow $\dot h_j=-1$,
with the configuration held fixed. Thus $\mathcal D_h$ acts only
on the phase variables, whereas $\B$ acts on the configuration
variables with the phases held fixed. Along the process, between
renewals, the full time derivative is represented by
$\mathcal D_h+\B$. The preceding formulas give
\begin{equation}\label{eq:qproperties}
 \mathcal D_h q_{ab}=-h_b,\quad
 q_{bb}=\frac{h_b^2}{2},\quad
 q_{ab}=0\text{ if }h_ah_b=0,\quad
 0\le q_{ab}\le
 \frac{\kappa^2}{\la_b\max(\la_a,\la_b)}.
\end{equation}
For $a\ne b$, the derivative identity follows by applying
$-\partial_{h_a}-\partial_{h_b}$ to either branch; for $a=b$,
it follows directly from $q_{bb}=h_b^2/2$. Thus, along the
process and away from renewals, $\dot q_{ab}=-h_b$.
The upper bound follows from
\[
q_{ab}\le h_b\min(h_a,h_b)
\le\tau_b\min(\tau_a,\tau_b)
\le\frac{\kappa^2}{\la_b\max(\la_a,\la_b)}.
\]

Consequently, differentiating the correction
$\sum_{a,b}q_{ab}\B_a\B_b\Phi_b$ between renewals produces
$-\sum_{a,b}h_b\B_a\B_b\Phi_b$, which cancels the second-order
terms above, together with the third-order remainder
$\sum_{a,b,c}q_{ab}\B_c\B_a\B_b\Phi_b$.

\begin{lemma}[Uniform summation]\label{lem:sums}
Let $N_L=1+\log L$. There is a constant $C$, independent of $J$,
$\kappa$, and $L\ge2$, such that for every admissible configuration
$\mathcal S$ with separation $r$ and $S=S(r)$, and every phase vector
$h\in\prod_j[0,\tau_j]$, with all sums over $\{j_0,\ldots,J\}$ and all
operator expressions evaluated at $\mathcal S$ as in
Lemma~\ref{lem:operator},
\begin{align}
 \left|\sum_bh_b\B_b\Phi_b\right|&\le C\kappa N_L S^\beta,
 \label{eq:corrector1bound}\\
 \left|\sum_{a,b}q_{ab}\B_a\B_b\Phi_b\right|
 &\le C\kappa^2N_L^2 S^\beta,\label{eq:corrector2bound}\\
 \left|\sum_{a,b,c}q_{ab}\B_c\B_a\B_b\Phi_b\right|
 &\le C\kappa^2N_L^2.\label{eq:cubicbound}
\end{align}
\end{lemma}

Informally, the cutoff restricts the outer sum to scales satisfying
$\ell_b\le 3LS$. Of these, only $O(N_L)$ satisfy $\ell_b\ge S$.
There may still be arbitrarily many finer scales as $J$ increases,
but their contributions decay geometrically and have a sum bounded
independently of $J$. The factors $N_L$ and $N_L^2$ below arise
from summing over the surviving coarse scales, rather than from
bounding the total number of scales.

\begin{proof}
Set
\[
 p_b=\gamma_b/\la_b=\min(1,\ell_b/S),\quad
 \check\gamma_b=\sum_{a\ge b}\gamma_a,\quad
 \check p_b=\sum_{a\ge b}p_a,\quad
 \check d_b=\sum_{a\ge b}\check\gamma_a.
\]
All sums here and below are restricted to the available generations.
Since every available scale satisfies $\ell_a\ge\ell_J$, the definition
$S=\max(r,\ell_J)$ gives
\[
 \gamma_a=
 \begin{cases}
 \ell_a^{-\beta},&\ell_a\ge S,\\
 \ell_a^\alpha/S,&\ell_a<S.
 \end{cases}
\]
Indeed, when $r<\ell_J$, only the first case occurs. Since consecutive
scales differ by a factor of two, summing from the scale nearest $S$
towards the coarser or finer scales, respectively, gives
\[
 \sum_{\ell_a\ge S}\gamma_a
 \le S^{-\beta}\sum_{n\ge0}2^{-n\beta}
 =\frac{S^{-\beta}}{1-2^{-\beta}},
 \qquad
 \sum_{\ell_a<S}\gamma_a
 \le S^{-\beta}\sum_{n\ge0}2^{-n\alpha}
 =\frac{S^{-\beta}}{1-2^{-\alpha}}.
\]
Here and below, sums over scales are restricted to the available
generations; extending a geometric sum only increases its value.
Also, $\gamma_a\le\la_a$ and
\[
 \sum_{a<b}\la_a
 =\la_b\sum_{n=1}^{b-j_0}2^{-n\beta}
 \le\frac{\la_b}{2^\beta-1}.
\]
Combining this bound with the bound for the full sum of the
$\gamma_a$ yields
\begin{align}
 \sum_a\gamma_a&\le CS^{-\beta},&
 \sum_{a<b}\gamma_a&\le C\min(\la_b,S^{-\beta}),&
 \sum_{a<b}\la_a&\le C\la_b.
 \label{eq:sum-basic}
\end{align}

For the tail $\check\gamma_b$, if $\ell_b\ge S$, the full-sum
estimate gives $\check\gamma_b\le CS^{-\beta}$. If $\ell_b<S$,
every $a\ge b$ is also below $S$, and hence
\[
 \check\gamma_b
 =\frac{\ell_b^\alpha}{S}
   \sum_{n=0}^{J-b}2^{-n\alpha}
 \le\frac{S^{-\beta}}{1-2^{-\alpha}}
       \left(\frac{\ell_b}{S}\right)^\alpha,
\]
where $\alpha+\beta=1$ was used in the last expression.
Finally, splitting the sum defining $K_{cb}$ at $c=b$ gives
\[
 \sum_cK_{cb}
 =\gamma_b\sum_{c<b}\la_c
   +\la_b\sum_{c\ge b}\gamma_c
 \le C\la_b\gamma_b+\la_b\check\gamma_b
 \le C\la_b\check\gamma_b,
\]
since $\gamma_b\le\check\gamma_b$. Thus
\begin{align}
 \check\gamma_b
 &\le CS^{-\beta}\min\bigl(1,(\ell_b/S)^\alpha\bigr),&
 \sum_cK_{cb}
 &\le C\la_b\check\gamma_b.
 \label{eq:sum-K}
\end{align}

 We next estimate the tails $\check p_b$ and $\check d_b$.
Suppose first that $\ell_b\ge S$. The number of indices $a\ge b$
with $\ell_a\ge S$ is at most
\[
 1+\left\lfloor\log_2(\ell_b/S)\right\rfloor.
\]
For these indices, $p_a=1$ and
$\check\gamma_a\le CS^{-\beta}$ by \eqref{eq:sum-K}.
For the remaining indices, $\ell_a<S$, we have
\[
 p_a=\frac{\ell_a}{S},
 \qquad
 \check\gamma_a
 \le CS^{-\beta}\left(\frac{\ell_a}{S}\right)^\alpha.
\]
Starting at the first such index, the ratios $\ell_a/S$ are bounded
by $1,2^{-1},2^{-2},\ldots$. Consequently,
\[
 \sum_{\substack{a\ge b\\\ell_a<S}}p_a
 \le\sum_{n\ge0}2^{-n}=2,
 \qquad
 \sum_{\substack{a\ge b\\\ell_a<S}}\check\gamma_a
 \le CS^{-\beta}\sum_{n\ge0}2^{-n\alpha}
 \le CS^{-\beta}.
\]
Adding the contributions from scales above and below $S$ gives
\[
 \check p_b\le C\bigl(1+\log(\ell_b/S)\bigr),
 \qquad
 \check d_b\le CS^{-\beta}\bigl(1+\log(\ell_b/S)\bigr),
 \qquad \ell_b\ge S.
\]

If instead $\ell_b<S$, every index $a\ge b$ satisfies $\ell_a<S$.
Using $\ell_{b+n}=2^{-n}\ell_b$, we obtain
\[
 \check p_b
 =\frac{\ell_b}{S}\sum_{n=0}^{J-b}2^{-n}
 \le 2\frac{\ell_b}{S},
\]
and, by \eqref{eq:sum-K},
\[
 \check d_b
 =\sum_{a\ge b}\check\gamma_a
 \le CS^{-\beta}\left(\frac{\ell_b}{S}\right)^\alpha
       \sum_{n=0}^{J-b}2^{-n\alpha}
 \le CS^{-\beta}\left(\frac{\ell_b}{S}\right)^\alpha.
\]
Thus
\begin{equation}\label{eq:sum-ef}
 \check p_b\le C\begin{cases}
 1+\log(\ell_b/S),&\ell_b\ge S,\\
 \ell_b/S,&\ell_b<S,
 \end{cases}\quad
 \check d_b\le CS^{-\beta}\begin{cases}
 1+\log(\ell_b/S),&\ell_b\ge S,\\
 (\ell_b/S)^\alpha,&\ell_b<S.
 \end{cases}
\end{equation}
If $r<\ell_J$, then $S=\ell_J$ and every available scale satisfies
$\ell_a\ge S$; the sums over $\ell_a<S$ are empty, and the same
bounds apply.

For $\mathcal I=\{b:\ell_b\le3LS\}$, split the sums into
$S\le\ell_b\le3LS$ and $\ell_b<S$. There are at most
$1+\lfloor\log_2(3L)\rfloor\le CN_L$ indices in the first range.
For each such index,
\[
 p_b=1,\qquad
 \check p_b\le CN_L,\qquad
 \check\gamma_b\le CS^{-\beta},\qquad
 \check d_b\le CS^{-\beta}N_L,
\]
by \eqref{eq:sum-K} and \eqref{eq:sum-ef}. In the second range,
the same estimates give
\[
 p_b+\check p_b\le C\frac{\ell_b}{S},
 \qquad
 \check\gamma_b+\check d_b
 \le CS^{-\beta}\left(\frac{\ell_b}{S}\right)^\alpha.
\]
Summing these geometric tails contributes at most $C$ and
$CS^{-\beta}$, respectively. Combining the two ranges, and using
$N_L\ge1$, yields
\begin{equation}\label{eq:sum-final}
 \sum_{b\in\mathcal I}p_b\le CN_L,\quad
 \sum_{b\in\mathcal I}\check p_b\le CN_L^2,\quad
 \sum_{b\in\mathcal I}\check\gamma_b\le CS^{-\beta}N_L,\quad
 \sum_{b\in\mathcal I}\check d_b\le CS^{-\beta}N_L^2.
\end{equation}

By Lemma~\ref{lem:operator}(i), every summand in
\eqref{eq:corrector1bound}--\eqref{eq:cubicbound} with
$b\notin\mathcal I$ vanishes, so all sums over $b$ may be restricted to
$\mathcal I$. Using $h_b\le\kappa/\la_b$, \eqref{eq:B1}, and
\eqref{eq:sum-final}, we obtain
\[
 \left|\sum_bh_b\B_b\Phi_b\right|
 \le CS^\beta\sum_{b\in\mathcal I}h_b\gamma_b
 \le C\kappa S^\beta\sum_{b\in\mathcal I}
                         \frac{\gamma_b}{\la_b}
 =C\kappa S^\beta\sum_{b\in\mathcal I}p_b
 \le C\kappa N_LS^\beta.
\]
This proves \eqref{eq:corrector1bound}. Splitting into $a<b$ and $a\ge b$ in
\eqref{eq:qproperties} gives
\begin{equation}\label{eq:qKsum}
 \sum_a q_{ab}K_{ab}
 \le C\kappa^2(p_b+\check p_b)\le C\kappa^2\check p_b.
\end{equation}
Indeed the first part is bounded by
$\kappa^2\la_b^{-2}\gamma_b\sum_{a<b}\la_a$, and the second is
bounded by $\kappa^2\sum_{a\ge b}\gamma_a/\la_a$.
Together with \eqref{eq:B2} and the bound on
$\sum_{b\in\mathcal I}\check p_b$ in \eqref{eq:sum-final}, this proves
\eqref{eq:corrector2bound}.

To verify the cubic bound, fix $b$ and sum the four terms of
\eqref{eq:B3} separately. We use the consequence of
\eqref{eq:qproperties}
\[
 q_{ab}\le
 \begin{cases}
 \kappa^2/\la_b^2,&a<b,\\
 \kappa^2/(\la_b\la_a),&a\ge b.
 \end{cases}
\]
For the first term, \eqref{eq:sum-basic} and \eqref{eq:qKsum}
give directly
\[
 \sum_{a,c}q_{ab}\gamma_cK_{ab}
 =\left(\sum_c\gamma_c\right)
  \left(\sum_aq_{ab}K_{ab}\right)
 \le C\kappa^2S^{-\beta}\check p_b.
\]

For the second term, first sum in $c$ using \eqref{eq:sum-K},
and then split the $a$-sum:
\begin{align*}
 \sum_{a,c}q_{ab}\gamma_aK_{cb}
 &\le C\la_b\check\gamma_b\sum_aq_{ab}\gamma_a\\
 &\le C\kappa^2\check\gamma_b
 \left(
 \frac{1}{\la_b}\sum_{a<b}\gamma_a
 +\sum_{a\ge b}\frac{\gamma_a}{\la_a}
 \right)\le C\kappa^2\check\gamma_b(1+\check p_b).
\end{align*}
Here $\sum_{a<b}\gamma_a\le C\la_b$ by
\eqref{eq:sum-basic}, and the second sum is $\check p_b$.

For the third term, \eqref{eq:sum-K} gives
$\sum_cK_{ca}\le C\la_a\check\gamma_a$. Thus
\begin{align*}
 \sum_{a,c}q_{ab}\gamma_bK_{ca}
 &\le C\gamma_b\sum_aq_{ab}\la_a\check\gamma_a\\
 &\le C\kappa^2\gamma_b
 \left(
 \frac{1}{\la_b^2}\sum_{a<b}\la_a\check\gamma_a
 +\frac{1}{\la_b}\sum_{a\ge b}\check\gamma_a
 \right)\\
 &\le C\kappa^2
       \left(S^{-\beta}p_b+p_b\check d_b\right).
\end{align*}
In the last step we used
$\check\gamma_a\le CS^{-\beta}$,
$\sum_{a<b}\la_a\le C\la_b$,
$\gamma_b/\la_b=p_b$, and
$\sum_{a\ge b}\check\gamma_a=\check d_b$.

Finally, splitting the fourth term at $a=b$,
\begin{align*}
 &\sum_{a,c}q_{ab}\la_{\min(a,b)}K_{c,\max(a,b)}\\
 &\qquad=
 \sum_{a<b}q_{ab}\la_a\sum_cK_{cb}
 +\sum_{a\ge b}q_{ab}\la_b\sum_cK_{ca}\\
 &\qquad\le
 C\sum_{a<b}\frac{\kappa^2}{\la_b^2}
                 \la_a\la_b\check\gamma_b
 +C\sum_{a\ge b}\frac{\kappa^2}{\la_b\la_a}
                 \la_b\la_a\check\gamma_a\\
 &\qquad=
 C\kappa^2\left(
 \frac{\check\gamma_b}{\la_b}\sum_{a<b}\la_a
 +\sum_{a\ge b}\check\gamma_a
 \right)\\
 &\qquad\le C\kappa^2(\check\gamma_b+\check d_b).
\end{align*}
The first inequality uses $q_{ab}\le\kappa^2/\la_b^2$ when
$a<b$ and $q_{ab}\le\kappa^2/(\la_b\la_a)$ when $a\ge b$,
together with
$\sum_cK_{cb}\le C\la_b\check\gamma_b$ and
$\sum_cK_{ca}\le C\la_a\check\gamma_a$.
The last inequality follows from
$\sum_{a<b}\la_a\le C\la_b$ and the definition
$\check d_b=\sum_{a\ge b}\check\gamma_a$.

Sum in $b\in\mathcal I$, use $\check\gamma_b\le CS^{-\beta}$, $p_b\le1$,
and \eqref{eq:sum-final}. The result is
$C\kappa^2S^{-\beta}N_L^2$. Multiplication by $S^\beta$ in
\eqref{eq:B3} cancels $S^{-\beta}$, proving \eqref{eq:cubicbound}.
\end{proof}

The weights $h_b$ and $q_{ab}$ satisfy $\mathcal D_h h_b=-1$ and
$\mathcal D_h q_{ab}=-h_b$, which produces the cancellation of the
first- and second-order terms in Proposition~\ref{prop:drift}, and
$q_{ab}$ vanishes whenever one of the two windows has just ended,
which gives the jumps at a refresh the form of
Proposition~\ref{prop:renewal}.
Define the corrected observable, a function of an admissible
configuration $\mathcal S$ and a phase vector $h$,
\begin{equation}\label{eq:Q}
 Q(\mathcal S,h)=\Phi+\sum_bh_b\B_b\Phi_b
             +\sum_{a,b}q_{ab}\B_a\B_b\Phi_b,
\end{equation}
where $\Phi=\Phi(z(x,\widetilde x))$, the operator expressions are
evaluated at $\mathcal S$ as in Lemma~\ref{lem:operator}, and
$q_{ab}=q_{ab}(h)$. Along the process we write
$Q_t=Q(\mathcal S_t,h(t))$, with $\mathcal S_t$ the configuration of
Remark~\ref{rem:attained} and $h(t)$ the current phases; $Q_t$ is
defined at every $t$ at which $R_t\le1/4$, and $Q_{t-}$ denotes the
same expression with the left limits of the window data and phases. At $t=1$, set $Q_1=Q_{1-}$; no renewal is performed at $1$.

\begin{proposition}[Evolution between renewals]\label{prop:drift}
As functions of an admissible configuration $\mathcal S$ and a phase
vector $h$, with $\mathcal D_h$ the phase-flow derivative defined in
\eqref{eq:phaseflow} and $\B$ the operator of
Definition~\ref{def:operators},
\begin{equation}\label{eq:exactevolution}
 (\mathcal D_h+\B)Q=E_0+E_3,
 \quad E_0=\sum_j\B_j(\Phi-\Phi_j),\quad
 E_3=\sum_{a,b,c}q_{ab}\B_c\B_a\B_b\Phi_b.
\end{equation}
Moreover, writing $\mathcal S=(x,\widetilde x,\Theta)$ and
$r=\dist(x,\widetilde x)$,
\begin{align}
 |E_0|&\le CL^{-\beta},\quad |E_3|\le C\kappa^2N_L^2,
 \label{eq:errors}\\
 |Q-\Phi|&\le C(\kappa N_L+\kappa^2N_L^2)
                     (r^\beta+\ell_J^\beta).
 \label{eq:endpointbound}
\end{align}
Consequently, on every open interval $(s_1,s_2)$ containing no renewal
time of any scale and on which $R_t\le1/4$,
\begin{equation}\label{eq:pathwiseevolution}
 \frac{\mathrm d}{\mathrm dt}Q_t=(E_0+E_3)(\mathcal S_t,h(t)),
 \qquad s_1<t<s_2 .
\end{equation}
\end{proposition}

\begin{proof}
The identity \eqref{eq:exactevolution} is algebraic. The operator $\B$
acts on the configuration variables only and $\mathcal D_h$ on the
phases only; $\Phi$ and the operator expressions in \eqref{eq:Q} do
not depend on $h$, and $h_b$, $q_{ab}$ do not depend on $\mathcal S$.
Applying $\mathcal D_h$ to \eqref{eq:Q} with
$\mathcal D_h h_b=-1$ and $\mathcal D_h q_{ab}=-h_b$ from
\eqref{eq:qproperties}, and applying $\B=\sum_a\B_a$ to
\eqref{eq:Q}, gives
\[
-\sum_b\B_b\Phi_b -\sum_{a,b}h_b\B_a\B_b\Phi_b+ \B\Phi
 +\sum_{a,b}h_b\B_a\B_b\Phi_b
+E_3,
\]
which is \eqref{eq:exactevolution}. In particular all mixed
second-order terms cancel exactly.
Since $\Phi_j'=w_j\Phi'$, the first identity in \eqref{eq:Bchain}
(applied to $\Phi_j$ and to $\Phi$) gives $\B_j\Phi_j=w_j(r)\B_j\Phi$,
so the first-order error is $E_0=\sum_j(1-w_j)\B_j\Phi$.
Its nonzero summands have $\ell_j>L\widetilde S\ge LS\ge r$, so
$\gamma_j=\la_j$, and Lemma~\ref{lem:operator}(ii) gives
\[
 |E_0|\le CS^\beta\sum_{\ell_j>LS}\ell_j^{-\beta}
          \le CL^{-\beta}.
\]
The bound on $E_3$ is \eqref{eq:cubicbound}, and
\eqref{eq:endpointbound} follows from
\eqref{eq:corrector1bound}--\eqref{eq:corrector2bound} and
$S^\beta\le r^\beta+\ell_J^\beta$.

For \eqref{eq:pathwiseevolution}, fix such an interval. On this interval, the active coefficients are constant, and
$\mathcal S_t$ is admissible by Remark~\ref{rem:attained}.
Writing $P_i^m(t)$ for the active center of generation $m$, the
particle and center equations read
\[
 \dot x_t=\sum_a u_a(x_t;\Theta_t),\qquad
 \dot{\widetilde x}_t=\sum_a u_a(\widetilde x_t;\Theta_t),\qquad
 \dot P_i^m(t)=\sum_{a<m}u_a(P_i^m(t);\Theta_t).
\]
The fields depend on time through the moving centers. Between
renewals, however, the current particle positions, center positions,
and fixed active coefficients determine all velocities: there is
no additional dependence on the time variable. In other words,
the configuration equations have the form
$\dot{\mathcal S}_t=V(\mathcal S_t)$ for a fixed vector field $V$
on each interval between renewals.

For any $C^1$ function $F$ of the configuration, the chain rule
therefore gives
\begin{align*}
 \frac{\mathrm d}{\mathrm dt}F(\mathcal S_t)
 ={}&
 \sum_a u_a(x_t;\Theta_t)\cdot\nabla_xF(\mathcal S_t)
 +\sum_a u_a(\widetilde x_t;\Theta_t)
                \cdot\nabla_{\widetilde x}F(\mathcal S_t)\\
 &+\sum_m\sum_{i\in I_m}\sum_{a<m}
 u_a(P_i^m(t);\Theta_t)\cdot\nabla_{P_i^m}F(\mathcal S_t)\\
 ={}&\sum_a\Bigg[
 u_a(x_t;\Theta_t)\cdot\nabla_xF(\mathcal S_t)
 +u_a(\widetilde x_t;\Theta_t)
                \cdot\nabla_{\widetilde x}F(\mathcal S_t)\\
 &\hspace{35mm}
 +\sum_{m>a}\sum_{i\in I_m}
 u_a(P_i^m(t);\Theta_t)\cdot\nabla_{P_i^m}F(\mathcal S_t)
 \Bigg]\\
 ={}&\sum_a\B_aF(\mathcal S_t)
 =\B F(\mathcal S_t).
\end{align*}
The penultimate equality is precisely
Definition~\ref{def:operators}, applied to a function independent
of the query point. There are no contributions from derivatives
with respect to the coefficients, since those are constant on
this interval. The phases satisfy $\dot h_j=-1$.

Since $Q$ is $C^1$ in the configuration (the operator expressions in
\eqref{eq:Q} are at least $C^2$, because $\B_au_b\in C^2$ by
Definition~\ref{def:operators}) and $C^1$ in $h$, the chain rule gives
$\frac{\mathrm d}{\mathrm dt}Q(\mathcal S_t,h(t))
 =\big((\mathcal D_h+\B)Q\big)(\mathcal S_t,h(t))$,
and \eqref{eq:exactevolution} completes the proof.
\end{proof}

\subsection{Correctors for the quadratic-core observable}

Let $Q^\sharp$ be the observable \eqref{eq:Q} with
$(\Phi^\sharp,\Phi^\sharp_j)$ in place of $(\Phi,\Phi_j)$.
The common corrector estimates apply to these functions as well.

\begin{corollary}\label{cor:coretransfer}
Lemma~\ref{lem:operator}, Lemma~\ref{lem:sums} and
Proposition~\ref{prop:drift} hold with $(\Phi^\sharp,\Phi^\sharp_j,Q^\sharp)$
in place of $(\Phi,\Phi_j,Q)$, with constants depending in addition on
$\Upsilon$, at every admissible configuration with $0<r\le1/4$, for the
forward and reversed evolutions. All correction terms of
$Q^\sharp$, as well as $E_0$ and $E_3$, are continuous up to
$x=\widetilde x$ and vanish there; moreover, $Q^\sharp-\Phi^\sharp$ is
$C^2$ in $(x,\widetilde x)$ on $\{\dist(x,\widetilde x)<1/4\}$, the
diagonal included.
\end{corollary}

\begin{proof}
The proofs of these results use the test functions only through
\eqref{eq:phiderivatives} with $k=1$, \eqref{eq:phijderivatives} with
$k\le3$, the identity $\Phi_j'=w_j\Phi'$ and the support property (i)
of Lemma~\ref{lem:operator}. Since ${m_\sharp}\ge r$, \eqref{eq:corederivs}
implies these bounds, and $(\Phi^\sharp_j)'=w_j(\Phi^\sharp)'$ by
definition. The phases satisfy $0\le h_b\le\tau_b\le\kappa/\la_b$ in both
evolutions. Every term in \eqref{eq:Bchain}, \eqref{eq:Bchain3} and in
$E_0=\sum_j(1-w_j)\nabla\Phi^\sharp(z)\cdot\delta u_j$ is a product of
derivatives of $\Phi^\sharp_b$ or $\Phi^\sharp$ at $z$, which are
continuous by Lemma~\ref{lem:core}(b), and of at least one two-point
difference, which is continuous and vanishes when $x=\widetilde x$.
Finally, each correction term is a product of derivatives of
$\Phi^\sharp_b$ at $z$, which are smooth on $\{|z|<1/4\}$ by
Lemma~\ref{lem:core}(b), with values at $x$ or $\widetilde x$ of the
fields $u_b\in C^3$ and $\B_au_b\in C^2$ (Definition~\ref{def:operators});
hence $Q^\sharp-\Phi^\sharp$ is $C^2$ in $(x,\widetilde x)$ up to the
diagonal.
\end{proof}

\subsection{Diffusion acting on the correctors}

Write $Y=(x,\widetilde x)$ and
$\Delta_Y=\Delta_x+\Delta_{\widetilde x}$. The principal diffusive term, $2\mu\Delta\Phi^\sharp$, is nonnegative.
We must show that diffusion applied to the remaining terms in
$Q^\sharp$ is smaller than the gain of order $\kappa$. The following
estimate is uniform in the finest scale, including at the diagonal.

\begin{lemma}\label{lem:diffcorrector}
Assume \eqref{eq:eps0small}, $\varepsilon_0\le\kappa$ and $\kappa N_L\le1$.
There is $C_9$, independent of $J,L,\kappa,\varepsilon_0$, such that at
every configuration satisfying (C1), (C2) and $0\le r\le1/4$ (the pair
$Y=(y,y')$ playing the role of the particles $(x,\widetilde x)$),
\[
 \mu\,|\Delta_Y(Q^\sharp-\Phi^\sharp)|\le C_9\kappa^2N_L .
\]
\end{lemma}

\begin{proof}
We estimate the first and second correctors separately, sum over
the scales, and use the definition of $r_\sharp$ to control the result
near the diagonal.

\emph{Step 1: product formulas for the two-point Laplacian.} For a function $G$ of the query
point and a function $f$ of $z$,
\begin{gather*}
 \Delta_Y\,\delta G=\delta(\Delta G),\qquad
 \Delta_Yf(z)=2\Delta f(z),\qquad
 \nabla_{\widetilde x}f(z)=-\nabla_xf(z),\\
 \nabla_x\delta G=\nabla G(x),\qquad
 \nabla_{\widetilde x}\delta G=-\nabla G(\widetilde x).
\end{gather*}
The product rule gives, for a vector field $f(z)$ and a bilinear form
$\Lambda(z)$,
\begin{align*}
 \Delta_Y[f\cdot\delta G]={}&2\Delta f\cdot\delta G+f\cdot\delta(\Delta G)
 +2\textstyle\sum_k\partial_kf\cdot(\partial_kG(x)+\partial_kG(\widetilde x)),\\
 \Delta_Y\big(\Lambda[\delta G_1,\delta G_2]\big)={}&2\Delta \Lambda[\delta G_1,\delta G_2]
 +\Lambda[\delta\Delta G_1,\delta G_2]+\Lambda[\delta G_1,\delta\Delta G_2]\\
 &+2\textstyle\sum_k\partial_k\Lambda\big[\partial_kG_1(x)+\partial_kG_1(\widetilde x),\delta G_2\big]\\
 &+2\textstyle\sum_k\partial_k\Lambda\big[\delta G_1,\partial_kG_2(x)+\partial_kG_2(\widetilde x)\big]\\
 &+2\textstyle\sum_k\big(\Lambda(z)[\partial_kG_1(x),\partial_kG_2(x)]
 +\Lambda(z)[\partial_kG_1(\widetilde x),\partial_kG_2(\widetilde x)]\big).
\end{align*}

\emph{Step 2: derivative bounds for the velocity terms.} Write
$\mathsf m_j=\min(1,r/\ell_j)$, as in the
proof of Lemma~\ref{lem:material}, and let
$\ell_\vee=\min(\ell_a,\ell_b)$. From
\eqref{eq:frozenderivatives}, \eqref{eq:spatialincrements} and
\eqref{eq:material1}, we have
\[
\begin{aligned}
 \|\nabla u_a\|&\le C\la_a,&
 |\delta u_a|&\le Cr\gamma_a\le Cr\la_a,\\
 |\delta\Delta u_a|&\le C\la_a\ell_a^{-1}\mathsf m_a,&
 |\delta F_{ab}|&\le CrK_{ab}\le Cr\la_a\la_b.
\end{aligned}
\]
We also need the first two spatial derivatives of $F_{ab}$. If $a<b$,
\eqref{eq:Fcoarse} gives
\[
 \|\nabla F_{ab}\|\le C\la_a\la_b,\qquad
 \|D^2F_{ab}\|\le C\la_a\la_b\ell_b^{-1}.
\]
If $a\ge b$, differentiate $F_{ab}=(u_a\cdot\nabla)u_b$ and use
\eqref{eq:frozenderivatives} with $q\le3$. Since $\ell_a\le\ell_b$,
\[
\begin{aligned}
 \|\nabla F_{ab}\|
 &\le C(\la_a\la_b+A_a\la_b\ell_b^{-1})\le C\la_a\la_b,\\
 \|D^2F_{ab}\|
 &\le C(\la_a\ell_a^{-1}\la_b+\la_a\la_b\ell_b^{-1}
       +A_a\la_b\ell_b^{-2})
 \le C\la_a\la_b\ell_a^{-1}.
\end{aligned}
\]
Thus both cases give
\[
 \|D^2F_{ab}\|\le C\la_a\la_b\ell_\vee^{-1},\qquad
 |\delta\Delta F_{ab}|\le2\|\Delta F_{ab}\|
       \le C\la_a\la_b\ell_\vee^{-1}.
\]
The test functions satisfy \eqref{eq:corederivs}. We write
$\mathcal I=\{b:\ell_b\le3LS\}$; terms with $b\notin\mathcal I$ vanish
identically near the configuration.

\emph{Step 3: the first corrector.} Step 1 with $f=\nabla\Phi^\sharp_b$ and
$G=u_b$ gives
\[
\begin{aligned}
 |\Delta_Y(h_b\B_b\Phi^\sharp_b)|
 &\le Ch_b\big(S^\beta{m_\sharp}^{-3}r\la_b
 +S^\beta r{m_\sharp}^{-2}\la_b\ell_b^{-1}\mathsf m_b
 +S^\beta{m_\sharp}^{-2}\la_b\big)\\
 &\le C\kappa S^\beta{m_\sharp}^{-2}(1+r\ell_b^{-1}\mathsf m_b),
\end{aligned}
\]
 using $r\le{m_\sharp}$
and $h_b\la_b\le\kappa$.

\emph{Step 4: the second corrector.} Recall from \eqref{eq:Bchain} that
\[
 \B_a\B_b\Phi^\sharp_b
 =D^2\Phi^\sharp_b[\delta u_a,\delta u_b]
  +\nabla\Phi^\sharp_b\cdot\delta F_{ab}.
\]
For the first term, apply the second product formula in Step~1 with
$\Lambda=D^2\Phi^\sharp_b$. After factoring out
$CS^\beta{m_\sharp}^{-2}$, its six contributions are bounded by
\[
\begin{gathered}
 {m_\sharp}^{-2}r^2\la_a\la_b,\qquad
 \la_a\ell_a^{-1}\mathsf m_a\,r\la_b,\qquad
 r\la_a\,\la_b\ell_b^{-1}\mathsf m_b,\\
 {m_\sharp}^{-1}\la_a r\la_b,\qquad
 {m_\sharp}^{-1}r\la_a\la_b,\qquad \la_a\la_b.
\end{gathered}
\]
For the second term, apply the first product formula with
$f=\nabla\Phi^\sharp_b$ and $G=F_{ab}$. With the same factor removed,
the three contributions are bounded by
\[
 {m_\sharp}^{-1}r\la_a\la_b,\qquad
 r\la_a\la_b\ell_\vee^{-1},\qquad \la_a\la_b.
\]
Since $r\le m_\sharp$, $\ell_\vee\le\ell_a,\ell_b$, and
$\mathsf m_a,\mathsf m_b\le1$, all nine terms are bounded by
\[
 CS^\beta{m_\sharp}^{-2}\la_a\la_b(1+r/\ell_\vee).
\]
The coefficient is estimated by \eqref{eq:qproperties}:
\[
 q_{ab}\la_a\la_b\le\kappa^2\la_a/\max(\la_a,\la_b).
\]

\emph{Step 5: summation over the scales.} Put $\mathcal N:=N_L+\log(S/\ell_J)$. Since
$\ell_J\le\ell_b\le3LS$ for $b\in\mathcal I$,
$\#\mathcal I\le C\mathcal N$ and $\#\{a:a\ge b\}\le C\mathcal N$. Also
\[
 \sum_br\ell_b^{-1}\mathsf m_b\le\sum_{\ell_b\ge r}r^2\ell_b^{-2}
+\sum_{\ell_J\le\ell_b<r}r/\ell_b\le C(1+r/\ell_J).
\]
 For $b\in\mathcal I$,
by \eqref{eq:sum-basic},
\[
 \sum_a\frac{\la_a}{\max(\la_a,\la_b)}\Big(1+\frac r{\ell_\vee}\Big)
 \le\sum_{a<b}\frac{\la_a}{\la_b}\Big(1+\frac r{\ell_b}\Big)
 +\sum_{a\ge b}\Big(1+\frac r{\ell_a}\Big)
 \le C\Big(\mathcal N+\frac r{\ell_J}\Big).
\]
Hence
\[
\begin{aligned}
 \mu|\Delta_Y(Q^\sharp-\Phi^\sharp)|
 &\le C\mu S^\beta{m_\sharp}^{-2}
 \big[\kappa(\mathcal N+r/\ell_J)\\
 &\hspace{38mm}+\kappa^2\mathcal N(\mathcal N+r/\ell_J)\big].
\end{aligned}
\]

\emph{Step 6: comparison with the diffusivity.} We distinguish the
regions above and below the finest scale. If $r\ge\ell_J$, then $S={m_\sharp}=r$,
$\mathcal N=N_L+\log(r/\ell_J)$ with $r/\ell_J\ge1$, and
$\mu S^\beta{m_\sharp}^{-2}=\varepsilon_0(r/\ell_J)^{-(1+\alpha)}$; since
$\alpha>0$,
$(r/\ell_J)^{-(1+\alpha)}\big(1+\log(r/\ell_J)\big)^2(1+r/\ell_J)\le C$
for $r/\ell_J\ge1$, so the
bracket contributes at most $C\varepsilon_0(\kappa N_L+\kappa^2N_L^2)$. If
$r<\ell_J$, then $S=\ell_J$, $\mathcal N=N_L$, ${m_\sharp}\ge r_\sharp$, and
$\mu\ell_J^\beta r_\sharp^{-2}=\kappa/(4\beta\rho^\beta)$ by
\eqref{eq:diffusivity}, so the contribution is at most
$C(\kappa^2N_L+\kappa^3N_L^2)$. With $\varepsilon_0\le\kappa$ and
$\kappa N_L\le1$ both are at most $C\kappa^2N_L$. Steps 2--6 use
Lemma~\ref{lem:material}, which assumes $r>0$; the case $r=0$ follows by
continuity, since $Q^\sharp-\Phi^\sharp$ is $C^2$ in $Y$ up to the diagonal
(Corollary~\ref{cor:coretransfer}).
\end{proof}

\section{Renewal averages and time reversal}\label{sec:renewals}

\subsection{Averages over one renewal window}

The following identities are averages over the random data of one
window, that is, integrals over the shift $\xi\in Q_j$ and the finitely
many coefficients $\eta_i^\nu$. The role of the random lattice shift is
visible in the first display of the proof: averaging over $\xi$ turns
a lattice sum into an integral over $\T^d$, so that one scale of the
model has, on average, the translation-invariant structure function of
a single wavelet.

\begin{lemma}[Fresh-window identities]\label{lem:fresh}
For a freshly drawn scale-$j$ field
$$U(x)=A_j\sum_i\sum_{\nu\in\mathcal A}\eta_i^\nu W_j^\nu(x-y_i^j-\xi),$$ and fixed
$x,y\in\T^d$ with shortest difference $z=x-y$ satisfying $|z|\le1/4$,
\begin{align}
& \E U(x)=0, \qquad 
 \E[(U\cdot\nabla)U](x)=0,\label{eq:freshdrift}\\
& \E[(U(x)-U(y))\otimes (U(x)-U(y))]=\sigma^2A_j^2D(z/\ell_j).
 \label{eq:freshcov}
\end{align}
Here $\E$ is the average over $(\xi,\eta)$ only. The identities
therefore remain true when $x$ and $y$ are any points that do not
depend on $(\xi,\eta)$, for instance the particle positions at the
start of the window, which are functions of the data of earlier
windows only.
\end{lemma}

\begin{proof}
Centering proves the first identity. In the quadratic expressions,
independence removes terms with different spatial or mode indices.
For any periodic integrable function $F$,
\[
 \E_\xi\sum_{i\in I_j}F(y_i^j+\xi)
 =\ell_j^{-d}\int_{\T^d}F(\widetilde y)\dd\widetilde y.
\]
The one-point drift is therefore
$\sigma^2A_j^2\ell_j^{-1}\sum_\nu
\int_{\R^d}(W^\nu\cdot\nabla)W^\nu=0$, because each integrand
is $\nabla\cdot(W^\nu\otimes W^\nu)$ and has integral zero. 

The same unfolding gives \eqref{eq:freshcov}. Averaging first over the centered independent
coefficients and then over the lattice shift gives
\[
 \E[(U(x)-U(y))^{\otimes2}]
 =\sigma^2A_j^2\ell_j^{-d}\sum_\nu\int_{\T^d}
 [W_j^\nu(\widetilde x)-W_j^\nu(\widetilde x-z)]^{\otimes2}\dd\widetilde x.
\]
Set $\Psi^\nu(\widetilde x)=W^\nu(\widetilde x/\ell_j)
-W^\nu((\widetilde x-z)/\ell_j)$. Its support lies in two balls
whose union has diameter at most $2R\ell_j+|z|\le3/4<1$, so
distinct integer translates of its support are disjoint. Thus the
torus integral unfolds to $\int_{\R^d}(\Psi^\nu)^{\otimes2}$.
Changing variables $\widetilde x=\ell_j\widetilde y$ cancels
$\ell_j^{-d}$, and summing over $\nu$ gives
$\sigma^2A_j^2D(z/\ell_j)$, as claimed. 
\end{proof}

\begin{lemma}[Fresh-window identities after a measure-preserving
deformation]\label{lem:freshPsi}
Let $\Psi:\T^d\to\T^d$ be a Lebesgue-measure-preserving bijection, and
$U(x)=A_j\sum_{i,\nu}\eta_i^\nu W_j^\nu(x-\Psi(y_i^j+\xi))$ with $\xi$
uniform on $Q_j$ and $(\eta_i^\nu)$ as in \eqref{eq:coefficients}, all
independent. Then \eqref{eq:freshdrift} and \eqref{eq:freshcov} hold
for $U$.
\end{lemma}

\begin{proof}
The identity $\E U(x)=0$ follows directly from the centering of the
coefficients. For the quadratic identities, we first average over
the coefficients with the shift $\xi$ fixed. Independence and
centering eliminate all mixed terms, while each diagonal term has
coefficient $\sigma^2$. Thus, writing
\[
 F_\nu(p):=
 [W_j^\nu(x-p)-W_j^\nu(y-p)]^{\otimes2},
\]
we obtain
\[
 \E_\eta[(U(x)-U(y))^{\otimes2}]
 =\sigma^2A_j^2\sum_{i,\nu}
 F_\nu\bigl(\Psi(y_i^j+\xi)\bigr).
\]
The deformation $\Psi$ is still present after this first average.

We now average over the shift. For every periodic integrable
function $F$, the uniform distribution of $\xi$ on $Q_j$ gives
\[
\begin{aligned}
 \E_\xi\sum_i F\bigl(\Psi(y_i^j+\xi)\bigr)
 &=\ell_j^{-d}\sum_i\int_{Q_j}
       F\bigl(\Psi(y_i^j+\xi)\bigr)\,\dd\xi\\
 &=\ell_j^{-d}\int_{\T^d}F(\Psi(q))\,\dd q=\ell_j^{-d}\int_{\T^d}F(p)\,\dd p.
\end{aligned}
\]
Here the second equality uses the tiling of $\T^d$ by the cells
$y_i^j+Q_j$, and the last equality uses that $\Psi$ preserves
Lebesgue measure. Applying this identity componentwise to $F_\nu$
removes the deformation and yields
\[
 \E[(U(x)-U(y))^{\otimes2}]
 =\sigma^2A_j^2\ell_j^{-d}\sum_\nu
 \int_{\T^d}
 [W_j^\nu(x-p)-W_j^\nu(y-p)]^{\otimes2}\,\dd p.
\]
This is precisely the covariance integral in the proof of
Lemma~\ref{lem:fresh}; the same unfolding and rescaling give
$\sigma^2A_j^2D(z/\ell_j)$, where $z=x-y$ is the shortest
displacement. This proves \eqref{eq:freshcov}.

For the quadratic drift, the coefficient average similarly gives
\[
 \E_\eta[(U\cdot\nabla)U](x)
 =\sigma^2A_j^2\sum_{i,\nu}
 \bigl[(W_j^\nu\cdot\nabla)W_j^\nu\bigr]
 \bigl(x-\Psi(y_i^j+\xi)\bigr).
\]
The same shift average therefore implies
\[
 \E[(U\cdot\nabla)U](x)
 =\sigma^2A_j^2\ell_j^{-d}\sum_\nu
 \int_{\T^d}(W_j^\nu\cdot\nabla)W_j^\nu(q)\,\dd q=0.
\]
Indeed, each wavelet is divergence-free, so the integrand equals
$\nabla\cdot(W_j^\nu\otimes W_j^\nu)$, whose integral vanishes
by periodicity. This completes the proof of \eqref{eq:freshdrift}.
\end{proof}

\subsection{The reversed field and its renewal data}\label{sec:reversal}

The dissipation of the scalar up to time $t_0$ is expressed through an
equation driven by the time-reversed field
(Lemma~\ref{lem:pairidentity}). This section describes the reversed
field and verifies the two properties needed for the same corrector
and jump identities: the structure
between renewals, and the independence of the fresh data at a renewal
time from the past of the reversed field.

Fix $J\ge j_0$ and $t_0\in(0,1]$. By Lemma~\ref{lem:betweenrenewals}, for
each realization, $t\mapsto u_j(x,t)$ has left limits $u_j(x,t-)$ for
$t\in(0,1]$.

Forward renewal times are denoted by $\zeta$, and their reversed
counterparts by $\overleftarrow{\zeta}=t_0-\zeta$.
In statements applying to either evolution, $\tilde{\zeta}$ denotes
a renewal time in the time variable of the chosen evolution; left
limits and jumps are always taken in that time variable.

\begin{definition}[Reversed field]\label{def:reversed}
For $s\in[0,t_0)$ and $j_0\le j\le J$ let $k$ be such that
$t_0-s\in(k\tau_j,(k+1)\tau_j]$, and set
\[
 \overleftarrow u_j(x,s):=-u_j\big(x,(t_0-s)-\big),\qquad
 \overleftarrow u:=\sum_{j=j_0}^J\overleftarrow u_j,\qquad
 \overleftarrow P{}^j_i(s):=P^j_{i,k}\big((t_0-s)-\big),
\]
where $P^j_{i,k}((k+1)\tau_j-)$ is the left limit within window $k$; at
$s=t_0$ all three are defined by their left limits as $s\uparrow t_0$.
The \emph{reversed renewal times} of scale $j$ are the points
$\overleftarrow{\zeta}=t_0-\zeta$ with
$\zeta\in\tau_j\mathbb N\cap(0,t_0)$; the reversed windows of scale $j$
are the intervals between consecutive reversed renewal times, together
with the initial window starting at $0$. The \emph{reversed phase} is
$\overleftarrow h_j(s):=(t_0-s)-k\tau_j$ for $s<t_0$, and
$\overleftarrow h_j(t_0):=\overleftarrow h_j(t_0-)$.
\end{definition}

Changing the field on a finite set of times changes neither the
solutions of the center equations nor the energy solutions of
Section~\ref{sec:weak}, so the left-limit conventions are immaterial;
in particular $\overleftarrow u(x,s)=-u(x,t_0-s)$ for all but finitely
many $s$. By construction $0<\overleftarrow h_j(s)\le\tau_j$ for $s<t_0$
and $\overleftarrow h_j(t_0)=0$; $\overleftarrow h_j$ decreases with slope
$-1$ between reversed renewals; and at a reversed renewal time $\overleftarrow{\zeta}$ of
scale $j$, $\overleftarrow h_j(\overleftarrow{\zeta}-)=0$ and $\overleftarrow h_j(\overleftarrow{\zeta})=\tau_j$.
Every reversed window starting at a time $\overleftarrow{\zeta}>0$ has length exactly
$\tau_j$ (it corresponds to the forward window
$[t_0-\overleftarrow{\zeta}-\tau_j,t_0-\overleftarrow{\zeta})$), and the initial reversed window has length
at most $\tau_j$. Thus the reversed phases play exactly the role of the
phases of Section~\ref{sec:corrector}.

\begin{lemma}[Structure of the reversed field]\label{lem:reversed}
\begin{enumerate}[(a)]
\item On each reversed window of scale $j$,
 $\overleftarrow u_j(x,s)=A_j\sum_{i,\nu}(-\eta^{j,\nu}_{i,k})
 W_j^\nu(x-\overleftarrow P{}^j_i(s))$ with constant coefficients, and
 $\frac{d}{ds}\overleftarrow P{}^j_i
 =\overleftarrow u_{<j}(\overleftarrow P{}^j_i,s)$ for a.e.\ $s$.
\item For $\kappa\le\kappa_*$, the conclusions of
 Lemma~\ref{lem:regularity} (separation, \eqref{eq:derivativebounds},
 \eqref{eq:holder}) hold for $\overleftarrow u$ with the same constants,
 and the configurations $(x,\widetilde x,\overleftarrow\Theta_s)$ with
 $0<\dist(x,\widetilde x)\le1/4$ are admissible
 (Definition~\ref{def:admissible}). Between consecutive reversed
 renewal times of all scales, the reversed configuration evolves with
 the generator $\overleftarrow{\B}=\sum_a\overleftarrow{\B}_a$ of
 Definition~\ref{def:operators}, built from the frozen fields with
 centers $\overleftarrow P$ and coefficients $-\eta$: explicitly, the
 reversed frozen field of generation $j$ at time $s$ is the field
 \eqref{eq:frozenfield} formed with the centers
 $\overleftarrow P{}^j_i(s)$ and the coefficients $-\eta^{j,\nu}_{i,k}$ of
 the current reversed window, and it equals $\overleftarrow u_j(\cdot,s)$.
 Since $|-\eta|=|\eta|\le M$, condition (C1) holds for these
 coefficients.
\item (Freshness.) Let $\overleftarrow{\zeta}\in(0,t_0)$ be a reversed renewal time
 of at least one scale, put $\zeta:=t_0-\overleftarrow{\zeta}$, and let
 $j_1$ be the smallest index with
 $\zeta\in\tau_{j_1}\mathbb N$. The set of scales renewed at $\overleftarrow{\zeta}$ is
 $\mathcal J(\overleftarrow{\zeta})=\{j_1,\dots,J\}$. Let $\omega_{\overleftarrow{\zeta}}$ be the data
 $(\xi^m,\eta^m)$ of the forward windows $[\zeta-\tau_m,\zeta)$,
 $m\in\mathcal J(\overleftarrow{\zeta})$, and $\omega_{\overleftarrow{\zeta}}^c$ all other data. Then:
 \begin{enumerate}[(c1)]
 \item $\overleftarrow u$ restricted to $[0,\overleftarrow{\zeta})$ is a function of
 $\omega_{\overleftarrow{\zeta}}^c$;
 \item for $m\in\mathcal J(\overleftarrow{\zeta})$,
 $\overleftarrow P{}^m_i(\overleftarrow{\zeta})=\Psi_m(y_i^m+\xi^m)$, where $\Psi_m$ is the
 forward flow map of $u_{<m}$ from time $\zeta-\tau_m$ to $\zeta$; $\Psi_m$ is a
 Lebesgue-measure-preserving $C^1$ diffeomorphism of $\T^d$ and a
 function of $\omega_{\overleftarrow{\zeta}}^c$ and of
 $(\xi^a,\eta^a)_{a\in\mathcal J(\overleftarrow{\zeta}),\,a<m}$.
 \end{enumerate}
\end{enumerate}
\end{lemma}

\begin{proof}
(a) On a reversed window the forward time $t_0-s$ runs through a
single forward window $(k\tau_j,(k+1)\tau_j]$, where the coefficients
are $\eta_{i,k}^{j,\nu}$; for $s$ off the finite set of renewal times,
$\frac{d}{ds}P(t_0-s)=-u_{<j}(P(t_0-s),t_0-s)$.

(b) The fields $\overleftarrow u_j(\cdot,s)$ are, up to sign, the
forward fields at time $(t_0-s)-$, and the reversed centers are forward
center positions or their left limits. The separation of the centers,
\eqref{eq:derivativebounds} and \eqref{eq:holder} are statements about
the field and the centers at one time, and they pass to left limits
(left limits of $\ell_j/2$-separated configurations are
$\ell_j/2$-separated). The statement about the generator is (a)
together with the computation in the proof of
Proposition~\ref{prop:drift}.

(c) By Lemma~\ref{lem:nested}, $\zeta\in\tau_m\mathbb N$ exactly for
$m\ge j_1$. We claim that for $f\in[\zeta,t_0)$, $u_a(\cdot,f)$ is a
function of the data
\begin{multline*}
 \mathsf D:=\{\text{data of all windows of scales }a'<j_1\}\\
 \cup\{\text{data of the windows of scales }m\ge j_1\text{ starting at times}\ge \zeta\}.
\end{multline*}
We argue by induction on $a$. For $a<j_1$, $u_a$ depends only on data
of scales $\le a<j_1$. For $a\ge j_1$, the scale-$a$ window containing
$f$ starts at some $\zeta'\in[\zeta,f]$ (as $\zeta\in\tau_a\mathbb N$), and
$u_a(\cdot,f)$ is determined by the data of that window, which lie in
$\mathsf D$, and by its centers, which solve $\dot P=u_{<a}(P,\cdot)$ on
$[\zeta',f]\subset[\zeta,t_0)$ from $y_i+\xi$; by the induction hypothesis
$u_{<a}$ on $[\zeta,t_0)$ is a function of $\mathsf D$. Left limits at
$f\in(\zeta,t_0]$ are limits of values on $[\zeta,f)$. Hence $\overleftarrow u$
on $[0,\overleftarrow{\zeta})$ is a function of $\mathsf D$, and
$\mathsf D\cap\omega_{\overleftarrow{\zeta}}=\emptyset$ because the windows in $\omega_{\overleftarrow{\zeta}}$
start at $\zeta-\tau_m<\zeta$. This is (c1).

For (c2): $\overleftarrow P{}^m_i(\overleftarrow{\zeta})=P^m_{i,k}(\zeta-)$ with
$k=\zeta/\tau_m-1$, and $P^m_{i,k}(\zeta-)=\Psi_m(y_i^m+\xi^m)$ by
\eqref{eq:centers}. We claim that, for $j_1\le a\le m$, the
restriction of $u_a$ to $[\zeta-\tau_a,\zeta)$ is a function of
$\omega^c_{\overleftarrow{\zeta}}$ and of $(\xi^{a'},\eta^{a'})_{j_1\le a'\le a}$; since
$[\zeta-\tau_m,\zeta)\subset[\zeta-\tau_a,\zeta)$ for $a\le m$, it follows that
$u_{<m}$ on $[\zeta-\tau_m,\zeta)$ is a function of
$\omega^c_{\overleftarrow{\zeta}}$ and of the fresh data of the scales in $[j_1,m)$.
We argue by induction on $a$. A scale $a'<j_1$ is determined, at all
times, by the data of scales $<j_1$, which belong to
$\omega^c_{\overleftarrow{\zeta}}$. A scale $j_1\le a\le m$ is, on
$[\zeta-\tau_a,\zeta)$, in its window $[\zeta-\tau_a,\zeta)$, whose data
$(\xi^a,\eta^a)$ belong to $\omega_{\overleftarrow{\zeta}}$, and its centers are
transported on this interval by $u_{<a}$, which by the induction
hypothesis is there a function of $\omega^c_{\overleftarrow{\zeta}}$ and of
$(\xi^{a'},\eta^{a'})_{j_1\le a'<a}$ (the restriction of $u_{a'}$ to
$[\zeta-\tau_a,\zeta)\subset[\zeta-\tau_{a'},\zeta)$ being a function of the
restriction to the larger interval). This proves the claim. On each subinterval between points of
$\tau_J\mathbb N$ this field is $C^3$ in $x$, continuous in time and
divergence-free (Lemma~\ref{lem:betweenrenewals}), so its flow is a
$C^1$ diffeomorphism with Jacobian determinant
$\exp(\int\nabla\cdot u_{<m})=1$; composing finitely many such maps
gives (c2).
\end{proof}

\subsection{Jump identity in both directions of time}

The same jump calculation applies to both families of observables.
Let $(\Xi,\Xi_j)$ be $(\Phi,\Phi_j)$ for a distinct pair, or
$(\Phi^\sharp,\Phi^\sharp_j)$ for a pair that may lie on the diagonal.
Define $Q$ by \eqref{eq:Q} with $(\Xi,\Xi_j)$ in place of
$(\Phi,\Phi_j)$, and put
\[
 \Gamma^\Xi_j(r):=\frac{\sigma^2A_j^2}{2}D(z/\ell_j):D^2\Xi_j(z),\qquad|z|=r;
\]
For the quadratic-core functions, $\Gamma^\Xi_j(0)=0$ since $D(0)=0$.

\begin{proposition}[Averaged jump in both directions of time]\label{prop:jumpboth}
Let $\tilde{\zeta}\in(0,t_0)$ be a renewal time of the forward or reversed field, $\mathcal J(\tilde{\zeta})$ the set of renewed scales,
and $\omega_{\tilde{\zeta}}$ the fresh data (the data of the windows starting at
$\tilde{\zeta}$ in the forward case; as in Lemma~\ref{lem:reversed}(c) in the
reversed case). Let $x,\widetilde x\in\T^d$ be fixed points with
$\dist(x,\widetilde x)\le1/4$ (strictly positive for the logarithmic functions), and let $\mathcal S_{\tilde{\zeta}-},\mathcal S_{\tilde{\zeta}}$
be the configurations with particles $x,\widetilde x$ and the window
data before and after $\tilde{\zeta}$. Then
\[
 \E_{\tilde{\zeta}}\big[Q(\mathcal S_{\tilde{\zeta}},h(\tilde{\zeta}))-Q(\mathcal S_{\tilde{\zeta}-},h(\tilde{\zeta}-))\big]
 =\sum_{j\in\mathcal J(\tilde{\zeta})}\tau_j^2\Gamma^\Xi_j(\dist(x,\widetilde x)),
\]
where $\E_{\tilde{\zeta}}$ integrates over $\omega_{\tilde{\zeta}}$ with all other data fixed.
\end{proposition}

\begin{proof}
Write $\mathcal J=\mathcal J(\tilde{\zeta})=\{j_1,\dots,J\}$ and drop the arrows.
By \eqref{eq:Bchain} and \eqref{eq:Fformula}, each first-order
term uses only the data of its scale and each second-order term
only the data of its two scales. Terms with all indices outside
$\mathcal J$ are therefore unchanged. Old terms with an index in
$\mathcal J$ vanish because $h_j(\tilde{\zeta}-)=0$ and $q_{ab}=0$
when $h_ah_b=0$. In the reversed case the new centers of scale
$m\in\mathcal J$ are $\Psi_m(y^m+\xi^m)$ with $\Psi_m$ as in
Lemma~\ref{lem:reversed}(c2); in the forward case
$\Psi_m=\mathrm{id}$. In both cases $\Psi_m$, the particles, and all data
of scales outside $\mathcal J$ are independent of $(\xi^m,\eta^m)$ and of
the fresh data of the scales finer than $m$. We treat the new terms.

\emph{First order,} $b\in\mathcal J$: $\tau_b\nabla\Xi_b(z)\cdot\delta u_b$ is
linear in $\eta^b$, and every other quantity in it is independent of
$\eta^b$. Integrating first over $\eta^b$ gives zero.

\emph{Mixed,} $a\ne b$ with $\{a,b\}\cap\mathcal J\neq\emptyset$: then
$m:=\max(a,b)\in\mathcal J$, since $\mathcal J$ is an upper set. By
\eqref{eq:Bchain} and \eqref{eq:Fformula}, $\B_a\B_b\Xi_b$ is linear in
$\eta^m$ when all centers and the other coefficients are fixed, and the
centers of the scales $a,b$ do not depend on $\eta^m$ (those of scale
$m$ depend on $\xi^m$ and on $\Psi_m$; those of the coarser scale do
not involve scale $m$). Integrating first over $\eta^m$ gives zero.

\emph{Diagonal,} $a=b\in\mathcal J$: the new term is
$\frac{\tau_b^2}{2}\big(D^2\Xi_b(z)[\delta u_b,\delta u_b]
+\nabla\Xi_b(z)\cdot\delta((u_b\cdot\nabla)u_b)\big)$. Integrating over
$(\xi^b,\eta^b)$ with $\Psi_b$ and the particles fixed,
Lemma~\ref{lem:freshPsi} (the sign of the coefficients is irrelevant)
gives $\frac{\tau_b^2}{2}\sigma^2A_b^2D(z/\ell_b):D^2\Xi_b(z)
=\tau_b^2\Gamma^\Xi_b$.

Each term is bounded, so Fubini's theorem justifies integrating in the
stated order; the remaining fresh variables are then integrated
trivially.
\end{proof}

\subsection{Jumps along the forward trajectories}

Along the forward trajectories, we integrate the random data one window at a time. For a renewal time $s$, write the random data as $(\omega_{<s},\omega_s,\omega_{>s})$, the data of the
windows starting before, at, and after $s$. Every quantity attached
to the process on $[0,s)$ or at $s-$, in particular $Q_{s-}$, $R_s$,
and the configuration $\mathcal S_{s-}$, is a function of
$\omega_{<s}$ alone, because trajectories up to time $s$ are determined
by the data of the windows active before $s$; $Q_s$ depends in
addition on $\omega_s$. For every bounded measurable $G$ define the partial average
\begin{equation}\label{eq:partialaverage}
 \E_s[G]:=\int G\,(\omega_{<s},\omega_s,\omega_{>s})\dd\mathbb P_s(\omega_s),
\end{equation}
the integral over $\omega_s$ against its law $\mathbb P_s$, a finite product
of uniform measures on cubes and of coefficient laws, with all other
variables held fixed. By independence and Fubini's theorem,
\[
 \E[G]
 =\int\!\!\int
   (\E_sG)(\omega_{<s},\omega_{>s})\,
   \dd\mathbb P_{<s}(\omega_{<s})\,\dd\mathbb P_{>s}(\omega_{>s})
 =\E[\E_sG].
\]
In the last expression, $\E$ still denotes expectation over all
random data; integration over $\omega_s$ contributes a factor one
because $\E_sG$ does not depend on that coordinate.

Retain the definition of $\Gamma_j$ from \eqref{eq:Gamma}.

\begin{proposition}[Averaged jump at a refresh]\label{prop:renewal}
At a renewal time $s<1$ with $0<R_s\le1/4$, let $\mathcal J(s)$ be the
set of renewed generations. Then
\begin{equation}\label{eq:jump}
 \E_s[Q_s-Q_{s-}]
    =\sum_{j\in\mathcal J(s)}\tau_j^2\Gamma_j(R_s),
\end{equation}
an identity between functions of $\omega_{<s}$.
At $s=0$ use $Q_{0-}=\Phi(r_0)$ and regard all initial windows as
renewals; then both sides are deterministic.
\end{proposition}

\begin{proof}
For $s>0$, the particle positions are functions of $\omega_{<s}$
and are independent of the fresh data. Proposition~\ref{prop:jumpboth}
in the forward direction, with $t_0=1$ and $(\Xi,\Xi_j)=(\Phi,\Phi_j)$, therefore
gives \eqref{eq:jump} after fixing $\omega_{<s}$.
At $s=0$ all scales are fresh and $Q_{0-}=\Phi(r_0)$.
The same calculation gives zero means for the first-order and mixed
terms, while the scale-$j$ diagonal term contributes
$\tau_j^2\Gamma_j(r_0)$.
\end{proof}

\section{Richardson dispersion of trajectories}\label{sec:trajectories}

We now count the gains. The lemma below follows from the fact that at every time some scale is comparable to the
current separation, and the window of that scale contributes
$\tau_j\Gamma_j\ge c\kappa$, while all other scales contribute
nonnegatively by Lemma~\ref{lem:gain}. 

\begin{lemma}[Accumulation of trajectory gains]\label{lem:count}
For sufficiently small $\kappa$, independently of $J$ and $L\ge2$,
there is $c_0>0$ such that, for every $t\in[0,1]$ and on every
realization with $0<R_s\le1/4$ for all $s\in[0,t]$,
\begin{equation}\label{eq:count}
 \sum_{j=j_0}^J\ 
 \sum_{\substack{0\le k<K_j\\k\tau_j\le t}}
       \tau_j^2\Gamma_j(R_{k\tau_j})\ge c_0\kappa t.
\end{equation}
At $t=1$ the sum includes only window starts strictly below $1$, as
already enforced by $k<K_j$.
\end{lemma}

\begin{proof}
For $s\in[0,t)$ let $s_j(s)=\tau_j\lfloor s/\tau_j\rfloor$,
the start of the active scale-$j$ window. Suppose first that
$R_s\ge\ell_J$. Choose $j$ with $\ell_j\le R_s<2\ell_j$ if
$R_s\le\ell_{j_0}$ (then $j_0\le j\le J$, because
$\ell_J\le R_s\le\ell_{j_0}$); if $R_s>\ell_{j_0}$ choose $j=j_0$.
In either case
\[
 1\le R_s/\ell_j\le C_1,
 \qquad C_1=\max(2,(4\ell_{j_0})^{-1}).
\]
By \eqref{eq:radiusdet},
\[
 |R_{s_j(s)}^\beta-R_s^\beta|
 \le C_H\tau_j\le C_H\kappa R_s^\beta.
\]
Choose $\kappa$ so small that this implies
$R_{s_j(s)}/R_s\in[1/2,2]$. It follows that
$R_{s_j(s)}/\ell_j\in[1/2,2C_1]$ and
$R_{s_j(s)}\ge\ell_j/2\ge\ell_J/2>2\varepsilon$.
Also $\ell_j\le2R_{s_j(s)}\le L\widetilde S(R_{s_j(s)})$.
%The cutoff is one there, with the derivatives needed in \eqref{eq:comparablegain} equal to those of the constant function even at the plateau endpoint. 
Thus, with $r=R_{s_j(s)}$, we have
\[
 2\varepsilon<r\le1/4,\qquad
 r/\ell_j\in[1/2,2C_1],\qquad w_j(r)=1.
\]
Applying \eqref{eq:comparablegain} with $a=1/2$ and $b=2C_1$
therefore gives
\[
 \tau_j\Gamma_j(R_{s_j(s)})
 \ge c_{1/2,\,2C_1}\kappa.
\]
Since $C_1$ depends only on the fixed model parameters, this
constant is independent of $J,L,\kappa$, and $s$; we denote it by $c$.

If $R_s<\ell_J$, choose $j=J$. The Lipschitz bound in
\eqref{eq:holder}, applied forward and backward on this window, gives
\[
 e^{-C\kappa}\le R_{s_J(s)}/R_s\le e^{C\kappa}.
\]
For small enough $\kappa$, $R_{s_J(s)}\le2\ell_J$, and
\eqref{eq:finestgain} gives the same positive bound. Every radius
used is at most $1/4$ by the local condition. Since the remaining
$\Gamma_j$ are nonnegative, in both cases
\[
 \sum_j\tau_j\Gamma_j(R_{s_j(s)})\ge c_0\kappa.
\]
Integrating the preceding inequality gives
\[
 c_0\kappa t
 \le \int_0^t\sum_j\tau_j\Gamma_j(R_{s_j(s)})\,\dd s.
\]
For each fixed $j$, the integrand
$\tau_j\Gamma_j(R_{s_j(s)})$ is constant on $I_{jk}$, with value
$\tau_j\Gamma_j(R_{k\tau_j})$. Consequently,
\begin{align*}
 \int_0^t\sum_j\tau_j\Gamma_j(R_{s_j(s)})\,\dd s
 &=
 \sum_j\sum_{\substack{0\le k<K_j\\k\tau_j\le t}}
 |I_{jk}\cap[0,t]|\,
 \tau_j\Gamma_j(R_{k\tau_j})\\
 &\le
 \sum_j\sum_{\substack{0\le k<K_j\\k\tau_j\le t}}
 \tau_j^2\Gamma_j(R_{k\tau_j}),
\end{align*}
where $|I_{jk}\cap[0,t]|$ denotes the length of the intersection.
The inequality uses $|I_{jk}\cap[0,t]|\le\tau_j$ and
$\Gamma_j\ge0$. This proves \eqref{eq:count}, including windows
that have started but have not yet ended at time $t$.
\end{proof}

\subsection{Proof of Theorem~\ref{thm:main}}\label{sec:proof}

\begin{proof}[Proof of Theorem~\ref{thm:main}]
First restrict $\kappa$ by all the fixed smallness conditions in
Lemmas~\ref{lem:regularity} and \ref{lem:count}. Fix a time satisfying
\eqref{eq:localtime}. Then $R_s\le1/4$ for all $s\in[0,t]$ by
Lemma~\ref{lem:regularity}, so by Remark~\ref{rem:attained} the
configuration $\mathcal S_s$ is admissible, and $Q_s$ is defined, for
every $s\in[0,t]$ and every realization.

For fixed $J$, the bound
$R_s\ge r_0e^{-C\lambda_Js}$ established in the proof of
Lemma~\ref{lem:regularity} gives
$\log R_s\ge\log r_0-C\lambda_Jt$ for every $s\in[0,t]$. Thus $\Phi(R_s)$
is bounded on $[0,t]$, uniformly in the realization. The corrector $Q-\Phi$ is
bounded by \eqref{eq:endpointbound}, and the set of renewal times
$\{k\tau_j\}$ in $[0,t]$ is finite and deterministic. On each
closed interval between consecutive renewal times, $s\mapsto Q_s$ is
continuous, with the window data and phases of that interval
(one-sided limits at the endpoints): the particles and centers are
continuous in time by Lemma~\ref{lem:betweenrenewals}, the phases are
affine, and $Q$ is a continuous function of the configuration and the
phases. Therefore, we can
integrate the deterministic identity
\eqref{eq:pathwiseevolution} on each interval between consecutive
renewal times and add the jumps at the renewal times, obtaining, on
every realization,
\[
 Q_t=\Phi(r_0)+\sum_{\zeta\le t}(Q_\zeta-Q_{\zeta-})+\int_0^t(E_0(s)+E_3(s))\dd s,
\]
with the sum over the renewal times in $[0,t]$. All terms are bounded
uniformly in the realization, so we may integrate over the random
data as well. For each renewal time $\zeta$, Fubini's theorem gives
$\E[Q_\zeta-Q_{\zeta-}]=\E\big[\E_\zeta[Q_\zeta-Q_{\zeta-}]\big]$, and
Proposition~\ref{prop:renewal} evaluates the inner average.
At $t=1$ use the left limits $Q_{1-}$ if necessary. Then
\begin{equation}\label{eq:expectationidentity}
 \E Q_t =\Phi(r_0) +\E\sum_{j,k:\,k\tau_j\le t}
              \tau_j^2\Gamma_j(R_{k\tau_j})+\E\int_0^t(E_0(s)+E_3(s))\dd s,
\end{equation}
where the window-index restrictions are as in
Lemma~\ref{lem:count}. The time-zero renewals are included in the
sum. Since $r_0\ge\ell_J>2\varepsilon$, $\Phi(r_0)=r_0^\beta$.
It follows from \eqref{eq:errors} and \eqref{eq:count} that
\begin{equation}\label{eq:prechoice}
 \E Q_t\ge r_0^\beta+
       [c_0\kappa-CL^{-\beta}-C\kappa^2(1+\log L)^2]t.
\end{equation}

Choose a fixed $K\ge1$ sufficiently large that $C/K\le c_0/4$,
and set
\begin{equation}\label{eq:Lchoice}
 L=(K/\kappa)^{1/\beta}.
\end{equation}
Decrease the permissible $\kappa_0$ so that $L\ge2$ and
$C\kappa(1+\log L)^2\le c_0/4$ whenever
$\kappa\le\kappa_0$. This is possible because
$\kappa(1+|\log\kappa|)^2\to0$. Equation \eqref{eq:prechoice}
then becomes
\begin{equation}\label{eq:Qlower}
 \E Q_t\ge r_0^\beta+(c_0/2)\kappa t.
\end{equation}

For this choice of $L$, \eqref{eq:endpointbound} and
\eqref{eq:phiproperties} imply
\[
 Q_t\le(1+\delta_\kappa)R_t^\beta
                     +\delta_\kappa\ell_J^\beta,
 \qquad \delta_\kappa=C\kappa(1+|\log\kappa|),
\]
after increasing the fixed constant $C$ and decreasing $\kappa_0$
if necessary. Combining this with \eqref{eq:Qlower} proves
\eqref{eq:maincorrected}. Make a final decrease of $\kappa_0$ to
ensure $\delta_\kappa\le1/4$. Since $\ell_J\le r_0$, it follows
that $\E R_t^\beta\ge c(r_0^\beta+\kappa t)$.
Convexity of $q\mapsto q^{2/\beta}$ now gives
\[
 \E R_t^2=\E[(R_t^\beta)^{2/\beta}]
 \ge(\E R_t^\beta)^{2/\beta}
 \ge c(r_0^\beta+\kappa t)^{2/\beta}.
\]
Finally, integrating \eqref{eq:radiusdet} gives
$R_t^\beta\le r_0^\beta+C_Ht$ on every realization, which implies
the upper bound in \eqref{eq:momentbounds}. This completes the proof.
\end{proof}

\begin{proof}[Proof of Corollary~\ref{cor:cubic}]
When $\alpha=1/3$, $\beta=2/3$ and $2/\beta=3$.
If $t\ge r_0^\beta/\kappa$, then
$r_0^\beta\le\kappa t\le t$, since $\kappa_0<1$.
Apply \eqref{eq:momentbounds}.
\end{proof}

\section{Advection--diffusion equations with bounded incompressible drift}\label{sec:weak}

We
work on a torus $\T^n$ ($n=d$ or $n=2d$), where integrations by parts
produce no boundary terms. We recall two standard facts and prove the
consequences used later. Only the boundedness and the incompressibility
of the drift enter, and nothing in this section is quantitative.

\emph{Notation.} Functions are real-valued; $(\cdot,\cdot)$ and
$\|\cdot\|$ are the inner product and norm of $L^2(\T^n)$,
$H^1=H^1(\T^n)$, $H^{-1}$ is its dual, and $\langle\cdot,\cdot\rangle$ the
duality pairing. Let $I=[T_0,T_1]$ be a compact interval and $\mu>0$.
A \emph{drift} is a Borel map $\mathbf b:I\times\T^n\to\R^n$ such that
\begin{equation}\tag{D}\label{eq:D}
 |\mathbf b|\le B_0,\qquad
 \int_{\T^n}\mathbf b(s,Y)\cdot\nabla\varphi(Y)\dd Y=0\quad
 \text{for all }\varphi\in C^\infty(\T^n)\text{ and a.e. }s\in I .
\end{equation}
Since $C^\infty(\T^n)$ contains a countable set that is dense in $C^1$, there
is one null set $N\subset I$ such that the identity in \eqref{eq:D} holds
for all $\varphi\in C^\infty$ and all $s\notin N$; for such $s$ it holds for
all $\varphi\in W^{1,1}(\T^n)$, since $C^\infty$ is dense in $W^{1,1}$ and
$\mathbf b(s)\in L^\infty$. If
$v,\varphi\in H^1$, then $v\varphi\in W^{1,1}$ with
$\nabla(v\varphi)=\varphi\nabla v+v\nabla\varphi$, and therefore, for $s\notin N$,
\begin{equation}\label{eq:antisym}
 (\mathbf b(s)\cdot\nabla v,\varphi)=-(\mathbf b(s)\cdot\nabla\varphi,v),
 \qquad(\mathbf b(s)\cdot\nabla v,v)=0 .
\end{equation}

We consider the equation
\begin{equation}\label{eq:adv}
 \partial_sw+\mathbf b\cdot\nabla w=\mu\Delta w\qquad\text{on }I\times\T^n .
\end{equation}

\begin{definition}[Energy solutions]\label{def:energy}
An \emph{energy solution} of \eqref{eq:adv} on $I$ is a function
$w\in L^2(I;H^1)$ whose distributional derivative satisfies
$\partial_sw\in L^2(I;H^{-1})$ and such that, for a.e.\ $s\in I$,
\begin{equation}\label{eq:weakform}
 \langle\partial_sw(s),\varphi\rangle+(\mathbf b(s)\cdot\nabla w(s),\varphi)
 +\mu(\nabla w(s),\nabla\varphi)=0\qquad\text{for all }\varphi\in H^1 .
\end{equation}
\end{definition}

We use two standard facts.

\smallskip\noindent
(E1) \emph{The Lions--Magenes lemma.} If $w,\widetilde w\in L^2(I;H^1)$
and $\partial_sw,\partial_s\widetilde w\in L^2(I;H^{-1})$, then, after
modification on null sets, $w,\widetilde w\in C(I;L^2)$; there are
$w_k\in C^1(I;H^1)$ with $w_k\to w$ in $L^2(I;H^1)$ and in $C(I;L^2)$,
and $\partial_sw_k\to\partial_sw$ in $L^2(I;H^{-1})$; and
$s\mapsto(w(s),\widetilde w(s))$ is absolutely continuous with derivative
$\langle\partial_sw,\widetilde w\rangle+\langle\partial_s\widetilde w,w\rangle$.
See \cite[\S5.9.2, Theorem~3 and its proof]{Eva10}, where the case
$\widetilde w=w$ is stated; the general case follows by polarization. The
approximations $w_k$ are obtained there by extending $w$ to a larger
interval and mollifying in time; extending by even reflection about $T_0$
and $T_1$ preserves the hypotheses, which gives the stated convergence up
to the endpoints. The proof is written for $H^1_0$ of a bounded domain
and applies verbatim on $\T^n$.

\smallskip\noindent
(E2) \emph{Existence.} For every $w_0\in L^2$ there is an energy solution
of \eqref{eq:adv} on $I$ with $w(T_0)=w_0$. This is proved by Galerkin
approximation, as in \cite[\S7.1.2, Theorems~1--3]{Eva10}, with the
spaces spanned by $\{e^{2\pi\mathrm ik\cdot Y}:|k|\le N\}$ in place of
those spanned by the Dirichlet eigenfunctions (the orthogonal projections
onto these spaces are contractions in $H^1$). The energy estimate holds
because $(\mathbf b\cdot\nabla v,v)=0$ by \eqref{eq:antisym}, and only the
boundedness of $\mathbf b$ is used. The weak formulation obtained there
(for each $\varphi\in H^1$, for a.e.\ $s$) implies \eqref{eq:weakform}
since $H^1$ is separable.

\smallskip
We always use the representative in $C(I;L^2)$.
The proofs of the five lemmas of this section are standard; they are
collected in Appendix~\ref{sec:appendix}.

\begin{lemma}[Basic properties]\label{lem:weakbasic}
Let $\mathbf b$ satisfy \eqref{eq:D} and let $w$ be an energy solution
of \eqref{eq:adv} on $I$.
\begin{enumerate}[(a)]
\item $\|w(s)\|^2+2\mu\int_{T_0}^s\|\nabla w\|^2=\|w(T_0)\|^2$ for
 $s\in I$. Consequently the energy solution with a given initial value
 is unique, and it depends linearly on the initial value.
\item $\int w(s)=\int w(T_0)$ for $s\in I$, and if $w(T_0)\ge0$ a.e., then
 $w(s)\ge0$ a.e.\ for every $s\in I$.
\item (Evolution identity.) Let $T_0\le s_1<s_2\le T_1$ and let
 $G:[s_1,s_2]\times\T^n\to\R$ be such that $G$, $\partial_sG$, $\nabla G$
 and $\nabla^2G$ exist and are continuous on $[s_1,s_2]\times\T^n$
 (one-sided in $s$ at the endpoints). Then
 \[
 \int G(s_2)w(s_2)-\int G(s_1)w(s_1)
 =\int_{s_1}^{s_2}\!\!\int\big(\partial_sG+\mathbf b\cdot\nabla G
 +\mu\Delta G\big)w\dd Y\dd s .
 \]
\end{enumerate}
\end{lemma}

For an energy solution $w$ with $w(T_0)\ge0$ and $\int w(T_0)=1$, each
$w(s)$ is therefore a probability density. For $T_0\le s_1<s_2\le T_1$
and Borel sets $A,B\subset\T^n$ we put
\begin{equation}\label{eq:Pi}
 \Pi_{s_1,s_2}(A;B):=\sup\Big\{\int_Bw(s_2)\Big\},
\end{equation}
the supremum over all energy solutions $w$ on $[s_1,s_2]$ whose
initial value $w(s_1)\in L^2$ is a probability density vanishing a.e.\
outside $A$ (and $\Pi_{s_1,s_2}(A;B):=0$ if $A$ is a null set).

\begin{lemma}[Transfer]\label{lem:transfer}
Let $w$ be an energy solution on $[s_1,s_2]$ with $w(s_1)\ge0$. Then for
Borel sets $A,B\subset\T^n$, with $A^c=\T^n\setminus A$,
\[
 \int_Bw(s_2)\le\int_Aw(s_1)+\Pi_{s_1,s_2}(A^c;B)\int_{A^c}w(s_1) .
\]
\end{lemma}

\begin{lemma}[Duality]\label{lem:duality}
Let $\mathbf b$ satisfy \eqref{eq:D} on $[0,T]$ and put
$\overleftarrow{\mathbf b}(s,Y):=-\mathbf b(T-s,Y)$. Let $w$ be an energy
solution of \eqref{eq:adv} with drift $\mathbf b$ on $[0,T]$, and
$\varrho$ an energy solution with drift $\overleftarrow{\mathbf b}$ on
$[0,T]$. Then $(w(T),\varrho(0))=(w(0),\varrho(T))$.
\end{lemma}

\begin{lemma}[Stability]\label{lem:stability}
Let $\mathbf b_k,\mathbf b$ satisfy \eqref{eq:D} on $I$ with the same
$B_0$, and $\mathbf b_k\to\mathbf b$ a.e.\ on $I\times\T^n$. Let $w_0\in L^2$
and let $w_k,w$ be the energy solutions with drifts $\mathbf b_k,\mathbf b$
and initial value $w_0$. Then $w_k(s)\to w(s)$ weakly in $L^2$ for
every $s\in I$.
\end{lemma}

\begin{lemma}[Products and isometries]\label{lem:products}
\begin{enumerate}[(a)]
\item Let $b$ satisfy \eqref{eq:D} on $I\times\T^d$, and let
 $\theta_1,\theta_2$ be energy solutions of \eqref{eq:adv} on $\T^d$ with
 drift $b$. Then $\mathbf b(s,y,y'):=(b(s,y),b(s,y'))$ satisfies
 \eqref{eq:D} on $I\times\T^{2d}$, with $B_0$ replaced by $\sqrt2B_0$ (the
 gradient in \eqref{eq:D} and \eqref{eq:adv} being that of $\T^{2d}$, so
 that $\mathbf b\cdot\nabla=b(s,y)\cdot\nabla_y+b(s,y')\cdot\nabla_{y'}$), and
 $(\theta_1\otimes\theta_2)(s,y,y'):=\theta_1(s,y)\theta_2(s,y')$ is the energy
 solution on $\T^{2d}$ with drift $\mathbf b$ and initial value
 $\theta_1(T_0)\otimes\theta_2(T_0)$. The constant $1$ is an energy
 solution for every drift.
\item Let $\iota(Y)=gY+a$ with $g\in O(n)$, $g\mathbb Z^n=\mathbb Z^n$, and
 $a\in\R^n$, an isometry of $\T^n$. If $w$ is an energy solution with
 drift $\mathbf b$, then $w\circ\iota^{-1}$ is the energy solution with
 drift $(\iota_*\mathbf b)(s,Y):=g\,\mathbf b(s,\iota^{-1}Y)$, which satisfies
 \eqref{eq:D}, and initial value $w(T_0)\circ\iota^{-1}$.
\end{enumerate}
\end{lemma}

\section{Dispersion of diffusing pairs}\label{sec:pair}

For anomalous dissipation we need a separation estimate for solutions
of the two-point advection--diffusion equation, in both directions of
time. We use the same corrector and jump identities as in
Section~\ref{sec:trajectories}; exponential moment bounds replace the
deterministic control of separation within a window.

\subsection{The two-point equation and separation growth}
Fix $J\ge j_0$, $t_0\in(0,1]$ and $\kappa\le\kappa_*$, and let
$\mu=\mu_J$ be as in \eqref{eq:diffusivity}. Let $v$ denote either
\begin{enumerate}
\item[(F)] the forward field $u=u^{(J)}$ on $[0,t_0]$, with the
 renewal times $k\tau_j$ and the phases $h_j$ of Section~\ref{sec:corrector}; or
\item[(B)] the reversed field $\overleftarrow u$ on $[0,t_0]$ of
 Definition~\ref{def:reversed}, with its renewal times and phases.
\end{enumerate}
We use the same notation for renewal times in the two cases. Let
$\mathcal Z$ be the (finite, deterministic) set of renewal times of
all scales in $(0,t_0)$. No renewal is performed at $t_0$: the window
data and phases at $t_0$ are their left limits. For $s\in(0,t_0)$ let
$s_j(s)$ be the largest renewal time of scale $j$ in $(0,s]$, or $0$ if
there is none, and $s_j(t_0):=s_j(t_0-)$; in both cases
$0\le s-s_j(s)\le\tau_j$. All renewal times lie in $\tau_J\mathbb N$ in
case (F) and in $t_0-\tau_J\mathbb N$ in case (B) (Lemma~\ref{lem:nested});
we call the open intervals between consecutive points of
$\mathcal Z\cup\{0,t_0\}$ \emph{slabs}. On each slab, $v$ is $C^3$ in $x$
and $C^1$ in time, with one-sided limits at the endpoints
(Lemma~\ref{lem:betweenrenewals} and Definition~\ref{def:reversed}).

We write $Y=(y,y')\in\T^d\times\T^d=\T^{2d}$, $R(Y)=\dist(y,y')$, and, when
$R(Y)<1/2$, $z=z(Y)$ for the shortest representative of $y-y'$, so
$|z|=R$. Put
\begin{gather*}
 \mathbf v(s,Y)=(v(y,s),v(y',s)),\qquad
 \nabla_Y=(\nabla_y,\nabla_{y'}),\qquad
 \Delta_Y=\Delta_y+\Delta_{y'},\\
 \mathcal T_sG=\mathbf v(s)\cdot\nabla_YG+\mu\Delta_YG .
\end{gather*}
Thus $\nabla_Y$ and $\Delta_Y$ are the gradient and the Laplacian of
$\T^{2d}$, so that
$\mathbf v\cdot\nabla_Y=v(y,s)\cdot\nabla_y+v(y',s)\cdot\nabla_{y'}$, and
$\mathbf v$ satisfies \eqref{eq:D} on $\T^{2d}$
(Lemmas~\ref{lem:betweenrenewals} and~\ref{lem:products}(a)). We study
energy solutions of the \emph{two-point equation}
\begin{equation}\label{eq:doubled}
 \partial_s\varrho+\mathbf v\cdot\nabla_Y\varrho=\mu\Delta_Y\varrho
 \qquad\text{on }[0,t_0]\times\T^{2d},
\end{equation}
which is the tracer equation with diffusion in the $2d$ variables $Y$,
the field acting on $y$ and on $y'$ separately. For a function $G$ of $z$
alone, $\mathcal T_sG=\delta v\cdot\nabla_zG+2\mu\Delta_zG$ with
$\delta v=v(y,s)-v(y',s)$, because $\nabla_{y'}z=-I$.

For $0<\varsigma\le1/4$ let $\mathfrak D_\varsigma$ be the set of probability
densities $\varrho_0\in L^2(\T^{2d})$ that vanish a.e.\ on $\{R>\varsigma\}$.
Thus these densities describe initial pairs of points whose separation
is at most $\varsigma$: the two particles start close together.
For $\varrho_0\in\mathfrak D_\varsigma$ we denote by $\varrho$ the energy
solution of \eqref{eq:doubled} on $[0,t_0]$ with $\varrho(0)=\varrho_0$, and
write $\varrho_s=\varrho(s)$ and $\varrho_s(B)=\int_B\varrho_s$. By
Lemma~\ref{lem:weakbasic}, each $\varrho_s$ is a probability density, and
$\Pi_{s',s}$ of \eqref{eq:Pi} refers to \eqref{eq:doubled}.

Before taking expectations, we record the dependence of the solutions
on the random data. Besides measurability, we shall need the fact that
the density at a renewal time does not depend on the data newly
introduced at that time.

\begin{lemma}[Dependence on the data]\label{lem:datadependence}
Fix $J$, $t_0$, case (F) or (B), and $\varrho_0\in L^2(\T^{2d})$. For each
$s\in[0,t_0]$, $\varrho_s$ is determined by the restriction of $\mathbf v$
to $[0,s)$, and for every $G\in L^2(\T^{2d})$, $\int G\varrho_s$ is a
continuous function of the finitely many data variables of the windows
meeting $[0,t_0]$. In particular, the function is Borel. The same holds for the energy
solutions of \eqref{eq:scalar}, and $\mathbf D_J(\theta_0)$ is a Borel
function of the data. All of this remains true, with the same proof,
when the diffusivity $\mu_J$ is replaced by an arbitrary $\mu>0$.
\end{lemma}

\begin{proof}
The restriction of $\varrho$ to $[0,s]$ is the energy solution on $[0,s]$
(Lemma~\ref{lem:weakbasic}(a)); hence two drifts that agree a.e.\ on
$[0,s)\times\T^{2d}$ give the same $\varrho_s$. The data space
of these windows is a compact metric space, a finite product of
cubes $Q_j$ and intervals $[-M,M]$ (Section~\ref{sec:model}), so sequential
continuity suffices, and by Lemma~\ref{lem:stability} it suffices to show
that $\mathbf v(s,Y)$ is a continuous function of the data for every $Y$
and every $s\notin\mathcal Z\cup\{0,t_0\}$; the bound in \eqref{eq:D} is
uniform by Lemma~\ref{lem:betweenrenewals}. On each window of scale $j$
the centers solve \eqref{eq:centers}, whose right-hand side $u_{<j}$ is
Lipschitz in $x$ with a constant valid for every datum
(Lemma~\ref{lem:regularity}) and, by induction on $j$,
continuous in the data uniformly in $(x,t)$. By Gronwall's inequality the centers depend
continuously on the data, uniformly in time, and since the $W_j^\nu$ are
uniformly continuous, $u_j(x,t)$ depends continuously on the data,
uniformly in $(x,t)$. In case (B), $\overleftarrow u(x,s)=-u(x,t_0-s)$ for
$s\notin\mathcal Z\cup\{0,t_0\}$. The same argument applies to
\eqref{eq:scalar}; then $\theta^{(J)}(t_0)$ depends weakly continuously on the
data, so $\|\theta^{(J)}(t_0)\|^2$ is lower semicontinuous and, by
\eqref{eq:energy}, $\mathbf D_J(\theta_0)$ is Borel.
\end{proof}

\begin{theorem}[Richardson growth near the diagonal]\label{thm:pair}
Let $d\ge2$ and assume \eqref{eq:ND} and \eqref{eq:coefficients}. There
exist $t_1\in(0,1]$, $c_E>0$, $\kappa_1\in(0,\kappa_*]$ and
$c_5,C,C_\delta>0$, depending only on $d,\alpha,\psi,R,M,\sigma$,
and for each $\kappa\in(0,\kappa_1]$ a number
$\varepsilon_1(\kappa)\in(0,1]$, such that the following holds for all
$\kappa\in(0,\kappa_1]$, $\varepsilon_0\in(0,\varepsilon_1(\kappa)]$,
$J\ge J_0:=\max(j_0,5)$ and $t_0\in(0,1]$, in both cases (F) and (B). If
$\varrho_0\in\mathfrak D_\varsigma$ with $\varsigma\le r_\sharp/8$ does not depend on the
data, then for every
$t\in(0,\min(t_0,t_1)]$
\begin{equation}\label{eq:pairmain}
\begin{gathered}
 (1+\delta_\kappa)\,\E\!\int R^\beta\varrho_t\ \ge\ c_5\kappa t-\mathrm{Err}_J,\\
 \mathrm{Err}_J:=C\big(1+\log(2\varepsilon/r_\sharp)\big)\ell_J^\beta+Ce^{-c_E/\mu},
\end{gathered}
\end{equation}
where $\delta_\kappa:=C_\delta\kappa(1+|\log\kappa|)\le\frac14$. Since
$\varepsilon/r_\sharp$ does not depend on $J$, $\mathrm{Err}_J\to0$ as
$J\to\infty$, uniformly in $t_0$, $t$ and $\varrho_0$. Whenever
$\mathrm{Err}_J\le c_5\kappa t/2$,
\begin{equation}\label{eq:pairmoments}
 \E\!\int R^\beta\varrho_t\ge\tfrac25c_5\kappa t,\qquad
 \E\!\int R^2\varrho_t\ge\big(\tfrac25c_5\kappa t\big)^{2/\beta}.
\end{equation}
\end{theorem}

We prove the theorem in the remainder of this section. First we
derive the evolution identity and estimate the errors introduced by
the spatial cutoff and by diffusion. We then compare the mass in each
separation range with its mass at the preceding renewal, and sum the
resulting gains. Throughout, the cutoff parameter from
Section~\ref{sec:testfunctions} is $L=(K/\kappa)^{1/\beta}$ and
$N_L=1+\log L$; the constant $K$ will be fixed in the proof.

\subsection{The cut-off observable and its evolution}

The corrected observable is defined near the diagonal. To use it as
a test function on the torus, we first localize it in that region.
Fix a nonincreasing $\chi_\circ\in C^\infty([0,\infty))$ with
$\chi_\circ=1$ on $[0,\frac18]$ and $\chi_\circ=0$ on $[\frac3{16},\infty)$;
then $Y\mapsto\chi_\circ(R(Y))$ is smooth on $\T^{2d}$. For $s\in[0,t_0]$
let $\Theta_s$ and $h(s)$ be the current window data and phases, and
$\Theta_{s-}$, $h(s-)$ their left limits. Identifying the pair
$Y=(y,y')$ with the particles $(x,\widetilde x)$ of a configuration, put
\[
 Q^\circ_s(Y):=
 \begin{cases}
 \chi_\circ(R(Y))\,Q^\sharp\big((Y,\Theta_s),h(s)\big),&R(Y)\le\tfrac14,\\
 0,&R(Y)>\tfrac14,
 \end{cases}
\]
and define $Q^\circ_{s-}$ in the same way with $(\Theta_{s-},h(s-))$. By
Corollary~\ref{cor:coretransfer}, $Q^\sharp$ is defined and continuous up
to the diagonal $y=y'$. The functions $Q^\sharp$, $\Phi^\sharp(z)$, $E_0$, $E_3$
and their derivatives are defined only on $\{R<\frac14\}$. Every product of
one of them with $\chi_\circ\circ R$ or with a derivative of $\chi_\circ\circ R$
(as in $Q^\circ_s$ and $\mathcal K_s$ below) is extended by $0$ to
$\{R\ge\frac14\}$; since $\chi_\circ\circ R$ vanishes on $\{R\ge\frac3{16}\}$, this
does not affect smoothness.

\begin{lemma}[Evolution of the cut-off observable]\label{lem:evolution}
Assume \eqref{eq:eps0small}. For fixed data, every $\varrho_0\in L^2(\T^{2d})$
and $t\in(0,t_0]$,
\begin{equation}\label{eq:evolution}
\begin{aligned}
 \int Q^\circ_t\varrho_t={}&\int Q^\circ_0\varrho_0
 +\int_0^t\!\!\int\mathcal K_s\varrho_s\dd s\\
 &+\sum_{\tilde{\zeta}\in\mathcal Z,\ \tilde{\zeta}\le t}
 \int\big(Q^\circ_{\tilde{\zeta}}-Q^\circ_{\tilde{\zeta}-}\big)\varrho_{\tilde{\zeta}},
\end{aligned}
\end{equation}
where, with everything evaluated at $((Y,\Theta_s),h(s))$,
\[
 \mathcal K_s:=\chi_\circ(R)\big[E_0+E_3+2\mu\Delta\Phi^\sharp(z)
 +\mu\Delta_Y(Q^\sharp-\Phi^\sharp)\big]+\mathcal K^{\rm cut}_s,
\]
\[
 \mathcal K^{\rm cut}_s:=Q^\sharp\,\mathcal T_s\big(\chi_\circ\circ R\big)
 +2\mu\big(\nabla_y(\chi_\circ\circ R)\cdot\nabla_yQ^\sharp
 +\nabla_{y'}(\chi_\circ\circ R)\cdot\nabla_{y'}Q^\sharp\big),
\]
$E_0,E_3$ are those of \eqref{eq:exactevolution} for
$(\Phi^\sharp,\Phi^\sharp_j)$, and $\mathcal K^{\rm cut}_s=0$ unless
$\frac18\le R\le\frac3{16}$. Moreover, if $\varrho_0$ does not depend on
the data, then for every $\tilde{\zeta}\in\mathcal Z$
\begin{equation}\label{eq:jumpaverage}
 \E\int\big(Q^\circ_{\tilde{\zeta}}-Q^\circ_{\tilde{\zeta}-}\big)\varrho_{\tilde{\zeta}}
 =\E\int\chi_\circ(R)\sum_{j\in\mathcal J(\tilde{\zeta})}\tau_j^2\Gamma^\sharp_j(R)\,\varrho_{\tilde{\zeta}} .
\end{equation}
\end{lemma}

\begin{proof}
We apply the weak evolution identity separately between successive
renewals and then sum. We first justify this operation, compute the
continuous contribution, and finally average the jumps.

On $\{R<\frac14\}$, $Q^\sharp$ is $C^2$ in $Y$ and $C^1$
in the centers and phases, and these derivatives are jointly continuous:
$\Phi^\sharp_b(z)$ is smooth (Lemma~\ref{lem:core}(b)), $u_a\in C^3$,
$\B_au_b\in C^2(\mathcal E)$ (Definition~\ref{def:operators}), and $q_{ab}$
is $C^1$. On the closure of a slab, the coefficients and shifts are
constant, the centers are $C^1$ in time with one-sided limits, and the
phases are affine. Since $\chi_\circ\circ R$ is smooth and vanishes for
$R\ge\frac3{16}$, the function equal to $Q^\circ$ with the data of a slab,
extended to the closed slab by one-sided limits, satisfies the
hypotheses of Lemma~\ref{lem:weakbasic}(c).

Inside a slab, write
\[
 c(Y):=\chi_\circ(R(Y)),\qquad
 q_s(Y):=Q^\sharp\big((Y,\Theta_s),h(s)\big),
 \qquad Q^\circ_s=cq_s.
\]
At fixed $Y$, only the centers and phases vary with $s$; the
coefficients are constant. Since $\dot h_j=-1$, the chain rule gives
\[
 \partial_s q_s
 =\sum_{m,i}\dot P_i^m(s)\cdot\nabla_{P_i^m}Q^\sharp
   -\sum_j\partial_{h_j}Q^\sharp.
\]
Here and below the derivatives of $Q^\sharp$ are evaluated at
$((Y,\Theta_s),h(s))$. The phase term is $\mathcal D_hQ^\sharp$.
The center-motion term, together with the motion of the pair, is
the configuration generator:
\[
 \B Q^\sharp
 =\mathbf v\cdot\nabla_YQ^\sharp
  +\sum_{m,i}\dot P_i^m(s)\cdot\nabla_{P_i^m}Q^\sharp.
\]
In case (B), the centers and generator in these formulas are the
reversed ones of Lemma~\ref{lem:reversed}(b). Thus, exactly as in the
proof of \eqref{eq:pathwiseevolution},
\[
 (\partial_s+\mathbf v\cdot\nabla_Y)q_s
 =(\mathcal D_h+\B)Q^\sharp=E_0+E_3.
\]
The last equality follows from \eqref{eq:exactevolution} and
Corollary~\ref{cor:coretransfer} for $y\ne y'$. Both sides are
continuous up to $y=y'$, by the regularity just established and
Corollary~\ref{cor:coretransfer}, so the identity extends to the diagonal.

Since $c$ is independent of time, the product rule yields
\[
\begin{aligned}
 (\partial_s+\mathcal T_s)(cq_s)
 ={}&c(\partial_s+\mathbf v\cdot\nabla_Y)q_s
       +q_s\,\mathbf v\cdot\nabla_Yc\\
    &+\mu c\,\Delta_Yq_s+\mu q_s\,\Delta_Yc
       +2\mu\nabla_Yc\cdot\nabla_Yq_s\\
 ={}&c\big(E_0+E_3+\mu\Delta_Yq_s\big)
       +q_s\mathcal T_sc
       +2\mu\nabla_Yc\cdot\nabla_Yq_s.
\end{aligned}
\]
To separate the main observable from its correctors, write
\[
 \Delta_Yq_s
 =\Delta_Y\Phi^\sharp(z)
   +\Delta_Y\big(q_s-\Phi^\sharp(z)\big).
\]
On functions of $z=y-y'$, $\nabla_y=\nabla_z$ and
$\nabla_{y'}=-\nabla_z$, so each Laplacian contributes the same term:
\[
 \Delta_Y\Phi^\sharp(z)
 =(\Delta_y+\Delta_{y'})\Phi^\sharp(z)
 =2\Delta_z\Phi^\sharp(z).
\]
Consequently,
\[
\begin{aligned}
 (\partial_s+\mathcal T_s)Q^\circ_s
 ={}&c\Big[E_0+E_3+2\mu\Delta_z\Phi^\sharp(z)
           +\mu\Delta_Y\big(q_s-\Phi^\sharp(z)\big)\Big]\\
    &+q_s\mathcal T_sc
       +2\mu\nabla_Yc\cdot\nabla_Yq_s
 =\mathcal K_s.
\end{aligned}
\]
Indeed, the last two terms are precisely $\mathcal K_s^{\rm cut}$,
since
\[
 \nabla_Yc\cdot\nabla_Yq_s
 =\nabla_yc\cdot\nabla_yq_s
  +\nabla_{y'}c\cdot\nabla_{y'}q_s.
\]
These are exactly the terms in which a spatial derivative falls on
the cutoff. All products are extended by zero as specified above.

We now sum over the slabs. Let $0=\tilde{\zeta}_0<\tilde{\zeta}_1<\dots<\tilde{\zeta}_n<\tilde{\zeta}_{n+1}=t$ be
the points of $\mathcal Z\cap(0,t)$ together with $0$ and $t$. Apply
Lemma~\ref{lem:weakbasic}(c) on each $[\tilde{\zeta}_k,\tilde{\zeta}_{k+1}]$ to the function
equal to $Q^\circ$ with the data of $(\tilde{\zeta}_k,\tilde{\zeta}_{k+1})$, which takes the
value $Q^\circ_{\tilde{\zeta}_k}$ at $\tilde{\zeta}_k$ and $Q^\circ_{\tilde{\zeta}_{k+1}-}$ at
$\tilde{\zeta}_{k+1}$. Summing, the boundary terms telescope up to the jumps at
the points of $\mathcal Z$; if $t\in\mathcal Z$, its jump is added, and
if $t\notin\mathcal Z$, then $Q^\circ_{t-}=Q^\circ_t$ (also for $t=t_0$,
by convention).

To obtain \eqref{eq:jumpaverage}, we condition on all data except
those refreshed at $\tilde{\zeta}$. By Lemma~\ref{lem:datadependence}, $\varrho_{\tilde{\zeta}}$
is determined by $\mathbf v$ on $[0,\tilde{\zeta})$: in case (F) by the data of
windows starting before $\tilde{\zeta}$, in case (B) by $\omega_{\tilde{\zeta}}^c$
(Lemma~\ref{lem:reversed}(c1)). In both cases it does not depend on the
fresh data $\omega_{\tilde{\zeta}}$. At a fixed $Y$ with $R(Y)\le\frac3{16}$,
Proposition~\ref{prop:jumpboth} gives
\[
 \E_{\tilde{\zeta}}\big[Q^\sharp_{\tilde{\zeta}}-Q^\sharp_{\tilde{\zeta}-}\big](Y)
 =\sum_{j\in\mathcal J(\tilde{\zeta})}\tau_j^2\Gamma^\sharp_j(R(Y)).
\]
Multiply by $\chi_\circ(R)\varrho_{\tilde{\zeta}}(Y)$, integrate in $Y$, and use Fubini's
theorem and $\E=\E\,\E_{\tilde{\zeta}}$. This is justified because the jump is
continuous in $Y$ and in the data and bounded uniformly in both (by
\eqref{eq:endpointbound}, for fixed $\kappa$ and $L$), because
$\|\varrho_{\tilde{\zeta}}\|_{L^1}\le\|\varrho_{\tilde{\zeta}}\|\le\|\varrho_0\|$ (Lemma~\ref{lem:weakbasic}(a)),
and because the map from the data to $\varrho_{\tilde{\zeta}}\in L^2(\T^{2d})$ is Borel
(Lemma~\ref{lem:datadependence}), so that
$(\omega,Y)\mapsto\varrho_{\tilde{\zeta}}[\omega](Y)$ has a jointly measurable version.
\end{proof}

The cutoff contributes only on an annulus separated from the
diagonal. On this annulus the singular derivative bounds become
uniform, and the diffusivity absorbs the growth in the number of
active scales.

\begin{lemma}[Cut-off error]\label{lem:cutoff}
Assume \eqref{eq:eps0small}, $\kappa N_L\le1$ and $J\ge J_0$. There
is $C_B$, independent of $J,L,\kappa,\varepsilon_0$, with
$|\mathcal K^{\rm cut}_s|\le C_B\one_{\{1/8\le R\le3/16\}}$.
\end{lemma}

\begin{proof}
Let $\frac18\le R\le\frac3{16}$. Since $J\ge J_0$, $\ell_J\le\frac1{32}$, so
$S=R$, ${m_\sharp}=R$ (as $r_\sharp\le\ell_J/16$) and $R\ge2\varepsilon$ (as
$2\varepsilon<\ell_J/4$); thus $S$ and ${m_\sharp}$ are comparable to $1$. By
Lemma~\ref{lem:core}(a) and \eqref{eq:phiproperties}, $\Phi^\sharp=\Phi=R^\beta$
there, and by
\eqref{eq:endpointbound} $|Q^\sharp|\le C$. Since $|v|\le C$
(Lemma~\ref{lem:betweenrenewals}), $|\mathcal T_s(\chi_\circ\circ R)|\le C$. For
the gradient, $|\nabla_y\Phi^\sharp|\le C$ by \eqref{eq:corederivs}; for
$b\in\mathcal I=\{b:\ell_b\le3LS\}$ (the other terms vanish), by
\eqref{eq:corederivs}, $|\delta u_b|\le CA_b\le C\la_b$ and
$\|\nabla u_b\|\le C\la_b$,
\[
 |\nabla_y(h_b\nabla\Phi^\sharp_b(z)\cdot\delta u_b)|
 \le Ch_b\big(|\delta u_b|+\|\nabla u_b\|\big)\le C\kappa ,
\]
and, using also \eqref{eq:Fsup} and $\|\nabla F_{ab}\|\le C\la_a\la_b$
(Step~2 of the proof of Lemma~\ref{lem:diffcorrector}),
$|\nabla_y(q_{ab}\B_a\B_b\Phi^\sharp_b)|\le Cq_{ab}\la_a\la_b
\le C\kappa^2\la_a/\max(\la_a,\la_b)$. With $\mathcal N=N_L+\log(S/\ell_J)\le N_L+J$
(as $S<1$ and $\ell_J=2^{-J}$, so $\log(S/\ell_J)\le J\log2\le J$),
$\#\mathcal I\le C\mathcal N$ and $\sum_a\la_a/\max(\la_a,\la_b)\le C\mathcal N$, so
$|\nabla_yQ^\sharp|\le C(1+\kappa\mathcal N+\kappa^2\mathcal N^2)\le C(1+J)^2$,
using $\kappa N_L\le1$. The same holds for $\nabla_{y'}$. Hence
$\mu|\nabla_YQ^\sharp|\le C\varepsilon_0\ell_J^{1+\alpha}(1+J)^2\le C$.
\end{proof}

\subsection{Exponential moment bounds}

We next control the transfer of mass between separation ranges.
This is needed because the gain is evaluated at renewal times,
whereas the density evolves throughout each window. By
Lemma~\ref{lem:core}(f), the rate $\tau_j\Gamma^\sharp_j$ is bounded below
by a multiple of $\kappa$ on the appropriate range of separations.
To use this bound at later times, we must exclude rapid changes of
separation, apart from a small exceptional mass.

The estimates below are parabolic off-diagonal bounds. We test
\eqref{eq:doubled} against an exponential whose spatial part varies
only across the annulus being crossed. The H\"older bound controls the
drift term; a suitable exponential parameter and a linear correction
in time control the diffusion terms. This produces a supersolution
and bounds the mass that crosses the annulus in either direction.

Fix $\lambda_\uparrow,\lambda_{\downarrow,2},\lambda_{\downarrow,3/2}\in C^\infty([0,\infty))$ with
$\lambda_\uparrow$ nondecreasing, $\lambda_\uparrow=1$ on $[0,1]$,
$\lambda_\uparrow=2^\beta$ on $[2,\infty)$, and, for $\theta\in\{\frac32,2\}$,
$\lambda_{\downarrow,\theta}$ nonincreasing, $\lambda_{\downarrow,\theta}=1$ on $[0,1]$, $\lambda_{\downarrow,\theta}=0$ on
$[\theta,\infty)$. Let $C_\lambda\ge1$ bound the first two derivatives of
these three functions; once the profiles are fixed as functions of
$\alpha$ (through $\beta=1-\alpha$), this constant depends only on
$\alpha$. For $a>0$ put
\[
 m_a:=\min\big(1,(a/\ell_J)^\beta\big),\qquad k_a:=2C_*C_\lambda m_a .
\]

\begin{lemma}[Exponential moment bounds]\label{lem:expmoment}
There is a constant $c_{\rm em}>0$, depending only on $d$, $\alpha$,
the drift bound $C_*$ and the fixed cutoff profiles, with the following
property. Let $0\le s_0<s_1\le t_0$, and let $\varrho$ be an energy
solution of \eqref{eq:doubled} on $[s_0,s_1]$ such that $\varrho(s_0)$
is a probability density. The same constant works in both cases (F)
and (B), independently of $J$, $\kappa$, $\mu$, the times and the data.
\begin{enumerate}[(a)]
\item (Growth.) Let $0<a\le\frac18$ satisfy
 \[
 m_a(s_1-s_0)\le c_{\rm em}a^\beta,
 \qquad \mu a^{\beta-2}\le c_{\rm em}m_a.
 \]
 If $\varrho(s_0)$ vanishes a.e.\
 on $\{R>a\}$, then
 \[
 \varrho_{s_1}\big(\{R\ge2a\}\big)\le
 \exp\Big(-c_{\rm em}\frac{m_a a^{2-\beta}}{\mu}\Big).
 \]
\item (Shrinking.) Let $\theta\in\{\frac32,2\}$ and $0<b\le\frac18$ satisfy
 \[
 m_b(s_1-s_0)\le c_{\rm em}b^\beta,
 \qquad \mu b^{\beta-2}\le c_{\rm em}m_b.
 \]
 If $\varrho(s_0)$ vanishes a.e.\ on $\{R<\theta b\}$, then
 \[
 \varrho_{s_1}\big(\{R\le b\}\big)\le
 \exp\Big(-c_{\rm em}\frac{m_b b^{2-\beta}}{\mu}\Big).
 \]
\end{enumerate}
\end{lemma}

\begin{proof}
We prove (a) by constructing an exponential test function which is
larger on the destination region than on the initial support.
Choose $c_{\rm em}$ sufficiently small, in terms of the stated
parameters, that the hypotheses in (a) imply
\[
 6k_a(s_1-s_0)\le(2^\beta-1)a^\beta,
 \qquad 2d\mu C_\lambda a^{\beta-2}\le k_a.
\]
The same choice makes the hypotheses in (b) imply the corresponding
inequalities with $a$ replaced by $b$ and $(2^\beta-1)a^\beta$
replaced by $b^\beta$.
Let $f(Y)=a^\beta\lambda_\uparrow(R(Y)/a)$. Since $f$ is constant where
$R\le a$ and where $R\ge2a$, and $2a\le\frac14$, $f$ is smooth on $\T^{2d}$
and is a function of $z$ where it is not constant. On $\{a\le R\le2a\}$,
radial differentiation in $z\in\R^d$ gives
\[
\begin{aligned}
 |\nabla_zf|&\le C_\lambda a^{\beta-1},\\
 |\Delta_zf|&\le C_\lambda a^{\beta-2}
 +(d-1)C_\lambda a^{\beta-1}/R\le dC_\lambda a^{\beta-2}.
\end{aligned}
\]
The relative drift satisfies \eqref{eq:holder}, also in case (B) by
Lemma~\ref{lem:reversed}(b). On the same annulus,
\[
 |\delta v|\le C_*\min(R^\alpha,\la_JR)
 \le C_*\min((2a)^\alpha,2\la_Ja)\le2C_*a^\alpha m_a,
\]
where we used $\la_Ja=a^\alpha(a/\ell_J)^\beta$. The choice
$\alpha+\beta=1$ therefore yields $|\delta v\cdot\nabla_zf|\le k_a$.
We compensate for this term and the two diffusion terms by setting
\[
 \mathsf M:=\frac{k_a\,a^{2-2\beta}}{2C_\lambda^2\mu},\qquad \mathsf K:=3k_a,\qquad
 U(s,Y):=\exp\big(\mathsf M(f(Y)+\mathsf K(s_1-s))\big).
\]
Since $\mathcal T_s$ acts on functions of $z$ as $\delta v\cdot\nabla_z+2\mu\Delta_z$,
\[
\begin{aligned}
 (\partial_s+\mathcal T_s)U
 &=\mathsf MU\big[-\mathsf K+\delta v\cdot\nabla_zf+2\mu\Delta_zf
 +2\mu\mathsf M|\nabla_zf|^2\big]\\
 &\le\mathsf MU\big[-3k_a+k_a+k_a+k_a\big]=0
\end{aligned}
\]
for a.e.\ $s$, using $2d\mu C_\lambda a^{\beta-2}\le k_a$ and
$2\mu\mathsf MC_\lambda^2a^{2\beta-2}=k_a$ on $\{a\le R\le2a\}$; elsewhere the
bracket equals $-\mathsf K$. By Lemma~\ref{lem:weakbasic}(c) and
$\varrho\ge0$,
\[
 \int U(s_1)\varrho_{s_1}\le\int U(s_0)\varrho_{s_0}=e^{\mathsf M(a^\beta+\mathsf K(s_1-s_0))},
\]
because $f=a^\beta$ on $\{R\le a\}$. On $\{R\ge2a\}$, $U(s_1)=e^{\mathsf M(2a)^\beta}$,
and $U>0$. Therefore
\[
\begin{aligned}
 \varrho_{s_1}(\{R\ge2a\})
 &\le\exp\big(-\mathsf M[(2^\beta-1)a^\beta-3k_a(s_1-s_0)]\big)\\
 &\le\exp\big(-\mathsf M(2^\beta-1)a^\beta/2\big),
\end{aligned}
\]
Since $k_a=2C_*C_\lambda m_a$, the last exponent is a fixed
positive multiple of $m_a a^{2-\beta}/\mu$. Decreasing $c_{\rm em}$
if necessary gives (a).

For (b), reverse the spatial cutoff while keeping the same time
correction. Let $g(Y)=b^\beta\lambda_{\downarrow,\theta}(R(Y)/b)$, which is smooth since $\theta b\le\frac14$;
$\nabla g$ is supported in $\{b\le R\le\theta b\}$, where the same bounds hold
with $a$ replaced by $b$ (as $\theta\le2$). With
$\mathsf M=k_bb^{2-2\beta}/(2C_\lambda^2\mu)$, $\mathsf K=3k_b$ and
$U=\exp(\mathsf M(g+\mathsf K(s_1-s)))$, again $(\partial_s+\mathcal T_s)U\le0$. Now
$g=0$ on $\{R\ge\theta b\}$ and $g=b^\beta$ on $\{R\le b\}$, so
\[
 \varrho_{s_1}(\{R\le b\})
 \le\exp(-\mathsf M[b^\beta-3k_b(s_1-s_0)])
 \le\exp(-\mathsf Mb^\beta/2).
\]
Again the exponent is a fixed positive multiple of
$m_b b^{2-\beta}/\mu$. Choose $c_{\rm em}$ small enough for both
parts and fix it for the rest of the paper.
\end{proof}

\subsection{Exit and transfer bounds}

We use the preceding lemma in two ways: to keep most of the mass near
the diagonal, and to compare the mass in a given separation range
with its mass at the preceding renewal. For the second application,
we partition the separations into target sets and enlarge each of
them to a set on which the corresponding gain has a lower bound.

Let $C_1':=\max(2,(16\ell_{j_0})^{-1})$. For $R\in[\ell_J,\frac1{16})$ let
$j(R)$ be the index with $\ell_j\le R<2\ell_j$ if $R\le\ell_{j_0}$, and
$j(R)=j_0$ otherwise; then $1\le R/\ell_{j(R)}\le C_1'$. Define the
\emph{target sets}
\[
\begin{aligned}
 \Sigma_{\rm core}&=\{R\le r_\sharp/2\},\\
 \Sigma_\flat&=\{r_\sharp/2<R<\ell_J\},\\
 \Sigma_j&=\{\ell_J\le R<\tfrac1{16},\ j(R)=j\},
\end{aligned}
\]
which partition $\{R<\frac1{16}\}$ (for $J\ge J_0$, $\ell_J\le\frac1{32}$),
and, with $b_j:=\min(C_1'\ell_j,\frac1{16})$, the \emph{gain sets}
\[
 \Sigma_j^+=\{\tfrac12\ell_j\le R\le2b_j\},\qquad
 \Sigma_\flat^+=\{\tfrac14r_\sharp\le R\le\tfrac32\ell_J\}.
\]
Then $\Sigma_j\subset\{\ell_j\le R\le b_j\}\subset \Sigma_j^+$, $\Sigma_\flat\subset \Sigma_\flat^+$,
$\Sigma_j^+\cup \Sigma_\flat^+\subset\{R\le\frac18\}$, and $\Sigma_j\ne\emptyset$ only if
$\ell_j<\frac1{16}$, in which case $b_j\ge\ell_j$.

\begin{corollary}[Exit and transfer]\label{cor:exit}
There exist $t_1\in(0,1]$ and positive constants
$c_E,c_T,c_\flat,\kappa_{\rm ex},\varepsilon_{\rm ex}$, depending only on
$d,\alpha,\psi,R,M,\sigma$, such
that the following holds. Let $J\ge J_0$, assume
\eqref{eq:eps0small}, and suppose
\[
 0<\kappa\le\kappa_{\rm ex},\qquad
 0<\varepsilon_0\le\varepsilon_{\rm ex}.
\]
These constants are independent of $J$, $t_0$, the realization and
the initial probability density.
\begin{enumerate}[(a)]
\item (Exit.) If $\varrho_0$ vanishes a.e.\ on $\{R>\frac1{32}\}$, then
 \[
 \varrho_s(\{R\ge\tfrac1{16}\})\le e^{-c_E/\mu},
 \qquad s\in(0,\min(t_0,t_1)].
 \]
\item (Transfer.) Let $0\le s'<s\le t_0$. For each $j$ with $\Sigma_j\ne\emptyset$,
 if $s-s'\le\tau_j$, then
 \[
 \Pi_{s',s}((\Sigma_j^+)^c;\Sigma_j)
 \le\epsilon_j:=\exp\big(-c_T(\ell_j/\ell_J)^{1+\alpha}/\varepsilon_0\big).
 \]
 If $s-s'\le\tau_J$, then
 \[
 \Pi_{s',s}((\Sigma_\flat^+)^c;\Sigma_\flat)
 \le\epsilon_\flat:=e^{-c_\flat/\kappa}+e^{-c_T/\varepsilon_0}.
 \]
\end{enumerate}
\end{corollary}

\begin{proof}
Write $c_{\rm em}$ for the constant in Lemma~\ref{lem:expmoment}.
All constants below depend only on the parameters in the statement.
Choose $t_1\in(0,1]$ so that
\[
 t_1\le c_{\rm em}32^{-\beta}.
\]
For (a), apply Lemma~\ref{lem:expmoment}(a) with $a=1/32$,
$s_0=0$ and $s_1=s$. Here $m_a=1$. The time condition follows
from $s\le t_1$, and the diffusion condition follows from
$\mu\le\varepsilon_0$ once $\varepsilon_{\rm ex}$ is sufficiently
small. The resulting exponent is a fixed positive multiple of
$1/\mu$, giving the claimed estimate for a suitable $c_E>0$.

For (b), split an initial probability density supported in
$(\Sigma_j^+)^c$ into its parts on $\{R<\ell_j/2\}$ and
$\{R>2b_j\}$. By linearity and positivity, it suffices to bound
the mass reaching $\Sigma_j$ from either part after normalization.
The two bounds are then weighted by initial masses whose sum is one.

For the first part, use Lemma~\ref{lem:expmoment}(a) with
$a=\ell_j/2$; the target lies in $\{R\ge2a\}$. Since
$a\ge\ell_J/2$, we have
\[
 m_a\ge2^{-\beta},\qquad
 \frac{m_a(s-s')}{a^\beta}\le 2^\beta\kappa,\qquad
 \frac{\mu a^{\beta-2}}{m_a}\le 4\varepsilon_0
\]
(using $s-s'\le\tau_j\le\kappa\ell_j^\beta$ and
$\mu a^{\beta-2}/m_a\le4\varepsilon_0(\ell_J/\ell_j)^{1+\alpha}$).
Thus the time and diffusion conditions hold for sufficiently small
$\kappa_{\rm ex}$ and $\varepsilon_{\rm ex}$. Moreover,
\[
 \frac{m_a a^{2-\beta}}{\mu}
 \ge c\,\frac{(\ell_j/\ell_J)^{1+\alpha}}{\varepsilon_0},
\]
using $2-\beta=1+\alpha$ and
$\mu=\varepsilon_0\ell_J^{1+\alpha}$.
For the second part, use Lemma~\ref{lem:expmoment}(b) with
$\theta=2$ and $b=b_j$. Here $b\ge\ell_j\ge\ell_J$, so $m_b=1$,
\[
 \frac{s-s'}{b^\beta}\le\kappa,\qquad
 \mu b^{\beta-2}\le\varepsilon_0,\qquad
 \frac{b^{2-\beta}}{\mu}
 \ge\frac{(\ell_j/\ell_J)^{1+\alpha}}{\varepsilon_0}.
\]
The two applications therefore give the stated $\epsilon_j$ with
one constant $c_T>0$, decreased if necessary.

For $\Sigma_\flat$, split instead at $r_\sharp/4$ and
$3\ell_J/2$. For the lower part apply
Lemma~\ref{lem:expmoment}(a) with $a=r_\sharp/4$.
Now $m_a=(a/\ell_J)^\beta$, and \eqref{eq:diffusivity} gives
\[
 \frac{m_a(s-s')}{a^\beta}\le\kappa,\qquad
 \frac{\mu a^{\beta-2}}{m_a}
 =\frac{\mu\ell_J^\beta}{a^2}
 =\frac{4\kappa}{\beta\rho^\beta},
 \qquad
 \frac{m_a a^{2-\beta}}{\mu}
 =\frac{\beta\rho^\beta}{4\kappa}.
\]
Thus sufficiently small $\kappa_{\rm ex}$ gives the hypotheses
and a bound $e^{-c_\flat/\kappa}$, with $c_\flat>0$ fixed.
For the upper part use Lemma~\ref{lem:expmoment}(b) with
$\theta=3/2$ and $b=\ell_J$. The same estimates as above give
$e^{-c_T/\varepsilon_0}$, after decreasing $c_T$ if necessary.
Adding the bounds proves the claim. The conditions used above
are all implied by
\[
 2^\beta\kappa_{\rm ex}\le c_{\rm em},\qquad
 \frac{4\kappa_{\rm ex}}{\beta\rho^\beta}\le c_{\rm em},\qquad
 32^{2-\beta}\varepsilon_{\rm ex}\le c_{\rm em}
\]
(the first implies $\kappa_{\rm ex}\le c_{\rm em}$; the last is the
diffusion condition in (a), where $\mu a^{\beta-2}\le32^{2-\beta}\varepsilon_0$
for $a=1/32$, and it implies $4\varepsilon_{\rm ex}\le c_{\rm em}$);
we fix $\kappa_{\rm ex}$ and $\varepsilon_{\rm ex}$ accordingly.
\end{proof}

The transfer errors are summable uniformly in $J$. Indeed, since
$(\ell_j/\ell_J)^{1+\alpha}\ge2^{J-j}\ge1+J-j$, if
$e^{-c_T/\varepsilon_0}\le\frac12$, then
\begin{equation}\label{eq:epstot}
 \epsilon_{\rm tot}:=\sum_j\epsilon_j+\epsilon_\flat
 \le\sum_{m\ge1}e^{-mc_T/\varepsilon_0}+\epsilon_\flat
 \le3e^{-c_T/\varepsilon_0}+e^{-c_\flat/\kappa}.
\end{equation}

\subsection{Accumulation of the gains}

We now distribute each renewal contribution over the following
window. The resulting quantity $\mathcal G(s)$ can be integrated in
time, while the transfer estimates compare it with the density at
time $s$. Diffusion supplies the gain on the quadratic core, where no
renewal contribution is needed.

Let $c_3:=\min\big(\frac d2,\frac{c_2}2,\frac14c_{1/2,\,2C_1'}\big)$, with
$c_2$ from Lemma~\ref{lem:core}(f) and $c_{a,b}$ from
\eqref{eq:comparablegain}. The functions $\Gamma^\sharp_j(R)$ are defined
for $R\le\frac14$; products with $\chi_\circ(R)$ are extended by $0$.
Recall that $s_j(s)$ is the most recent renewal time of scale $j$
at or before $s$, or $0$ if there has been no renewal since time $0$,
with the left-limit convention at $t_0$. Thus $\varrho_{s_j(s)}$ is
the density at the start of the current window, and the indicator
$\one_{\{s_j(s)>0\}}$ excludes the initial window. For $s\in(0,t_0]$ put
\[
 \mathcal G(s):=\sum_j\tau_j\,\one_{\{s_j(s)>0\}}
 \int_{\{R\le1/8\}}\Gamma^\sharp_j\,\varrho_{s_j(s)} .
\]

\begin{lemma}[Gains]\label{lem:gains}
Assume \eqref{eq:eps0small}, $L\ge2$, $J\ge J_0$ and the hypotheses of
Corollary~\ref{cor:exit}, and let $\varrho_0\in\mathfrak D_\varsigma$ with
$\varsigma\le r_\sharp/8$. For fixed data:
\begin{enumerate}[(a)]
\item for $t\in(0,t_0]$,
 
\[
 \sum_{\tilde{\zeta}\in\mathcal Z,\,\tilde{\zeta}\le t}\ \sum_{j\in\mathcal J(\tilde{\zeta})}
 \tau_j^2\!\int\chi_\circ\Gamma^\sharp_j\,\varrho_{\tilde{\zeta}}
 \ge\int_0^t\mathcal G(s)\dd s;
\]

\item for $s\in(0,t_0]$,
 
\[
 \mathcal G(s)+\int\chi_\circ\,2\mu\Delta\Phi^\sharp\,\varrho_s
 \ge c_3\kappa\big(1-\varrho_s(\{R\ge\tfrac1{16}\})-\epsilon_{\rm tot}\big).
\]

\end{enumerate}
\end{lemma}

\begin{proof}
All integrands are nonnegative (Lemma~\ref{lem:core}(e),(f)), and
$\chi_\circ=1$ on $\{R\le\frac18\}$.

For (a), fix $j$. The set of $s\in(0,t]$ with $s_j(s)>0$ is the disjoint union,
over the renewal times $\tilde{\zeta}\in(0,t]$ of scale $j$, of the sets of
$s\in(0,t]$ in the window of scale $j$ starting at $\tilde{\zeta}$; each has length
at most $\tau_j$, and $s_j(s)=\tilde{\zeta}$ on it. Hence the contribution of $j$
to the right side is at most
\[
 \sum_{\tilde{\zeta}}\tau_j^2\int_{\{R\le1/8\}}\Gamma^\sharp_j\varrho_{\tilde{\zeta}}\le
\sum_{\tilde{\zeta}}\tau_j^2\int\chi_\circ\Gamma^\sharp_j\varrho_{\tilde{\zeta}},
\]
 and a renewal time
$\tilde{\zeta}\le t$ of scale $j$ is a point of $\mathcal Z$ with $j\in\mathcal J(\tilde{\zeta})$.

To prove (b), we first apply the lower bound for each gain on its
enlarged set. On $\Sigma_j^+$ with $j<J$, $\tau_j\Gamma^\sharp_j\ge\frac12c_{1/2,2C_1'}\kappa$ by
Lemma~\ref{lem:core}(f) and \eqref{eq:comparablegain}: there
$R/\ell_j\in[\frac12,2C_1']$, $R\le\frac18$, $R\ge\frac12\ell_J>2\varepsilon$, and
$w_j(R)=1$ since $\ell_j\le2R\le LR$. On $\Sigma_J^+\cup \Sigma_\flat^+$,
$\tau_J\Gamma^\sharp_J\ge\min(\frac12c_{1/2,2C_1'},c_2)\kappa$, and
$\one_{\Sigma_J^+\cup \Sigma_\flat^+}\ge\frac12(\one_{\Sigma_J^+}+\one_{\Sigma_\flat^+})$. On $\Sigma_{\rm core}$,
$2\mu\Delta\Phi^\sharp=\frac{d\kappa}2$ by Lemma~\ref{lem:core}(e). Hence
the left side is at least
\[
 c_3\kappa\Big[\sum_j\one_{\{s_j(s)>0\}}\varrho_{s_j(s)}(\Sigma_j^+)
 +\one_{\{s_J(s)>0\}}\varrho_{s_J(s)}(\Sigma_\flat^+)+\varrho_s(\Sigma_{\rm core})\Big].
\]
Terms with $\Sigma_j=\emptyset$ are nonnegative and may be discarded.
It remains to compare these masses at renewal times with the masses
of the target sets at time $s$. More precisely, we claim that
\[
 \one_{\{s_j(s)>0\}}\varrho_{s_j(s)}(\Sigma_j^+)
 \ge\varrho_s(\Sigma_j)-\epsilon_j
\] for every
$j$ with $\Sigma_j\ne\emptyset$, and
\[
 \one_{\{s_J(s)>0\}}\varrho_{s_J(s)}(\Sigma_\flat^+)
 \ge\varrho_s(\Sigma_\flat)-\epsilon_\flat.
\] If
$s_j(s)=s$, this holds because $\Sigma_j\subset \Sigma_j^+$. If $s_j(s)<s$, then
$s-s_j(s)\le\tau_j$, and Lemma~\ref{lem:transfer} with
Corollary~\ref{cor:exit}(b) gives
$\varrho_s(\Sigma_j)\le\varrho_{s_j(s)}(\Sigma_j^+)+\epsilon_j$; if moreover $s_j(s)=0$,
then $\varrho_0(\Sigma_j^+)=0$, because $\varrho_0$ vanishes on $\{R>r_\sharp/8\}$
and $\Sigma_j^+\subset\{R\ge\ell_J/2\}$. The same argument applies to $\Sigma_\flat^+$,
since $\Sigma_\flat^+\subset\{R\ge r_\sharp/4\}$. Finally, the sets $\Sigma_{\rm core}$,
$\Sigma_\flat$ and $\Sigma_j$ partition $\{R<\frac1{16}\}$.
\end{proof}

\subsection{Proof of Theorem~\ref{thm:pair}}

We choose the parameters so that the errors in the evolution identity
are smaller than the accumulated gain. We then compare the corrected
observable with $R^\beta$ and obtain the second moment by Jensen's
inequality.

\emph{Choice of parameters.} Let $C_0,C_3$ be the constants in $|E_0|\le C_0L^{-\beta}$ and
$|E_3|\le C_3\kappa^2N_L^2$ (Proposition~\ref{prop:drift} and
Corollary~\ref{cor:coretransfer}), $C_9$ that of
Lemma~\ref{lem:diffcorrector}, $C_B$ that of Lemma~\ref{lem:cutoff}, and
$C'$ that of \eqref{eq:endpointbound}. Fix $K\ge1$ (not necessarily the
$K$ of Subsection~\ref{sec:proof}) with $C_0/K\le c_3/8$, and
$L=(K/\kappa)^{1/\beta}$, so that $|E_0|\le c_3\kappa/8$. Since
$N_L\le C(1+|\log\kappa|)$, we can choose
$\kappa_1\le\min(\kappa_*,\kappa_{\rm ex})$, with $\kappa_{\rm ex}$
from Corollary~\ref{cor:exit}, such that, for $\kappa\le\kappa_1$,
\[
 L\ge2,\qquad \kappa N_L\le1,\qquad c_3\kappa\le1,
\]
\[
 C_3\kappa N_L^2\le c_3/8,\qquad C_9\kappa N_L\le c_3/8,
 \qquad e^{-c_\flat/\kappa}\le\tfrac1{16}.
\]
We also require
\[
 \delta_\kappa:=C_\delta\kappa(1+|\log\kappa|)\le\tfrac14,
 \qquad C'(\kappa N_L+\kappa^2N_L^2)\le\delta_\kappa,
\]
where $C_\delta$ is chosen to ensure the second inequality. Let
$\varepsilon_1(\kappa)>0$ be such that $\varepsilon_0\le\varepsilon_1(\kappa)$
implies \eqref{eq:eps0small}, $\varepsilon_0\le\kappa$,
$\varepsilon_0\le\varepsilon_{\rm ex}$ and $3e^{-c_T/\varepsilon_0}\le\frac1{16}$. Then
$\epsilon_{\rm tot}\le\frac18$ by \eqref{eq:epstot}.

\emph{The evolution estimate.} Let $\varrho_0\in\mathfrak D_\varsigma$ with
$\varsigma\le r_\sharp/8$ be independent of the data, and fix
$t\in(0,\min(t_0,t_1)]$. Lemma~\ref{lem:datadependence} gives
measurability in the data. Together with continuity of $\varrho_s$ in
$L^2$ and piecewise continuity of $\mathcal K_s$ and $\mathcal G(s)$,
this gives joint measurability in $(s,\text{data})$ and allows us to
interchange expectations and time integrals.

First, bound the initial value. On the support of $\varrho_0$,
$\chi_\circ=1$ and $\Phi^\sharp\ge\Phi^\sharp(0)$ by
Lemma~\ref{lem:core}(a). The endpoint estimate
\eqref{eq:endpointbound}, extended to the diagonal by continuity, gives
\[
 |Q^\sharp-\Phi^\sharp|
 \le\delta_\kappa(R^\beta+\ell_J^\beta)\le\ell_J^\beta,
 \qquad
 \int Q^\circ_0\varrho_0\ge\Phi^\sharp(0)-\ell_J^\beta.
\]

Next, separate the positive diffusion term from the errors in
$\mathcal K_s$. Corollary~\ref{cor:coretransfer} and
Lemma~\ref{lem:diffcorrector} give
\[
 \chi_\circ|E_0+E_3+\mu\Delta_Y(Q^\sharp-\Phi^\sharp)|
 \le\tfrac38c_3\kappa,
 \qquad
 |\mathcal K_s^{\rm cut}|
 \le C_B\one_{\{R\ge1/16\}}.
\]
Since $\int\varrho_s=1$, these cost at most
$\tfrac38c_3\kappa t$ and
$C_B\int_0^t\E\varrho_s(\{R\ge1/16\})\dd s$, respectively.
Take expectations in \eqref{eq:evolution}, use
\eqref{eq:jumpaverage}, and apply Lemma~\ref{lem:gains}(a) to
replace the sum of renewal gains by $\int_0^t\mathcal G(s)\dd s$.
We obtain
\[
\begin{aligned}
 \E\!\int Q^\circ_t\varrho_t
 &\ge\Phi^\sharp(0)-\ell_J^\beta
 +\int_0^t\E\Big[\mathcal G(s)
       +\int\chi_\circ\,2\mu\Delta\Phi^\sharp\,\varrho_s\Big]\dd s\\
 &\quad-\tfrac38c_3\kappa t
       -C_B\int_0^t\E\varrho_s(\{R\ge\tfrac1{16}\})\dd s.
\end{aligned}
\]

Now bound the combined renewal and diffusion gains.
Lemma~\ref{lem:gains}(b), $\epsilon_{\rm tot}\le1/8$, and the exit
estimate in Corollary~\ref{cor:exit}(a) give, for every $s\le t$,
\[
\begin{aligned}
 \E\Big[\mathcal G(s)
       +\int\chi_\circ\,2\mu\Delta\Phi^\sharp\,\varrho_s\Big]
 &\ge c_3\kappa\big[1-\E\varrho_s(\{R\ge\tfrac1{16}\})
                       -\epsilon_{\rm tot}\big]\\
 &\ge c_3\kappa\big(\tfrac78-e^{-c_E/\mu}\big).
\end{aligned}
\]
The same exit estimate bounds the cutoff error by
$C_Bte^{-c_E/\mu}$. Substituting these two bounds above leaves
\[
\begin{aligned}
 \E\!\int Q^\circ_t\varrho_t
 &\ge\Phi^\sharp(0)-\ell_J^\beta
   +\big(\tfrac78-\tfrac38\big)c_3\kappa t
   -(c_3\kappa+C_B)t e^{-c_E/\mu}\\
 &\ge\Phi^\sharp(0)-\ell_J^\beta
   +\tfrac12c_3\kappa t-(1+C_B)e^{-c_E/\mu},
\end{aligned}
\]
where the last inequality uses $t\le1$ and $c_3\kappa\le1$.

\emph{Return to the separation moments.} We now compare the corrected
observable with $R^\beta$. Where $Q^\circ>0$, we have $Q^\sharp>0$
and $0\le\chi_\circ\le1$, so Lemma~\ref{lem:core}(a) and
\eqref{eq:endpointbound} give
\[
\begin{aligned}
 Q^\circ
 &\le Q^\sharp
 \le\Phi^\sharp+|Q^\sharp-\Phi^\sharp|\\
 &\le R^\beta+\ell_J^\beta
       +\delta_\kappa(R^\beta+\ell_J^\beta)\\
 &\le(1+\delta_\kappa)R^\beta+2\ell_J^\beta.
\end{aligned}
\]
Where $Q^\circ\le0$, the last bound is immediate. Integrating it
against the probability density $\varrho_t$ and using the evolution
estimate therefore yields
\[
\begin{aligned}
 (1+\delta_\kappa)\E\!\int R^\beta\varrho_t
 &\ge\E\!\int Q^\circ_t\varrho_t-2\ell_J^\beta\\
 &\ge\tfrac12c_3\kappa t+\Phi^\sharp(0)-3\ell_J^\beta
       -(1+C_B)e^{-c_E/\mu}.
\end{aligned}
\]
Thus the remaining endpoint error comes from $\Phi^\sharp(0)$.
Lemma~\ref{lem:core}(c) and $\varepsilon\le\ell_J$ give
\[
 \Phi^\sharp(0)
 \ge-\beta\varepsilon^\beta
       \big(1+4^\beta\log(2\varepsilon/r_\sharp)\big)
 \ge-\big(\beta+\beta4^\beta\log(2\varepsilon/r_\sharp)\big)\ell_J^\beta.
\]
Substituting this bound above, we obtain
\[
\begin{aligned}
 (1+\delta_\kappa)\E\!\int R^\beta\varrho_t
 &\ge\tfrac12c_3\kappa t-(1+C_B)e^{-c_E/\mu}\\
 &\quad-\big(3+\beta+\beta4^\beta\log(2\varepsilon/r_\sharp)\big)\ell_J^\beta.
\end{aligned}
\]
This is \eqref{eq:pairmain} with $c_5=c_3/2$, after choosing the
constant in $\mathrm{Err}_J$ large enough.

If $\mathrm{Err}_J\le c_5\kappa t/2$, then $\delta_\kappa\le1/4$ gives
\[
 \E\!\int R^\beta\varrho_t
 \ge\frac{c_5\kappa t-\mathrm{Err}_J}{1+\delta_\kappa}
 \ge\frac45\,\frac{c_5\kappa t}{2}
 =\frac25c_5\kappa t.
\]
Finally, $R^2=(R^\beta)^{2/\beta}$ and $2/\beta>1$. Jensen's
inequality, first against $\varrho_t$ for each realization and then
against the probability law of the field, gives
\[
\begin{aligned}
 \E\!\int R^2\varrho_t
 &\ge\E\Big[\Big(\int R^\beta\varrho_t\Big)^{2/\beta}\Big]\\
 &\ge\Big(\E\!\int R^\beta\varrho_t\Big)^{2/\beta}
 \ge\Big(\frac25c_5\kappa t\Big)^{2/\beta}.
\end{aligned}
\]
This proves both bounds in \eqref{eq:pairmoments}. \qed

\section{Anomalous dissipation}\label{sec:AD}

We convert the separation estimate into a lower bound for scalar
dissipation. The first step is a duality identity expressing the
dissipation through the reversed two-point equation. Translation
invariance then diagonalizes the mean dissipation in Fourier space,
and the remaining symmetries turn the separation bound into a lower
bound for each nonzero mode. A Fourier truncation gives the result
for general initial data.

The last subsection passes from the truncated fields to the single
limiting field. This uses a deterministic energy comparison and does
not require a new dispersion argument.

\subsection{The dissipation as a two-point quantity}

To approximate two particles starting at the same point, we take
their initial separation to be less than $\varsigma$, while keeping each
marginal equal to Lebesgue measure. For $0<\varsigma\le\frac14$ let
\[
 \varrho^\varsigma_0(y,y'):=|B_\varsigma|^{-1}\one_{\{\dist(y,y')<\varsigma\}},\qquad
 \mathcal M_\varsigma f(y):=\int_{\T^d}f(y')\varrho^\varsigma_0(y,y')\dd y',
\]
where $|B_\varsigma|$ is the volume of the Euclidean ball of radius $\varsigma$ in
$\R^d$. Since $\varsigma\le\frac14<\frac12$, the set $\{y':\dist(y,y')<\varsigma\}$ is
isometric to that ball for every $y$, so both marginals of $\varrho^\varsigma_0$
are Lebesgue measure on $\T^d$; thus $\varrho^\varsigma_0\in\mathfrak D_\varsigma$, and $\mathcal M_\varsigma$ is
the average over balls of radius $\varsigma$. As $\varsigma\downarrow0$,
$\mathcal M_\varsigma f\to f$ in $L^2(\T^d)$ for every $f\in L^2(\T^d)$.

\begin{lemma}[Dissipation identity]\label{lem:pairidentity}
Fix $J$, $t_0\in(0,1]$ and the data. For $0<\varsigma\le\frac14$ let $\varrho^\varsigma$
be the energy solution of \eqref{eq:doubled} in case (B) (the reversed
field with final time $t_0$) with $\varrho^\varsigma(0)=\varrho^\varsigma_0$. Then for
every real $\theta_0\in L^\infty(\T^d)$, with $\theta=\theta^{(J)}$,
\begin{equation}\label{eq:pairidentity}
 \frac12\int\!\!\int\big(\theta_0(y)-\theta_0(y')\big)^2\varrho^\varsigma_{t_0}(y,y')\dd y\dd y'
 =\|\theta_0\|^2-\big(\theta(t_0),\mathcal M_\varsigma\theta(t_0)\big),
\end{equation}
the left side lies in $[0,2\|\theta_0\|_\infty^2]$, and consequently
\[
 \mathbf D_J(\theta_0)=\lim_{\varsigma\downarrow0}\frac12\int\!\!\int
 \big(\theta_0(y)-\theta_0(y')\big)^2\varrho^\varsigma_{t_0}(y,y')\dd y\dd y' .
\]
\end{lemma}

\begin{proof}
Duality moves the reversed evolution from the density to the scalar
increments. The product formula then reduces the result to the
one-point scalar equation.

Let 
\[
\begin{aligned}
 F(y,y')&:=(\theta_0(y)-\theta_0(y'))^2\\
 &=\theta_0^2\otimes1+1\otimes\theta_0^2
 -2\,\theta_0\otimes\theta_0\in L^2(\T^{2d}),
\end{aligned}
\]
 and let $w$ be the energy
solution of \eqref{eq:doubled} in case (F) on $[0,t_0]$ with $w(0)=F$. The
drift of case (B) is $-\mathbf u(t_0-s)$ for all but finitely many $s$,
where $\mathbf u$ is the drift of case (F)
(Definition~\ref{def:reversed}). By Lemma~\ref{lem:duality},
\[
 \int F\varrho^\varsigma_{t_0}=\int w(t_0)\varrho^\varsigma_0 .
\]
For $g\in L^2(\T^d)$ let $\mathcal P_{0,t_0} g$ be the value at time $t_0$ of the energy
solution of \eqref{eq:scalar} with initial value $g$ at time $0$.
Lemma~\ref{lem:products}(a) and linearity give
\[
\begin{aligned}
 w(t_0,y,y')&=\mathcal P_{0,t_0}(\theta_0^2)(y)+\mathcal P_{0,t_0}(\theta_0^2)(y')\\
 &\quad-2\mathcal P_{0,t_0}\theta_0(y)\mathcal P_{0,t_0}\theta_0(y').
\end{aligned}
\]
We integrate these three terms against $\varrho^\varsigma_0$ separately.
For the first term, the marginal identity and conservation of the
integral (Lemma~\ref{lem:weakbasic}(b)) give
\[
\begin{aligned}
 \iint \mathcal P_{0,t_0}(\theta_0^2)(y)\varrho^\varsigma_0(y,y')\dd y'\dd y
 &=\int \mathcal P_{0,t_0}(\theta_0^2)(y)
       \Big(\int\varrho^\varsigma_0(y,y')\dd y'\Big)\dd y\\
 &=\int \mathcal P_{0,t_0}(\theta_0^2)(y)\dd y
 =\int\theta_0^2(y)\dd y=\|\theta_0\|^2.
\end{aligned}
\]
The second term gives the same quantity, using the other marginal.
Here $\mathcal P_{0,t_0}(\theta_0^2)$ denotes the evolution of the initial square,
whose integral is conserved. For the cross term, the definition of
$\mathcal M_\varsigma$ yields
\[
\begin{aligned}
 &\iint \mathcal P_{0,t_0}\theta_0(y)\mathcal P_{0,t_0}\theta_0(y')\varrho^\varsigma_0(y,y')\dd y'\dd y\\
 &\qquad=\int \mathcal P_{0,t_0}\theta_0(y)
       \Big(\int \mathcal P_{0,t_0}\theta_0(y')\varrho^\varsigma_0(y,y')\dd y'\Big)\dd y\\
 &\qquad=\big(\mathcal P_{0,t_0}\theta_0,\mathcal M_\varsigma \mathcal P_{0,t_0}\theta_0\big).
\end{aligned}
\]
Adding the three contributions, we obtain
\[
 \int w(t_0)\varrho^\varsigma_0
 =2\|\theta_0\|^2-2\big(\mathcal P_{0,t_0}\theta_0,\mathcal M_\varsigma \mathcal P_{0,t_0}\theta_0\big),
\]
which is \eqref{eq:pairidentity}, as $\mathcal P_{0,t_0}\theta_0=\theta(t_0)$. The bounds
follow from $0\le F\le4\|\theta_0\|_\infty^2$ and the fact that $\varrho^\varsigma_{t_0}$ is a
probability density. Finally, $\mathcal M_\varsigma\theta(t_0)\to\theta(t_0)$
in $L^2$ as $\varsigma\downarrow0$, so
\[
\begin{aligned}
 \lim_{\varsigma\downarrow0}\frac12\iint
 (\theta_0(y)-\theta_0(y'))^2\varrho^\varsigma_{t_0}(y,y')\dd y\dd y'
 &=\|\theta_0\|^2-\|\theta(t_0)\|^2\\
 &=2\mu_J\int_0^{t_0}\|\nabla\theta(t)\|^2\dd t
 =\mathbf D_J(\theta_0),
\end{aligned}
\]
where the second equality is the energy identity \eqref{eq:energy}.
\end{proof}

\subsection{Symmetries and the mean dissipation of single modes}

We distinguish two uses of symmetry. Translations eliminate the
interaction between distinct Fourier modes. Signed permutations of
the coordinates will subsequently allow us to estimate a mode in
terms of the full separation, rather than its projection onto a
single direction.

Let $\mathsf G_d$ be the group of signed permutation matrices. For
$g\in\mathsf G_d$ and $a\in\R^d$ let $\iota(x)=gx+a$, an isometry of $\T^d$
(since $g\mathbb Z^d=\mathbb Z^d$), and $\boldsymbol\iota(y,y')=(\iota y,\iota y')$. For a
vector field $v$ put $(\iota_*v)(x,t):=g\,v(\iota^{-1}x,t)$. Write
$e_k(x)=e^{2\pi\mathrm ik\cdot x}$ and
$\widehat\theta_0(k)=\int\theta_0\overline{e_k}$.

\begin{lemma}[Symmetries of the law]\label{lem:symmetry}
Assume (S). For every $\iota$ as above there is a bijection $\mathsf T_\iota$ of
the data space, preserving its law, such that
$u_{\mathsf T_\iota\omega}=\iota_*u_\omega$ for every $J$; then also
$\overleftarrow u_{\mathsf T_\iota\omega}=\iota_*\overleftarrow u_\omega$ for every
final time. Consequently the solutions $\varrho^\varsigma$ of
Lemma~\ref{lem:pairidentity}, for any final time, satisfy
$\varrho^\varsigma_s[\mathsf T_\iota\omega]=\varrho^\varsigma_s[\omega]\circ\boldsymbol\iota^{-1}$
for all $s$.
\end{lemma}

\begin{proof}
We construct the transformation of the random data explicitly,
first in one renewal window and then across all scales. Its action
on the wavelets determines the necessary relabelling and signs.

Since $\psi$ is radial, $\nabla(\psi\circ g^{-1})=g\nabla\psi\circ g^{-1}$,
so $gW^{pq}(g^{-1}y)=(gJ^{pq}g^{\mathsf T})\nabla\psi(y)$. If $ge_p=s_pe_{\pi(p)}$
with $s_p=\pm1$, then
\[
\begin{aligned}
 gJ^{pq}g^{\mathsf T}
 &=s_ps_q(e_{\pi(q)}e_{\pi(p)}^{\mathsf T}-e_{\pi(p)}e_{\pi(q)}^{\mathsf T})\\
 &=\epsilon(g,\nu)J^{\pi_g\nu}
\end{aligned}
\]
 with $\epsilon(g,\nu)=\pm1$ and $\pi_g$ a
permutation of $\mathcal A$. Periodization commutes with $g$, so
$gW_j^\nu(g^{-1}y)=\epsilon(g,\nu)W_j^{\pi_g\nu}(y)$.

Consider one window of one scale $j$, with shift $\xi$ and
coefficient family $\boldsymbol\eta=(\eta_i^\nu)_{i\in I_j,\nu\in\mathcal A}$.
Define $\xi'\in Q_j$ and $n(\xi)\in\mathbb Z^d$ by $g\xi+a=\xi'+\ell_jn(\xi)$.
Since $g$ preserves the lattice $\ell_j\mathbb Z^d$, there is a permutation
$\sigma_\xi$ of $I_j$ with $g(y_i^j+\xi)+a=y^j_{\sigma_\xi(i)}+\xi'$ in $\T^d$;
it depends on $\xi$ only through $n(\xi)$, so it is piecewise constant and
Borel in $\xi$. Set $\eta'_{\sigma_\xi(i),\pi_g\nu}:=\epsilon(g,\nu)\eta_i^\nu$. The map
$\xi\mapsto\xi'$ is induced by an affine isometry of $\R^d/\ell_j\mathbb Z^d$,
so it preserves the uniform law on $Q_j$; for fixed $\xi$, the map
$\boldsymbol\eta\mapsto\boldsymbol\eta'$ is a signed permutation of an
independent, identically distributed and symmetric family, so it
preserves the law of $\boldsymbol\eta$. Since $\xi$ and $\boldsymbol\eta$ are
independent, Fubini's theorem gives, for bounded Borel $f$,
\[
 \E f(\xi',\boldsymbol\eta')=\E_\xi\E_{\boldsymbol\eta}f(\xi',\boldsymbol\eta)=\E f(\xi,\boldsymbol\eta).
\]
 Doing this
independently in every window defines $\mathsf T_\iota$, a law-preserving
bijection (its inverse is constructed from $\iota^{-1}$ in the same way).

We next verify that this transformation has the asserted effect on
the velocity. Proceed by induction on $j$. If $\iota_*u_{<j}=u'_{<j}$ (primes refer
to $\mathsf T_\iota\omega$), then $\iota P$ solves
$\frac{d}{dt}\iota P=(\iota_*u_{<j})(\iota P)$ and starts at
$\iota(y_i^j+\xi)=y^j_{\sigma_\xi(i)}+\xi'$, so
$\iota P^j_i=P'^j_{\sigma_\xi(i)}$; since
$\iota^{-1}x-P^j_i=g^{-1}(x-\iota P^j_i)$, the mode identity gives
$g\,u_j(\iota^{-1}x)=u'_j(x)$. The reversed field is obtained from $u$ by
a fixed operation in time, which commutes with $\iota_*$.

Finally, the transformation passes to the solutions by uniqueness.
The two-point drift of $\iota_*v$ is
$\boldsymbol\iota_*\mathbf v$, with the orthogonal matrix
$\mathrm{diag}(g,g)$, which preserves $\mathbb Z^{2d}$. By
Lemma~\ref{lem:products}(b), $\varrho^\varsigma[\omega]\circ\boldsymbol\iota^{-1}$ is the
energy solution with drift $\mathbf v_{\mathsf T_\iota\omega}$ and initial
value $\varrho^\varsigma_0\circ\boldsymbol\iota^{-1}=\varrho^\varsigma_0$; by uniqueness it equals
$\varrho^\varsigma[\mathsf T_\iota\omega]$.
\end{proof}

For $t\in(0,1]$ let
$\mathbf D^{(t)}_J(\theta_0):=2\mu_J\int_0^t\|\nabla\theta^{(J)}(s)\|^2\dd s$
be the dissipation up to time $t$. For $0<\varsigma\le1/4$, let
$\varrho^{\varsigma,t}$ denote the solution of Lemma~\ref{lem:pairidentity}
with final time $t$. More explicitly, $\varrho^{\varsigma,t}_s(y,y')$,
$0\le s\le t$, is the unique energy solution of the two-point
advection--diffusion equation with the velocity reversed from time $t$:
\[
\begin{cases}
 \partial_s\varrho^{\varsigma,t}
 -u^{(J)}(y,t-s)\cdot\nabla_y\varrho^{\varsigma,t}
 -u^{(J)}(y',t-s)\cdot\nabla_{y'}\varrho^{\varsigma,t}
 =\mu_J(\Delta_y+\Delta_{y'})\varrho^{\varsigma,t},\\[2pt]
 \displaystyle
 \varrho^{\varsigma,t}_0(y,y')
 =\varrho^\varsigma_0(y,y')
 =|B_\varsigma|^{-1}\one_{\{\dist(y,y')<\varsigma\}}.
\end{cases}
\]
The equation holds on $(0,t)\times\T^{2d}$; the reversed drift is
understood with the convention of Definition~\ref{def:reversed} at
renewal times. Thus the superscript $t$ specifies the final time
used to reverse the velocity, while the subscript $s$ is the
elapsed time in this reversed evolution.

The following proposition is a standard consequence of translation
invariance. By linearity of the scalar equation, the expected
dissipation is a quadratic form in the initial datum. Translation
invariance of the velocity law makes this form diagonal in the
Fourier basis: each mode contributes its initial energy multiplied
by its expected dissipated fraction, with no cross terms between
distinct modes. The energy identity places these fractions in
$[0,1]$, and Lemma~\ref{lem:pairidentity} expresses them through the
two-point density. Thus the result uses our construction only
through translation invariance of its law; neither coordinate
symmetries nor the Richardson estimates are needed here. The
field-specific estimates enter afterwards, when we prove positive
lower bounds for the dissipated fractions.

\begin{proposition}[Mean dissipation mode by mode]\label{prop:spectral}
Assume (S). Fix $J$ and $t\in(0,1]$. For every
$k\in\mathbb Z^d\setminus\{0\}$ the limit
\[
 m^{(t)}_J(k):=\lim_{\varsigma\downarrow0}\,2\,\E\!\int\sin^2\!\big(\pi k\cdot(y-y')\big)
 \varrho^{\varsigma,t}_t(y,y')\dd y\dd y'
\]
exists, $m^{(t)}_J(-k)=m^{(t)}_J(k)\in[0,1]$, and for every real
$\theta_0\in L^2(\T^d)$ with $\int\theta_0=0$,
\[
 \E\,\mathbf D^{(t)}_J(\theta_0)=\sum_{k\ne0}|\widehat\theta_0(k)|^2m^{(t)}_J(k).
\]
\end{proposition}

\begin{proof}
We first diagonalize the dissipation for trigonometric polynomials.
The existence of the coefficient for each mode then follows by
testing a single cosine, and $L^2$ contraction allows passage to
general data.

Let $\theta_0=\sum_{0<|k|\le K}a_ke_k$ be a real trigonometric polynomial.
For fixed data $\omega$, put
\[
 I_\varsigma(\omega):=\frac12\iint
 |\theta_0(y)-\theta_0(y')|^2
 \varrho^{\varsigma,t}_t[\omega](y,y')\dd y\dd y'.
\]
Lemma~\ref{lem:pairidentity}, with final time $t$, gives
$I_\varsigma(\omega)\to\mathbf D^{(t)}_J(\theta_0)(\omega)$ as
$\varsigma\downarrow0$. Since $\varrho^{\varsigma,t}_t$ is nonnegative and has
integral $1$,
\[
 0\le I_\varsigma(\omega)\le2\|\theta_0\|_\infty^2
 \qquad\text{for every $\varsigma$ and $\omega$}.
\]
Since $\omega\mapsto I_\varsigma(\omega)$ is Borel
(Lemma~\ref{lem:datadependence}), dominated convergence with respect to the law of the data therefore
allows us to interchange the limit and expectation:
\[
 \E\,\mathbf D^{(t)}_J(\theta_0)
 =\lim_{\varsigma\downarrow0}\E I_\varsigma.
\]
Next expand the finite Fourier sum. Writing
\[
 F_{kk'}(y,y'):=(e_k(y)-e_k(y'))
               \overline{(e_{k'}(y)-e_{k'}(y'))},
\]
we have $|F_{kk'}|\le4$ and
\[
\begin{aligned}
 |\theta_0(y)-\theta_0(y')|^2
 &=\left|\sum_k a_k(e_k(y)-e_k(y'))\right|^2\\
 &=\sum_{k,k'}a_k\overline{a_{k'}}F_{kk'}(y,y').
\end{aligned}
\]
The sum is finite, so integration and expectation can be taken term
by term. Thus, setting
$\mathcal G^\varsigma_{kk'}:=\E\!\int F_{kk'}\varrho^{\varsigma,t}_t$, we obtain
\[
\begin{aligned}
 \E I_\varsigma&=\frac12\sum_{k,k'}a_k\overline{a_{k'}}\mathcal G^\varsigma_{kk'},\\
 \E\,\mathbf D^{(t)}_J(\theta_0)
 &=\lim_{\varsigma\downarrow0}\frac12\sum_{k,k'}
 a_k\overline{a_{k'}}\mathcal G^\varsigma_{kk'}.
\end{aligned}
\]
At this stage we use convergence of the whole sum, not of its
individual terms.

By Lemma~\ref{lem:symmetry} with $g=I$ and $\iota(x)=x+a$, and since $\mathsf T_\iota$
preserves the law, $\E\int F\varrho^{\varsigma,t}_t=\E\int(F\circ\boldsymbol\iota)\varrho^{\varsigma,t}_t$
for bounded Borel $F$. As $F_{kk'}\circ\boldsymbol\iota=e^{2\pi\mathrm i(k-k')\cdot a}F_{kk'}$,
it follows that $\mathcal G^\varsigma_{kk'}=0$ for $k\ne k'$. Since
$|e_k(y)-e_k(y')|=2|\sin(\pi k\cdot(y-y'))|$,
$\frac12\mathcal G^\varsigma_{kk}=2\E\int\sin^2(\pi k\cdot(y-y'))\varrho^{\varsigma,t}_t=:M^\varsigma(k)$,
and $M^\varsigma(-k)=M^\varsigma(k)$. Hence
$\E\,\mathbf D^{(t)}_J(\theta_0)=\lim_\varsigma\sum_k|a_k|^2M^\varsigma(k)$. For
$\theta_0=\cos(2\pi k\cdot x)$ this reads
$2\E\,\mathbf D^{(t)}_J(\theta_0)=\lim_\varsigma M^\varsigma(k)$, so the limit
$m^{(t)}_J(k)$ exists and lies in $[0,1]$ (as $\mathbf D^{(t)}_J(\theta_0)\le\frac12$),
and the formula holds for trigonometric polynomials.

For general $\theta_0$ let $\theta_K$ be its Fourier truncations. Since the
solution map of \eqref{eq:scalar} is linear and contracts $L^2$
(Lemma~\ref{lem:weakbasic}(a)), \eqref{eq:energy} gives
$|\mathbf D^{(t)}_J(\theta_K)-\mathbf D^{(t)}_J(\theta_0)|\le4\|\theta_K-\theta_0\|\,\|\theta_0\|$
for every realization; let $K\to\infty$ and use monotone convergence on
the right side.
\end{proof}

The preceding formula reduces the problem to a lower bound for
$m^{(t)}_J(k)$. For modes up to frequency $\mathfrak K$, we work at a
time when the separation is still at most a small multiple of
$\mathfrak K^{-1}$, except for an exponentially small tail. In this
region the sine in the dissipation formula can be bounded below by
its argument.

\begin{proposition}[The modes up to a given frequency]\label{prop:modebound}
Let $d\ge2$ and assume \eqref{eq:ND}, \eqref{eq:coefficients} and (S).
Let $\kappa\in(0,\kappa_1]$ and $\varepsilon_0\in(0,\varepsilon_1(\kappa)]$, with
$\kappa_1,\varepsilon_1,c_5$ from Theorem~\ref{thm:pair} and $t_1$ from
Corollary~\ref{cor:exit}. Retain $c_{\rm em}$ from
Lemma~\ref{lem:expmoment}, and let $t_0\in(0,1]$. For $\mathfrak K\ge1$ put
$r_{\mathfrak K}:=1/(8\mathfrak K)$ and
\[
 t_{\mathfrak K}:=\min\big(t_0,\ t_1,\ c_{\rm em}(r_{\mathfrak K}/2)^\beta\big).
\]
Let $J_1(\mathfrak K)>J_0:=\max(j_0,5)$ be such that for all $J\ge J_1(\mathfrak K)$:
$\ell_J\le r_{\mathfrak K}/2$, $\mu_J\le c_{\rm em}(r_{\mathfrak K}/2)^{2-\beta}$,
$\mathrm{Err}_J\le c_5\kappa t_{\mathfrak K}/2$, and
\[
 \Big(\frac{\sqrt d}2\Big)^\beta\exp\Big(-c_{\rm em}\frac{(r_{\mathfrak K}/2)^{2-\beta}}{\mu_J}\Big)
 \le\frac{c_5\kappa t_{\mathfrak K}}5
\]
(all these quantities tend to $0$ as $J\to\infty$). Then for
$J\ge J_1(\mathfrak K)$ and $0<|k|\le\mathfrak K$,
\[
 m^{(t_{\mathfrak K})}_J(k)\ge\frac{8|k|^2}{d}\Big(\frac{c_5\kappa t_{\mathfrak K}}5\Big)^{2/\beta}.
\]
\end{proposition}

\begin{proof}
Let $0<\varsigma\le r_\sharp/8$, let $\varrho:=\varrho^{\varsigma,t_{\mathfrak K}}_{t_{\mathfrak K}}$
(reversed field with final time $t_{\mathfrak K}$), and
$\mathcal V:=\{R\le r_{\mathfrak K}\}$.

We first localize the separation moment to $\mathcal V$. Apply
Lemma~\ref{lem:expmoment}(a) with
\[
 a=r_{\mathfrak K}/2\le\tfrac1{16},\qquad s_0=0,
 \qquad s_1=t_{\mathfrak K}.
\]
The initial density is supported in $\{R\le a\}$, since
$\varsigma\le\ell_J\le a$. Moreover, $m_a=1$.
The required time and diffusion conditions are precisely
\[
 t_{\mathfrak K}\le c_{\rm em}(r_{\mathfrak K}/2)^\beta,
 \qquad \mu_J\le c_{\rm em}a^{2-\beta},
\]
which hold by our choices of $t_{\mathfrak K}$ and $J_1(\mathfrak K)$.
We obtain
\[
 \varrho(\mathcal V^c)
 \le\exp\big(-c_{\rm em}a^{2-\beta}/\mu_J\big).
\]
On the other hand, \eqref{eq:pairmoments} in case (B), with final time
$t_{\mathfrak K}\le t_1$, gives
\[
 \E\int R^\beta\varrho\ge\tfrac25c_5\kappa t_{\mathfrak K}.
\]
Since $R\le\sqrt d/2$, the tail condition in the proposition implies
\[
 \E\int_{\mathcal V}R^\beta\varrho
 \ge\tfrac15c_5\kappa t_{\mathfrak K}.
\]
Apply Jensen's inequality to $R^\beta\one_{\mathcal V}$, first with
respect to the probability density $\varrho$ and then with respect
to $\E$. It follows that
\[
 \E\int_{\mathcal V}R^2\varrho
 \ge\Big(\E\int_{\mathcal V}R^\beta\varrho\Big)^{2/\beta}.
\]

We now use the signed permutation symmetry to recover the full
second moment from its projection in the direction of $k$. For $Y\in\mathcal V$ write $z=z(Y)$, so $|z|=R\le\frac18$,
and note $z(\boldsymbol\iota Y)=gz(Y)$ when $\iota(x)=gx$. By
Lemma~\ref{lem:symmetry} with $\iota(x)=gx$, since $\mathsf T_\iota$ preserves the
law and $\mathcal V$ is invariant under $\boldsymbol\iota$, for every bounded Borel
$\phi$ and $g\in\mathsf G_d$,
\[
 \E\int_{\mathcal V}\phi(z)\varrho
 =\E\int_{\mathcal V}\phi(gz)\varrho.
\]
With $\langle\cdot\rangle$ the average over $\mathsf G_d$,
\[
 \E\!\int\sin^2(\pi k\cdot(y-y'))\varrho
 \ge\E\!\int_{\mathcal V}\big\langle\sin^2(\pi(g^{\mathsf T}k)\cdot z)\big\rangle\varrho .
\]
On $\mathcal V$, $|(g^{\mathsf T}k)\cdot z|\le\mathfrak Kr_{\mathfrak K}=\frac18$, and
$\sin^2(\pi s)\ge4s^2$ for $|s|\le\frac12$, so the average is at least
\[
 4z^{\mathsf T}\langle g^{\mathsf T}kk^{\mathsf T}g\rangle z
 =4|k|^2|z|^2/d.
\]
Indeed, sign changes remove the off-diagonal entries of the averaged
matrix, and permutations equalize its diagonal entries. Its trace is
$|k|^2$, so the matrix is $|k|^2I/d$. Consequently,
\[
 2\E\int\sin^2(\pi k\cdot(y-y'))\varrho
 \ge\frac{8|k|^2}d\E\int_{\mathcal V}R^2\varrho.
\]
This holds for every $\varsigma\le r_\sharp/8$. Letting $\varsigma\downarrow0$
and using the localized second moment proves the proposition.
\end{proof}

\begin{corollary}[Every mode at the final time]\label{cor:everymode}
Under the hypotheses of Proposition~\ref{prop:modebound}, for every
$k\in\mathbb Z^d\setminus\{0\}$ and every $J\ge J_1(|k|)$,
\[
 m^{(t_0)}_J(k)\ge c\,\big(\kappa\min(t_0,t_1)\big)^{2/\beta},\qquad
 c=\frac8d\Big(\frac{c_5}5\Big)^{2/\beta}.
\]
\end{corollary}

\begin{proof}
Dissipation is nondecreasing in time, so it suffices to use the
estimate at $t_{|k|}$. Indeed, the solution of \eqref{eq:scalar} on $[0,t]$ is the
restriction of the solution on $[0,t_0]$, so
$t\mapsto\mathbf D^{(t)}_J(\theta_0)$ is nondecreasing. By
Proposition~\ref{prop:spectral} with $\theta_0=\cos(2\pi k\cdot x)$,
$m^{(t)}_J(k)=2\E\,\mathbf D^{(t)}_J(\theta_0)$; hence
$m^{(t_0)}_J(k)\ge m^{(t_{|k|})}_J(k)$.

Apply Proposition~\ref{prop:modebound} with $\mathfrak K=|k|$.
If $t_{|k|}=\min(t_0,t_1)$, use $|k|^2\ge1$. Otherwise
$t_{|k|}=c_{\rm em}(16|k|)^{-\beta}$, so
\[
 |k|^2t_{|k|}^{2/\beta}
 =(c_{\rm em}16^{-\beta})^{2/\beta}\ge t_1^{2/\beta},
\]
by the choice $t_1\le c_{\rm em}32^{-\beta}$ in the proof of
Corollary~\ref{cor:exit}.
\end{proof}

\subsection{Proof of Theorem~\ref{thm:ADgeneral}}

The lower bound is now uniform over any fixed finite set of nonzero
modes. We first pass to arbitrary $L^2$ data, then obtain the
probability statement, and finally make the cutoff uniform under a
Sobolev bound.

Let $\kappa_1$, $\varepsilon_1(\kappa)$ and $c_5$ be as in
Theorem~\ref{thm:pair}, $t_1$ as in Corollary~\ref{cor:exit}, and
$c_{\rm mode}:=\frac8d(\frac{c_5}5)^{2/\beta}$, the constant in
Corollary~\ref{cor:everymode}. Set $c:=c_{\rm mode}t_1^{2/\beta}$.
Since $t_0,t_1\in(0,1]$, we have $\min(t_0,t_1)\ge t_1t_0$, and hence
\[
 c_{\rm mode}\big(\kappa\min(t_0,t_1)\big)^{2/\beta}
 \ge c(\kappa t_0)^{2/\beta}=c_{\rm AD}.
\]
Thus the auxiliary time $t_1$ is absorbed into the model-dependent
constant $c$. All constants have the dependencies listed in the theorem. Fix
$\kappa\in(0,\kappa_1]$, $\varepsilon_0\in(0,\varepsilon_1(\kappa)]$,
$t_0\in(0,1]$ and a real $\theta_0\in L^2(\T^d)$ with $\int\theta_0=0$.

\emph{Part (a): general initial data.} Let $J_2(\mathfrak K):=\max\{J_1(|k|):0<|k|\le\mathfrak K\}$. By
Proposition~\ref{prop:spectral} at $t=t_0$, $m^{(t_0)}_J\ge0$ and
Corollary~\ref{cor:everymode}, for every $\mathfrak K\ge1$ and
$J\ge J_2(\mathfrak K)$,
\[
 \E\,\mathbf D_J(\theta_0)\ge c_{\rm AD}\sum_{0<|k|\le\mathfrak K}|\widehat\theta_0(k)|^2 .
\]
Let $J\to\infty$ and then $\mathfrak K\to\infty$; recall $\widehat\theta_0(0)=0$.

\emph{Part (b): the probability bound.} All $u^{(J)}$ are defined on one probability space
(Section~\ref{sec:model}), and $X_J:=\mathbf D_J(\theta_0)/\|\theta_0\|^2\in[0,1]$
by \eqref{eq:energy}. Fatou's lemma applied to $1-X_J$ gives
\[
 \E\limsup_JX_J\ge\limsup_J\E X_J\ge c_{\rm AD},
\]
where the last inequality is (a). Since $\limsup_JX_J\le1$, this gives
\[
 \mathbb P(\limsup_JX_J>0)\ge c_{\rm AD}.
\]

\emph{Part (c): uniformity under a Sobolev bound.} Fix $\gamma>0$ and $A\ge1$, and put
$\mathfrak K:=\max\big(1,(2A^2)^{1/(2\gamma)}/(2\pi)\big)$.
For every real mean-zero $\theta_0\in H^\gamma(\T^d)$ with
$\|\theta_0\|_{\dot H^\gamma}\le A\|\theta_0\|$,
\[
 \sum_{|k|>\mathfrak K}|\widehat\theta_0(k)|^2
 \le(2\pi\mathfrak K)^{-2\gamma}\|\theta_0\|_{\dot H^\gamma}^2
 \le(2\pi\mathfrak K)^{-2\gamma}A^2\|\theta_0\|^2
 \le\tfrac12\|\theta_0\|^2.
\]
The finite-frequency bound in (a), with
$J_*:=J_2(\mathfrak K)$, therefore gives the claim, uniformly over
all such data. \qed

\subsection{Remarks}

\begin{remark}[Fixed cutoff]\label{rem:fixedJ}
For each fixed $J$, the field $u^{(J)}$ is spatially Lipschitz, so its
scalar dissipation vanishes as $\mu\downarrow0$: for
$\theta_0\in H^1(\T^d)$, the standard energy estimate for
$\nabla\theta^\mu$ (differentiate the equation, use
$\nabla\cdot u^{(J)}=0$ and \eqref{eq:holder}, and apply Gronwall's
inequality) gives
$\|\nabla\theta^\mu(t)\|\le e^{C_*\la_Jt}\|\nabla\theta_0\|$, hence
$2\mu\int_0^{t_0}\|\nabla\theta^\mu\|^2\dd t
\le2\mu t_0e^{2C_*\la_Jt_0}\|\nabla\theta_0\|^2\to0$ as
$\mu\downarrow0$; anomalous dissipation
arises here as the velocity cutoff is removed simultaneously with
the diffusivity.
\end{remark}

\begin{remark}
(1) For $\alpha=1/3$ the constant in Theorem~\ref{thm:ADgeneral} is
$c_{\rm AD}=c\kappa^3t_0^3$, with $c$ independent of $\kappa$ and
$t_0$; this reflects
the Richardson law \eqref{eq:momentbounds}.

(2) The probability bound in
Theorem~\ref{thm:ADgeneral}(b) is for each fixed datum, not simultaneous
over all data; we do not claim that the dissipation
event has probability one.

\end{remark}

\subsection{Anomalous dissipation for the fixed limiting field}
\label{sec:fixedfield}

We now compare the dissipation for $u^{(\infty)}$ with that of a
truncation whose matched diffusivity is comparable to $\mu$.
The velocity-tail estimate of Lemma~\ref{lem:tailpotential} supplies
the drift perturbation in the following comparison.

\subsubsection{Comparison under a perturbation of drift and diffusivity}

\begin{lemma}[Comparison of dissipations]\label{lem:disscomparison}
Let $u,v$ be drifts satisfying \eqref{eq:D} on $[0,t_0]\times\T^d$, and
suppose $u-v=\operatorname{div}B$ in the sense of distributions on
$\T^d$ for a.e.\ $t$, where $B\in L^\infty([0,t_0]\times\T^d)$ is
skew-symmetric. Write $b:=\|B\|_\infty$. For $\mu,\mu'>0$ let
$\theta,\phi$ be the energy solutions of \eqref{eq:adv} with drifts
$u,v$ and diffusivities $\mu,\mu'$, respectively, and common initial
value $\theta_0\in L^2$. Set
\[
 D_{u,\mu}:=2\mu\int_0^{t_0}\|\nabla\theta\|^2\dd t,
 \qquad D_{v,\mu'}:=2\mu'\int_0^{t_0}\|\nabla\phi\|^2\dd t.
\]
Then
\begin{equation}\label{eq:disscomparison}
 D_{u,\mu}\ge
 \frac{\mu'\mu}{(\mu'+|\mu-\mu'|+b)^2}\,D_{v,\mu'} .
\end{equation}
No smallness of $b/\mu$ is required. Exchanging the roles of $(u,\mu)$
and $(v,\mu')$, with $B$ replaced by $-B$, gives the reverse inequality;
in particular, for $\mu=\mu'$,
\begin{equation}\label{eq:twosided}
 (1+b/\mu)^{-2}D_{v,\mu}\le D_{u,\mu}\le(1+b/\mu)^{2}D_{v,\mu} .
\end{equation}
\end{lemma}

\begin{proof}
The skew-symmetry of $B$ lets us express the drift perturbation as
a bilinear form involving only first derivatives. We then estimate
the difference of the two scalar solutions in the energy norm.

\emph{Step 1: the skew-symmetry identity.} For $\theta,\varphi\in H^1(\T^d)$
and a.e.\ $t$,
\begin{equation}\label{eq:skewidentity}
 \big((\operatorname{div}B)\cdot\nabla\theta,\varphi\big)
 =\big(B\nabla\theta,\nabla\varphi\big).
\end{equation}
Indeed, for $\theta,\varphi\in C^\infty$, the definition of the
distributional divergence gives
\[
\begin{aligned}
 \big((\operatorname{div}B)\cdot\nabla\theta,\varphi\big)
 &=-\sum_{i,m}\big(B_{im},\partial_m(\varphi\,\partial_i\theta)\big)\\
 &=-\sum_{i,m}\big(B_{im}\partial_i\theta,\partial_m\varphi\big)
 -\sum_{i,m}\big(B_{im},\varphi\,\partial_{im}\theta\big).
\end{aligned}
\]
The last sum vanishes, because $B$ is skew-symmetric and
$\partial_{im}\theta$ is symmetric in $(i,m)$, and the first is
$-(B^{\mathsf T}\nabla\theta,\nabla\varphi)=(B\nabla\theta,\nabla\varphi)$.
Both sides of \eqref{eq:skewidentity} are continuous in
$(\theta,\varphi)\in H^1\times H^1$, because $\operatorname{div}B=u-v$ and
$B$ are bounded, and $C^\infty$ is dense in $H^1$.

\emph{Step 2: the equation for the difference.} Put $w:=\phi-\theta$.
Then $w\in L^2(0,t_0;H^1)$, $\partial_tw\in L^2(0,t_0;H^{-1})$ and
$w(0)=0$. Subtracting the identities \eqref{eq:weakform} for $\phi$
(drift $v$, diffusivity $\mu'$) and for $\theta$ (drift $u$,
diffusivity $\mu$), and using
$v\cdot\nabla\phi-u\cdot\nabla\theta=v\cdot\nabla w-(u-v)\cdot\nabla\theta$,
$u-v=\operatorname{div}B$ and \eqref{eq:skewidentity}, we get, for
a.e.\ $t$ and all $\varphi\in H^1$,
\begin{equation}\label{eq:diffeq}
 \langle\partial_tw,\varphi\rangle+(v\cdot\nabla w,\varphi)
 +\mu'(\nabla w,\nabla\varphi)
 =\big(((\mu-\mu')I+B)\nabla\theta,\nabla\varphi\big).
\end{equation}

\emph{Step 3: energy estimate.} Take $\varphi=w(t)$ in
\eqref{eq:diffeq}. By (E1) and \eqref{eq:antisym}, and since
$|B\nabla\theta|\le b|\nabla\theta|$ pointwise, for a.e.\ $t$
\[
 \tfrac12\tfrac{d}{dt}\|w\|^2+\mu'\|\nabla w\|^2
 \le(|\mu-\mu'|+b)\,\|\nabla\theta\|\,\|\nabla w\| .
\]
Integrating over $[0,t_0]$, using $w(0)=0$ and the Cauchy--Schwarz
inequality in $t$,
\[
 \tfrac12\|w(t_0)\|^2+\mu'\|\nabla w\|_{L^2_{t,x}}^2
 \le(|\mu-\mu'|+b)\,\|\nabla\theta\|_{L^2_{t,x}}\,\|\nabla w\|_{L^2_{t,x}},
\]
so $\|\nabla w\|_{L^2_{t,x}}\le\frac{|\mu-\mu'|+b}{\mu'}\|\nabla\theta\|_{L^2_{t,x}}$.
Since $\phi=\theta+w$,
\[
 \|\nabla\phi\|_{L^2_{t,x}}
 \le\|\nabla\theta\|_{L^2_{t,x}}+\|\nabla w\|_{L^2_{t,x}}
 \le\frac{\mu'+|\mu-\mu'|+b}{\mu'}\,\|\nabla\theta\|_{L^2_{t,x}} .
\]
Squaring and multiplying by $2\mu'$ gives
$D_{v,\mu'}\le\frac{(\mu'+|\mu-\mu'|+b)^2}{\mu'\mu}D_{u,\mu}$, which is
\eqref{eq:disscomparison}. The hypotheses are symmetric in
$(u,\mu,B)$ and $(v,\mu',-B)$, which gives the last assertion.
\end{proof}

\begin{remark}
Lemma~\ref{lem:disscomparison} is a time-dependent form of the
classical bound $\mu_{\rm eff}\le\mu+\|B\|_\infty^2/\mu$ for the
effective diffusivity of a periodic flow $u=\operatorname{div}B$ in
terms of its stream matrix $B$, compare \cite{FP94}: a perturbation
of the drift by $\operatorname{div}B$ with $\|B\|_\infty\lesssim\mu$,
that is, at bounded P\'eclet number, changes the dissipation by at
most a constant factor, in either direction.
\end{remark}

\subsubsection{Proof of Theorem~\ref{thm:ADfixed}}

\begin{proof}[Proof of Theorem~\ref{thm:ADfixed}]
We first justify measurability for the limiting field and then apply
the comparison lemma at a truncation adapted to $\mu$. Put
\begin{equation}\label{eq:astar}
 a_*:=\big(2^{(1+\alpha)/2}+C_{\rm pot}/\varepsilon_0\big)^{-2}\in(0,1).
\end{equation}
Since $u^{(N)}\to u^{(\infty)}$ uniformly and each $u^{(N)}$ is
bounded, Borel measurable, and divergence-free, the limiting field
satisfies \eqref{eq:D} on $[0,1]\times\T^d$. Hence
Section~\ref{sec:weak} gives existence, uniqueness and the energy
identity for \eqref{eq:scalarfixed}, and
$\mathbf D^\mu_\infty(\theta_0)=\|\theta_0\|^2-\|\theta^\mu(t_0)\|^2$.

\emph{Measurability with respect to the random data defining the velocity.}
To take the expectations in (a) and (c), we first show that, for each
fixed $\mu>0$, the dissipation $\mathbf D^\mu_\infty(\theta_0)$ is a
measurable function of the random data. For the probability statement
in (b), we also show continuity in $\mu>0$ for each realization; this
allows the $\limsup$ to be taken over rational diffusivities and
ensures that the event in (b) is measurable.

Fix $\mu>0$ and let $\theta^{\mu,N}$ be the energy
solution with drift $u^{(N)}$ and diffusivity $\mu$. Subtracting the
identities \eqref{eq:weakform} for $\theta^{\mu,N}$ and $\theta^\mu$ and
using \eqref{eq:antisym}, $w:=\theta^{\mu,N}-\theta^\mu$ satisfies, for
a.e.\ $t$ and all $\varphi\in H^1$,
\[
\begin{aligned}
 &\langle\partial_tw,\varphi\rangle+(u^{(N)}\cdot\nabla w,\varphi)
 +\mu(\nabla w,\nabla\varphi)\\
 &\qquad=\big((u^{(\infty)}-u^{(N)})\cdot\nabla\theta^\mu,\varphi\big)\\
 &\qquad=-\big((u^{(\infty)}-u^{(N)})\theta^\mu,\nabla\varphi\big).
\end{aligned}
\]
With $\varphi=w(t)$, (E1), Young's inequality and
$\|\theta^\mu(t)\|\le\|\theta_0\|$,
\[
 \sup_{t\le t_0}\|\theta^{\mu,N}(t)-\theta^\mu(t)\|^2
 \le\frac{t_0}{\mu}\,\|u^{(\infty)}-u^{(N)}\|_\infty^2\,\|\theta_0\|^2
 \longrightarrow0\qquad(N\to\infty).
\]
Hence $\mathbf D^\mu_\infty(\theta_0)$ is the pointwise limit of
$\|\theta_0\|^2-\|\theta^{\mu,N}(t_0)\|^2$, each of which is a Borel
function of the data by Lemma~\ref{lem:datadependence} (applied with
$J=N$ and diffusivity $\mu$); so it is Borel. Next, for a fixed
realization and $\mu,\mu'>0$, the same computation for
$w:=\theta^{\mu'}-\theta^\mu$, whose equation has drift $u^{(\infty)}$,
diffusivity $\mu'$ and right-hand side
$(\mu-\mu')(\nabla\theta^\mu,\nabla\varphi)$, gives
\[
 \sup_{t\le t_0}\|\theta^{\mu'}(t)-\theta^\mu(t)\|^2
 \le\frac{|\mu'-\mu|^2}{\mu'}\int_0^{t_0}\|\nabla\theta^\mu\|^2\dd t
 \le\frac{|\mu'-\mu|^2}{2\mu\mu'}\|\theta_0\|^2 .
\]
Thus $\mu\mapsto\mathbf D^\mu_\infty(\theta_0)$ is continuous on
$(0,\infty)$, the $\limsup$ in (b) may be taken over rational $\mu$,
and the event in (b) is measurable.

\emph{Comparison with a truncation.} Let $\mu_J=\varepsilon_0\ell_J^{1+\alpha}$ and
$q:=2^{1+\alpha}=\mu_{J}/\mu_{J+1}$. For $0<\mu<q\mu_{j_0}$ let
$J=J(\mu)\ge j_0$ be the unique index with $\mu_J\le\mu<q\mu_J$; then
$J(\mu)\to\infty$ as $\mu\downarrow0$. Lemmas~\ref{lem:tailpotential}
and~\ref{lem:disscomparison}, applied for each realization with
$(u,\mu)=(u^{(\infty)},\mu)$ and $(v,\mu')=(u^{(J)},\mu_J)$, so
that $D_{u,\mu}=\mathbf D^\mu_\infty(\theta_0)$ and
$D_{v,\mu'}=\mathbf D_J(\theta_0)$, with $B=B_{>J}$ (so that
$u-v=\operatorname{div}B_{>J}$) and
$b\le C_{\rm pot}\ell_J^{1+\alpha}=(C_{\rm pot}/\varepsilon_0)\mu_J$, give
\[
 \mathbf D_\infty^\mu(\theta_0)
 \ge\frac{\mu_J\mu}{(\mu+C_{\rm pot}\ell_J^{1+\alpha})^2}\,
       \mathbf D_J(\theta_0)
 \ge a_*\,\mathbf D_{J(\mu)}(\theta_0) .
\]
The first inequality uses $\mu\ge\mu'=\mu_J$, so that
$\mu'+|\mu-\mu'|=\mu$ in \eqref{eq:disscomparison}.
For the last inequality write $c':=\mu/\mu_J\in[1,q)$: the middle
factor is 
\[
\begin{aligned}
 c'/(c'+C_{\rm pot}/\varepsilon_0)^2
 &=\big(\sqrt{c'}+C_{\rm pot}/(\varepsilon_0\sqrt{c'})\big)^{-2}\\
 &\ge(\sqrt q+C_{\rm pot}/\varepsilon_0)^{-2}=a_*.
\end{aligned}
\]

Taking expectations and letting $\mu\downarrow0$, so that
$J(\mu)\to\infty$,
Theorem~\ref{thm:ADgeneral}(a) gives (a). For (c), fix
$\gamma>0$ and $A\ge1$, let
$J_*=J_*(\gamma,A,\kappa,\varepsilon_0,t_0)$ be as in
Theorem~\ref{thm:ADgeneral}(c), uniformly for real mean-zero
$\theta_0\in H^\gamma$ with
$\|\theta_0\|_{\dot H^\gamma}\le A\|\theta_0\|$, and set
$\mu_*:=\mu_{J_*}$; then
$\mu\le\mu_*$ implies $J(\mu)\ge J_*$, and (c) follows.

\emph{The probability statement.} For (b), apply \eqref{eq:twosided} at the matched diffusivities
$\mu=\mu'=\mu_J$:
\begin{equation}\label{eq:fixedmatched}
 (1+C_{\rm pot}/\varepsilon_0)^{-2}\mathbf D_J(\theta_0)
 \le\mathbf D_\infty^{\mu_J}(\theta_0)
 \le(1+C_{\rm pot}/\varepsilon_0)^{2}\mathbf D_J(\theta_0),
\end{equation}
valid for every realization and every $J\ge j_0$.
Hence the event $\{\limsup_J\mathbf D_J(\theta_0)>0\}$ is contained in
$\{\limsup_{\mu\downarrow0}\mathbf D^\mu_\infty(\theta_0)>0\}$, and
Theorem~\ref{thm:ADgeneral}(b) gives the probability bound.
\end{proof}

Part (b) says that, on an event of probability at least $c_{\rm AD}$,
a single realization of $u^{(\infty)}$ dissipates the variance of
$\theta_0$ anomalously. For the explicit model of
Proposition~\ref{prop:explicitwavelet} this gives a deterministic
statement.

\begin{corollary}\label{cor:deterministic}
Let $d\ge2$, $0<\alpha<1$ and $t_0\in(0,1]$. There is a divergence-free
vector field $u\in L^\infty([0,1];C^\alpha(\T^d))$ such that, for
$\theta_0(x)=\cos(2\pi x_1)$, the energy solutions $\theta^\mu$ of
$\partial_t\theta^\mu+u\cdot\nabla\theta^\mu=\mu\Delta\theta^\mu$,
$\theta^\mu(0)=\theta_0$, satisfy
\[
 \limsup_{\mu\downarrow0}\ \mu\int_0^{t_0}\|\nabla\theta^\mu(t)\|^2\dd t>0 .
\]
\end{corollary}

\begin{proof}
Take $\psi(x)=(1-|x|^2)_+^5$ and independent coefficients $\pm1$ with
probability $\frac12$. This model satisfies \eqref{eq:ND}
(Proposition~\ref{prop:explicitwavelet}), \eqref{eq:coefficients} with
$M=\sigma=1$, and (S). Fix $\kappa=\kappa_1$ and
$\varepsilon_0=\varepsilon_1(\kappa_1)$. By Theorem~\ref{thm:ADfixed}(b)
with this $t_0$ and $\theta_0=\cos(2\pi x_1)$, the event
$\{\limsup_{\mu\downarrow0}\mathbf D^\mu_\infty(\theta_0)>0\}$ has
probability at least $c\,(\kappa_1t_0)^{2/\beta}>0$; let $u$ be
the field $u^{(\infty)}$ of any realization in this event.
\end{proof}

\begin{remark}
(1) By \eqref{eq:fixedmatched}, at the diffusivity $\mu_J$ the
dissipation for the fixed field and the dissipation for the truncation
$u^{(J)}$ are comparable, with the constant $(1+C_{\rm pot}/\varepsilon_0)^2$,
for every realization: the scales below the diffusive scale change the
dissipation by at most a constant factor. The same comparison, with
the constant $(2^{2+\alpha}+C_{\rm pot}/\varepsilon_0)^2$, holds between
$\mathbf D^\mu_\infty(\theta_0)$ and $\mathbf D_{J(\mu)}(\theta_0)$ for
every $\mu<2^{1+\alpha}\mu_{j_0}$, by Lemma~\ref{lem:disscomparison}
and $\mu_{J(\mu)}\le\mu<2^{1+\alpha}\mu_{J(\mu)}$.

(2) Part (a) is a $\liminf$ of expectations over all diffusivities
tending to zero. In part (b), the $\limsup$ is taken for each fixed
realization of the velocity; it is positive with probability at least
$c_{\rm AD}$ for each fixed nonzero datum. We do not claim a positive
$\liminf$ for individual realizations, probability one, or a single
realization working simultaneously for all data: from $\E X_J\ge c$
and $X_J\in[0,1]$ one cannot conclude anything about $\liminf_JX_J$
without information on the joint law of the $X_J$.
\end{remark}

\appendix

\section{Proofs of the lemmas of Section~\ref{sec:weak}}\label{sec:appendix}

The notation, the facts (E1) and (E2), and the identities
\eqref{eq:antisym}--\eqref{eq:weakform} are those of
Section~\ref{sec:weak}.

\begin{proof}[Proof of Lemma~\ref{lem:weakbasic}]
(a) By (E1) with $\widetilde w=w$, \eqref{eq:weakform} with $\varphi=w(s)$
and \eqref{eq:antisym},
$\frac{d}{ds}\|w\|^2=2\langle\partial_sw,w\rangle=-2\mu\|\nabla w\|^2$ for
a.e.\ $s$. The difference of two energy solutions with the same
initial value is an energy solution with initial value $0$, hence
vanishes.

(b) The constant $1$ is in $L^2(I;H^1)$ with derivative $0$, so by (E1)
and \eqref{eq:weakform} with $\varphi=1$,
$\frac{d}{ds}\int w=\langle\partial_sw,1\rangle=-(\mathbf b\cdot\nabla w,1)=0$
by \eqref{eq:antisym}. For positivity write $v_-=\max(-v,0)$ and let
$w_k$ be as in (E1). The functional $v\mapsto\frac12\|v_-\|^2$ is
continuously differentiable on $L^2$ with derivative $-v_-$, since
$r\mapsto\frac12(r_-)^2$ has the $1$-Lipschitz derivative $-r_-$. As
$s\mapsto w_k(s)\in L^2$ is $C^1$,
\[
 \tfrac12\|(w_k)_-(s)\|^2-\tfrac12\|(w_k)_-(T_0)\|^2
 =-\int_{T_0}^s\langle\partial_sw_k,(w_k)_-\rangle .
\]
As $k\to\infty$ the left side converges to the same expression for $w$,
because $|a_--b_-|\le|a-b|$ and $w_k\to w$ in $C(I;L^2)$. On the right,
$\partial_sw_k\to\partial_sw$ in $L^2(I;H^{-1})$, while $(w_k)_-\to w_-$ in
$L^2(I;L^2)$ and $\|\nabla(w_k)_-\|\le\|\nabla w_k\|$, so $(w_k)_-$ is
bounded in $L^2(I;H^1)$ and converges to $w_-$ weakly there. Hence
$\frac12\|w_-(s)\|^2-\frac12\|w_-(T_0)\|^2=-\int_{T_0}^s\langle\partial_sw,w_-\rangle$.
Here we used that $v_-\in H^1$ with $\nabla v_-=-\one_{\{v<0\}}\nabla v$ for
$v\in H^1$ \cite[Lemma~7.6]{GT01}. By \eqref{eq:weakform}, this formula,
the pointwise identity $(\mathbf b\cdot\nabla w)w_-=-(\mathbf b\cdot\nabla w_-)w_-$
and \eqref{eq:antisym},
\[
 -\langle\partial_sw,w_-\rangle=(\mathbf b\cdot\nabla w,w_-)+\mu(\nabla w,\nabla w_-)
 =-(\mathbf b\cdot\nabla w_-,w_-)-\mu\|\nabla w_-\|^2\le0 .
\]
Thus $\|w_-(s)\|\le\|w_-(T_0)\|=0$.

(c) The function $G$ belongs to $C([s_1,s_2];H^1)$, and its classical
derivative $\partial_sG$, which is continuous, is its distributional
derivative. By (E1) on $[s_1,s_2]$, $s\mapsto\int Gw$ is absolutely
continuous with derivative $\langle\partial_sw,G\rangle+(\partial_sG,w)$. By
\eqref{eq:weakform}, \eqref{eq:antisym} and an integration by parts,
$\langle\partial_sw,G\rangle=-(\mathbf b\cdot\nabla w,G)-\mu(\nabla w,\nabla G)
=(\mathbf b\cdot\nabla G,w)+\mu(w,\Delta G)$. Integrate in $s$.
\end{proof}

\begin{proof}[Proof of Lemma~\ref{lem:transfer}]
Let $w^{(1)},w^{(2)}$ be the energy solutions with initial values
$\one_Aw(s_1)$ and $\one_{A^c}w(s_1)$. By Lemma~\ref{lem:weakbasic}(a),
$w=w^{(1)}+w^{(2)}$. By Lemma~\ref{lem:weakbasic}(b), $w^{(1)}\ge0$ and
$\int w^{(1)}(s_2)=\int_Aw(s_1)$, so $\int_Bw^{(1)}(s_2)\le\int_Aw(s_1)$.
If $m:=\int_{A^c}w(s_1)>0$, then $w^{(2)}/m$ is admissible in
\eqref{eq:Pi} for the set $A^c$, so $\int_Bw^{(2)}(s_2)\le m\,\Pi_{s_1,s_2}(A^c;B)$;
if $m=0$, then $w^{(2)}=0$.
\end{proof}

\begin{proof}[Proof of Lemma~\ref{lem:duality}]
The function $\widetilde w(s):=w(T-s)$ lies in $L^2(0,T;H^1)$, with
$\partial_s\widetilde w(s)=-(\partial_sw)(T-s)\in L^2(0,T;H^{-1})$. By (E1),
$I(s):=(\widetilde w(s),\varrho(s))$ is absolutely continuous, and for a.e.\
$s$, by \eqref{eq:weakform} for $w$ at time $T-s$ and for $\varrho$ at
time $s$, with $\mathbf b=\mathbf b(T-s)$,
\[
 I'(s)=\big[(\mathbf b\cdot\nabla\widetilde w,\varrho)+\mu(\nabla\widetilde w,\nabla\varrho)\big]
 +\big[(\mathbf b\cdot\nabla\varrho,\widetilde w)-\mu(\nabla\varrho,\nabla\widetilde w)\big]=0,
\]
by \eqref{eq:antisym}. Hence $I(0)=I(T)$.
\end{proof}

\begin{proof}[Proof of Lemma~\ref{lem:stability}]
By Lemma~\ref{lem:weakbasic}(a), $\|w_k(s)\|\le\|w_0\|$ and
$2\mu\int_I\|\nabla w_k\|^2\le\|w_0\|^2$. It suffices to show that every
subsequence has a further subsequence along which $w_k(s)\rightharpoonup w(s)$
for every $s$. Pass to a subsequence along which
$w_k\rightharpoonup\bar w$ weakly in $L^2(I;H^1)$. For
$\varphi\in C^\infty(\T^n)$, Lemma~\ref{lem:weakbasic}(c) with
$G(s,Y)=\varphi(Y)$ gives
\[
 (w_k(s),\varphi)=(w_0,\varphi)+\int_{T_0}^s(w_k,\mathbf b_k\cdot\nabla\varphi+\mu\Delta\varphi)\dd r .
\]
Since $\mathbf b_k\cdot\nabla\varphi\to\mathbf b\cdot\nabla\varphi$ in
$L^2(I\times\T^n)$ by dominated convergence, the right side converges for
every $s$ to
$\ell_s(\varphi):=(w_0,\varphi)+\int_{T_0}^s(\bar w,\mathbf b\cdot\nabla\varphi+\mu\Delta\varphi)\dd r$.
As $\|w_k(s)\|\le\|w_0\|$, it follows that for every $s$, $w_k(s)$
converges weakly in $L^2$ to some $\lambda_s$ with $(\lambda_s,\varphi)=\ell_s(\varphi)$.
Testing against $\chi\in L^2(I)$ and using dominated convergence shows
$(\bar w(s),\varphi)=\ell_s(\varphi)$ for a.e.\ $s$. By \eqref{eq:antisym},
$(\bar w(r),\mathbf b\cdot\nabla\varphi+\mu\Delta\varphi)=\langle a_r,\varphi\rangle$, where
$a_r:=-\mathbf b(r)\cdot\nabla\bar w(r)+\mu\Delta\bar w(r)\in H^{-1}$ has norm at
most $(B_0+\mu)\|\nabla\bar w(r)\|$. The map $r\mapsto a_r$ is weakly
measurable, hence strongly measurable by Pettis' theorem, and lies in
$L^2(I;H^{-1})$. Thus $\bar w(s)=w_0+\int_{T_0}^sa_r\dd r$ in $H^{-1}$ for
a.e.\ $s$, so $\partial_s\bar w=a\in L^2(I;H^{-1})$, and \eqref{eq:weakform}
holds: $\bar w$ is an energy solution. Its continuous representative satisfies
$(\bar w(s),\varphi)=\ell_s(\varphi)$ for all $s$, in particular
$\bar w(T_0)=w_0$. By uniqueness $\bar w=w$, so $\lambda_s=w(s)$ for every $s$.
\end{proof}

\begin{proof}[Proof of Lemma~\ref{lem:products}]
(a) Let $N$ be the null set of \eqref{eq:D} for $b$. For $s\notin N$ and
$\varphi\in C^\infty(\T^{2d})$,
$\int\!\!\int b(s,y)\cdot\nabla_y\varphi\dd y\dd y'=0$, since the inner integral
vanishes for each $y'$, and similarly in $y'$; so $\mathbf b$ satisfies
\eqref{eq:D}. Since
$\sup_s\|\theta_i(s)\|<\infty$ and $\theta_i\in L^2(I;H^1)$,
$\theta_1\otimes\theta_2\in L^2(I;H^1(\T^{2d}))$
and $\theta_1\otimes\theta_2\in C(I;L^2(\T^{2d}))$. For $\varphi_1,\varphi_2\in C^\infty(\T^d)$, by (E1)
and \eqref{eq:weakform} for $\theta_1,\theta_2$, the function
$((\theta_1\otimes\theta_2)(s),\varphi_1\otimes\varphi_2)=(\theta_1(s),\varphi_1)(\theta_2(s),\varphi_2)$ is absolutely
continuous with derivative $-a_s(\varphi_1\otimes\varphi_2)$, where, for
$\varphi\in H^1(\T^{2d})$,
\[
\begin{gathered}
 a_s(\varphi):=(\mathbf b(s)\cdot\nabla(\theta_1\otimes\theta_2)(s),\varphi)
 +\mu(\nabla(\theta_1\otimes\theta_2)(s),\nabla\varphi),\\
 |a_s(\varphi)|\le(\sqrt2B_0+\mu)\|\nabla(\theta_1\otimes\theta_2)(s)\|\,\|\varphi\|_{H^1}.
\end{gathered}
\]
Finite sums of such products are dense in $H^1(\T^{2d})$ (trigonometric
polynomials), so $\partial_s(\theta_1\otimes\theta_2)=-a_s\in L^2(I;H^{-1}(\T^{2d}))$ and
\eqref{eq:weakform} holds. The claim follows by uniqueness; the last
statement is obvious.

(b) The map $\iota$ preserves Lebesgue measure,
$\nabla(v\circ\iota^{-1})=g(\nabla v)\circ\iota^{-1}$, and $g$ is orthogonal. Hence
$(\iota_*\mathbf b\cdot\nabla(w\circ\iota^{-1}),\varphi)=(\mathbf b\cdot\nabla w,\varphi\circ\iota)$,
$(\nabla(w\circ\iota^{-1}),\nabla\varphi)=(\nabla w,\nabla(\varphi\circ\iota))$ and
$\langle\partial_s(w\circ\iota^{-1}),\varphi\rangle=\langle\partial_sw,\varphi\circ\iota\rangle$, so
\eqref{eq:weakform} for $w\circ\iota^{-1}$ and $\varphi$ is \eqref{eq:weakform}
for $w$ and $\varphi\circ\iota$. The same change of variables gives
\eqref{eq:D} for $\iota_*\mathbf b$.
\end{proof}

\fussy

\end{document}